\documentclass{article}

\usepackage[numbers,sort&compress]{natbib}

\usepackage{mathtools}
\usepackage{authblk}
\usepackage[title]{appendix}
\usepackage{geometry}
\usepackage{subfiles}
\usepackage{caption}
\usepackage{subcaption}
\usepackage{setspace}
\usepackage{amsthm}
\usepackage{amsmath}
\usepackage{amssymb}
\usepackage{ulem}
\usepackage[dvipsnames]{xcolor}

\usepackage{mathrsfs}

\mathtoolsset{showonlyrefs}

\usepackage{float}
\usepackage{fancyhdr}
\usepackage{graphicx}
\usepackage[colorlinks,            linkcolor=black,            anchorcolor=black,            citecolor=black            ]{hyperref}
\usepackage{listings}
 \usepackage{tikz}
 \usetikzlibrary{arrows,decorations.pathmorphing,backgrounds,positioning,fit,petri,arrows.meta,bending,positioning}

\graphicspath{ {./images/} }
\newtheorem{corollary}{Corollary}[section]

\newtheorem{proposition}{Proposition}[section]
\newtheorem{theorem}{Theorem}[section]
\newtheorem{lemma}{Lemma}[section]
\newtheorem{definition}{Definition}[section]

\newtheorem{remark}{Remark}[section]
\newtheorem{assumption}{Assumption}
\newtheorem*{theorem*}{Theorem}
\numberwithin{equation}{section}

\newcommand{\p}{\partial}
\newcommand{\beq}{\begin{equation}}
\newcommand{\eeq}{\end{equation}}

\newcommand{\tp}{\tilde{p}}

\newcommand{\bueq}{\begin{equation*}}
\newcommand{\eueq}{\end{equation*}}

\newcommand{\bthm}{\begin{theorem}}
	\newcommand{\ethm}{\end{theorem}}

\newcommand{\eps}{\epsilon}

\providecommand{\keywords}[1]
{
	\small	
	\textbf{{Keywords:}} #1
}
\providecommand{\msc}[1]
{
	\small	
	\textbf{{Mathematics Subject Classification:}} #1
}

\title{ Dilating blow-up time: A generalized solution of the NNLIF neuron model and its global well-posedness}%
\author{Xu'an Dou\thanks{Beijing International Center for Mathematical Research, Peking University, Beijing, 100871, China (dxa@pku.edu.cn)} \,\, and \,\, Zhennan Zhou\thanks{Beijing International Center for Mathematical Research, Peking University, Beijing, 100871, China (zhennan@bicmr.pku.edu.cn).}}

\date{\today}
\begin{document}
\maketitle
\begin{abstract}
The nonlinear noisy leaky integrate-and-fire (NNLIF) model is a popular mean-field description of a large number of interacting neurons, which has attracted mathematicians to study from various aspects. A core property of this model is the finite time blow-up of the firing rate, which scientifically corresponds to the synchronization of a neuron network, and mathematically prevents the existence of a global classical solution. In this work, we propose a new generalized solution based on reformulating the PDE model with a specific change of variable in time. A firing rate dependent timescale is introduced, in which the transformed equation can be shown to be globally well-posed for any connectivity parameter {even in the event of the blow-up.} The generalized solution is then defined via the backward change of timescale, and it may have a jump when the firing rate blows up. We establish properties of the generalized solution including the characterization of blow-ups and the global well-posedness in the original timescale. The generalized solution provides a new perspective to understand the dynamics when the firing rate blows up as well as the continuation of the solution after a blow-up.%
\end{abstract}


\msc{35Q84; 35Q92; 35B44; 35Dxx; 92B20}


\section{Introduction}

To describe a large number of interacting neurons, mean-field approaches have been successfully used in computational neuroscience \cite{brunel1999fast,omurtag2000simulation,renart2004mean,cai2004effective}, which lead to various novel models including nonlinear partial differential equations with unique structures. Recently such equations have attracted many mathematicians to study both from the PDE point of view \cite{caceres2011analysis,perthame2013voltage,mischler2016kinetic,dumont2020mean} and from the corresponding SDE point of view \cite{delarue2015global,delarue2015particle,cormier2020hopf,mna:8775}. One of the common features shared by these models is the synchronization, which means a positive fraction of neurons spike at the same time. Mathematically this leads to the finite-time blow-up of the firing rate which obstructs the existence of a global classical solution. Our goal of this work, is to propose a new generalized solution which admits blow-ups, and can be shown to be globally well-posed. 

\paragraph{1.1 Model description}
We consider the nonlinear noisy leaky integrate-and-fire (NNLIF for short) model proposed in \cite{amit1997model,brunel1999fast,brunel2000dynamics}, which has become a standard model in computational neuroscience (see e.g. \cite{burkitt2006review,burkitt2006review2,gerstner2014neuronal}). Here each neuron is characterized by its voltage (membrane potential) $v$, which is no greater than a threshold $V_F$. And the ensemble of neurons is described by $p(v,t)$, the probability density at each time $t$ to find a neuron with voltage $v$. The density function $p(v,t)$ is governed by the following PDE 
\begin{equation}\label{eq:NNLIF-intro}
    \p_tp+\p_v\bigl((-v+V_L+\mu(t))p\bigr)=a(t)\p_{vv}p+N(t)\delta_{v=V_R},\quad t>0,v\in(-\infty,V_F).
\end{equation}The term $-v+V_L$ in the drift models the leaky mechanism which drives the voltage to the equilibrium potential $V_L<V_F$. Each neuron is also subject to a fluctuating input current, of which the mean is $\mu (t)$ and the noise level is described by $a(t)$. Such an input can drive the voltage of a neuron away from $V_L$ and may exceed $V_F$ at a certain time.

The unique feature of \eqref{eq:NNLIF-intro} comes from the firing mechanism. A neuron spikes when its voltage arrives at $V_F$, the firing potential. In this macroscopic model, the spike activity is described by the mean firing rate $N(t)$, whose physical meaning is the number of spikes in expectation per unit time.

A spike has two consequences. First, after the spike, the voltage of the neuron is immediately reset to $V_R$, the reset potential. Such a reset gives the rise to the Dirac source term $N(t)\delta_{v=V_R}$ in \eqref{eq:NNLIF-intro}, where we use $\delta_{v=V_R}$ to denote the Dirac measure supported at $V_R$. Physically, we have $V_L\leq V_R<V_F$. Also the instantaneous reset gives an absorbing boundary condition at $V_F$
\begin{equation}\label{intro-bc}
    p(V_F,t)=0,\quad t>0,
\end{equation} since no neuron is allowed to stay at $V_F$.

Second, when a neuron spikes, it influences other neurons by contributing to their input currents. In this way, the input current consists of two parts: an external part and an intrinsic part from the spikes of other neurons within this network. Precisely, its mean $\mu(t)$ and noise level $a(t)$ are given by
\begin{equation}\label{intro-mu-a}
    \mu(t)=\mu_0+bN(t),\quad a(t)=a_0+a_1N(t).
\end{equation} Here $\mu_0\in\mathbb{R}$ and $a_0>0$ represents the mean and the noise level of the external part, which are assumed to be constant in time. The connectivity parameter $b\in\mathbb{R}$ determines the influence of one spike. When $b>0$, the spike of a neuron increases the voltage of other neurons, which is called the excitatory case. When $b<0$, the spike of a neuron decreases the voltage of other neurons, which is called the inhibitory case. Last but not least, $a_1\geq 0$ is the scaling parameter of the noise level induced by spikes. 

Finally, the firing rate $N(t)$ is determined in a self-consistent way 
\begin{equation}\label{intro-bd-flux}
    N(t):=-a(t)\p_vp(V_F,t)=-(a_0+a_1N(t))\p_vp(V_F,t).
\end{equation} Note the boundary condition $p(V_F,t)=0$ \eqref{intro-bc}. \eqref{intro-bd-flux} means that the firing rate is equal to the boundary flux crossing $V_F$, which gives the conservation of mass in \eqref{eq:NNLIF-intro}
\begin{equation*}
    \frac{d}{dt}\left(\int_{-\infty}^{V_F}p(v,t)dv\right)=0.
\end{equation*} In view of \eqref{intro-bd-flux}, the dependence of $N(t)$ in \eqref{intro-mu-a} makes \eqref{eq:NNLIF-intro} nonlinear. The system is completed with the boundary condition at the infinity $p(-\infty,t)=0,t>0,$ and an initial condition $p(v,0)=p_I(v),\, v\leq V_F$. {For more detailed derivations of this model, we refer to \cite{brunel1999fast,caceres2011analysis,liu2020rigorous}}.

The PDE \eqref{eq:NNLIF-intro} can be viewed as a Fokker-Planck equation of a corresponding Mckean-Vlasov SDE \cite{delarue2015global,delarue2015particle}. Similar mathematical structures have also been investigated in a wide range of modeling context, including portfolio or bank systems in mathematical finance \cite{hambly2019mckean,2019-AAP1403}  and the supercooled Stefan problem \cite{delarue2022global}.

\paragraph{1.2 Literature Review}
Mathematical analysis on \eqref{eq:NNLIF-intro} started relatively recently albeit its modeling success. A remarkable phenomenon of \eqref{eq:NNLIF-intro} is the possible finite-time blow-up of the firing rate $N(t)$, first discovered in \cite{caceres2011analysis}. Such a blow-up not only becomes a core concept in the mathematical studies of \eqref{eq:NNLIF-intro} but is indeed scientifically relevant. 

Mathematically, the blow-up of $N(t)$ is exactly the obstacle for the existence of a global solution in the classical theory, while \cite{carrillo2013CPDEclassical} has established the global well-posedness of a classical solution as long as $N(t)$ remains finite. Many studies have been devoted to the emergence of the blow-up and the existence of a global solution in the absence of the blow-up. In \cite{carrillo2013CPDEclassical} (and via a different method in \cite{Antonio_Carrillo_2015}), the blow-up of $N(t)$ is excluded when $b\leq 0,a_1=0$, and thus in this case the global well-posedness is established.  However, as long as $b>0$, there exists smooth initial data such that the blow-up occurs in finite time \cite{caceres2011analysis,roux2021towards}. Conditions on the initial data and the connectivity parameter $b>0$ to exclude or ensure the blow-up have been studied from a stochastic perspective in \cite{delarue2015global,hambly2019mckean} and via deterministic methods in \cite{roux2021towards}. Besides, if a time-delayed effect is added to \eqref{eq:NNLIF-intro}, then the resulting modified model has been shown to be globally well-posed \cite{delarue2015particle,caceres2019global}.

Scientifically, a blow-up of the firing rate $N(t)$ corresponds to synchronization of a neuron network. {Here synchronization refers to the scenario when a portion of the neurons within the network fire simultaneously}, {a phenomenon ubiquitous in a brain and of fundamental importance in neuroscience}. We remark that this definition of synchronization shares many similarities to the so-called multi-firing event, {a relatively recently proposed and less understood concept in computational neuroscience \cite{rangan2013dynamics,rangan2013emergent,zhang2014coarse}}.

Quite interesting it is to consider {{generalized solutions}} of \eqref{eq:NNLIF-intro}, which allows the blow-up of $N(t)$ to be part of the dynamics. Such a generalized solution might yield global well-posedness for a much larger parameter regime than the classical counterpart, and it can describe the dynamics in the presence of synchronization. A pioneering work in this direction is \cite{delarue2015particle} which defines a generalized solution called the ``physical solution'' for the corresponding SDE of \eqref{eq:NNLIF-intro}. In their definition, the blow-up of $N(t)$ gives a Dirac measure in time, which induces a jump of the solution. The jump size, which corresponds to the proportion of neurons synchronized at that time, is characterized by a ``physical'' condition. They proved the global existence of this ``physical solution'' but its uniqueness was left open. Such a physical solution has been studied subsequently in \cite{2019-AAP1403,hambly2019mckean,2020BLMSLedgerAndreas}. More recently, a detailed characterization for the regularity of $N(t)$ around the blow-up time is given and the global uniqueness is finally established in \cite{delarue2022global}. Actually, \cite{2019-AAP1403,hambly2019mckean,2020BLMSLedgerAndreas,delarue2022global} study different but mathematically closely related models rather than \eqref{eq:NNLIF-intro} itself. For those models blow-up is also important. Up to now, most theoretical contributions to generalized solutions of \eqref{eq:NNLIF-intro} works with the SDE formulation with probability techniques, see also a recent numerical study \cite{CiCP-30-820}.

Last but not least, all global existence results aforementioned deals with the case $a_1=0$. By \eqref{intro-bd-flux} the firing rate is given by
\begin{equation}\label{intro-Nt}
N(t)=\frac{-a_0\p_vp(V_F,t)}{1+a_1\p_vp(V_F,t)}.
\end{equation}When $a_1=0$, $N(t)<+\infty$ is equivalent to $-\p_vp(V_F,t)<+\infty$. But when $a_1>0$, to have a finite and non-negative $N(t)$ we need to require additionally $-\p_vp(V_F,t)<1/a_1$. Thus in terms of the classical solution, the case $a_1>0$ is more ``ill-posed''. Indeed, \cite[Appendix A]{Antonio_Carrillo_2015} has shown that when $a_1>0$ for any $b\in\mathbb{R}$, including the inhibitory case ($b\leq 0$), there is no global classical solution if the initial data is concentrated enough near $V_F$.

\paragraph{1.3 Results and contributions}In this work, we propose a generalized solution of \eqref{eq:NNLIF-intro} from a new perspective, and throughout this paper we always take $a_1>0$. We formally derive from \eqref{eq:NNLIF-intro} an equation in a new timescale linked with the firing rate $N(t)$. This equation can be shown to be globally well-posed for any $b\in\mathbb{R}$, {even in the event of the blow-up of $N(t)$.} The generalized solution of \eqref{eq:NNLIF-intro} is then defined via the backward change of timescale, {and in particular} it may have a jump when the firing rate blows up. We establish properties of the generalized solution including the characterization of blow-ups (Proposition \ref{prop:jump}) and the global well-posedness in the original timescale (Theorem \ref{thm:gwp-t}).

The idea can be described simply as dividing everything by $N(t)+1$ in the equation. More precisely, we introduce a new 
timescale $\tau$ called the dilated timescale given by 
\begin{equation}\label{intro-dtau-dt}
d\tau=(N(t)+1)dt.
\end{equation}
First we derive formally an equation in timescale $\tau$, which allows $N=+\infty$ and is globally well-posed in the classical sense. The generalized solution in the original timescale $t$ is then defined via a (non-trivial) timescale transform from the classical solution in timescale $\tau$. It turns out when the firing rate blows up, i.e. $N=+\infty$, the equation in timescale $\tau$ may persist for a finite or infinite time interval before the firing rate returns to a finite value. And by \eqref{intro-dtau-dt}  such an interval in the dilated time scale is mapped onto a single point in the original time variable. Therefore the generalized solution in timescale $t$ might have a jump at a blow-up time although the evolution in timescale $\tau$ is always continuous. Such a jump at a single time is ``dilated'' to a continuous dynamics in a time interval in timescale $\tau$ which explains our name for this dilated timescale. .

One of the main results on this generalized solution is the global well-posedness when $b<V_F-V_R$ with counter-examples for $b\geq V_F-V_R$ (Theorem \ref{thm:gwp-tau}). While in the dilated timescale $\tau$ the solution is always global, the global well-posedness in the original timescale $t$ boils down to investigating the long time behavior of the equation in timescale $\tau$. Roughly speaking, the solution fails to be global in timescale $t$ if it is trapped in the blow-up regime in the dilated timescale $\tau$. 

We stress that $a_1>0$ is crucial for the well-posedness of our generalized solution. It ensures the equation in timescale $\tau$ is uniformly parabolic. Hence, we view $a_1>0$ as an essential parameter that yields the well-posedness of the generalized solution. This is in contrast to the fact that $a_1>0$ deteriorates the well-posedness in the classical theory \cite[Appendix A]{Antonio_Carrillo_2015}.

As a new perspective, this generalized solution gives several immediate implications on understanding the NNLIF dynamics as well as opens new questions. {First, it provides a way to understand the dynamics when the firing rate blows up as well as the continuation of the solution after a blow-up}: first solve the equation in $\tau$ and then obtain the solution in $t$ by the inverse change of timescale. Second, the equation in $\tau$ may be a suitable platform to study the long time behavior of the NNLIF dynamics \eqref{eq:NNLIF-intro}. It is globally well-posed in the classical sense and has at least one steady state for all $b\in\mathbb{R}$. Third, the idea of this generalized solution, namely introducing the dilated timescale, might be applied to other integrate-and-fire models, including the model with a refractory state \cite{CACERES201481}, population density models of integrate-and-fire neurons with jumps \cite{Dumont2013Population,dumont2013synchronization,dumont2020mean}, and a model derived via a fast conductance limit \cite{carrillo2022simplified}.

 There are also other possible ways to define generalized solutions and it is interesting to investigate relations between various generalized solutions. In particular, in {some} recent works \cite{taillefumier2022characterization,sadun2022global} the authors have studied a similar model to \eqref{eq:NNLIF-intro}, and also used the idea of introducing a new timescale related with the firing rate $N(t)$. However, there are some subtle but essential differences between their generalized solution and ours. We shall give a discussion on issues in defining a generalized solution in Section \ref{sec:discussion}. 
 
The rest of this paper is arranged as follows. Section \ref{sec:2} is a detailed introduction to our generalized solution and main results. Proofs are given in Section \ref{sc:tau-gwp}, \ref{sec:4e} and \ref{sec:5i}. Precisely, in Section \ref{sc:tau-gwp} we prove the global well-posedness in timescale $\tau$. And the global well-posedness in timescale $t$ is treated in Section \ref{sec:4e} and \ref{sec:5i} for the excitatory and the inhibitory case, respectively. Finally, additional discussions are given in Section \ref{sec:discussion}.

\section{The generalized solution and main results}\label{sec:2}
This section is a detailed introduction to the generalized solution of the NNLIF model \eqref{eq:NNLIF-intro}. Heuristic explanations and precise definitions are given in Section \ref{sec:2.1}. Section \ref{sec:2.2gwp} and Section \ref{sec:2.3blow} are devoted to the global well-posedness result and the characterization of the blow-up for the generalized solution, respectively.

We recall the NNLIF model \eqref{eq:NNLIF-intro}. Substituting the definitions of $\mu(t)$ and $a(t)$ \eqref{intro-mu-a}, we have
\begin{equation}\label{eq:NNLIF-original}
\begin{aligned}
    \p_tp+\p_v[(-v+b_0+bN(t))p]=(a_0+a_1N(t))\p_{vv}p+N(t)\delta_{v=V_R},\quad v\in(-\infty,V_F),t>0,\\
\end{aligned}
\end{equation}
where the constant drift $b_0$ is a sum of the leaky potential and the external input
\begin{equation}\label{def-b0}
    b_0:=V_L+\mu_0.
\end{equation} 
The boundary and initial conditions of \eqref{eq:NNLIF-original} are given by
\begin{equation}\label{eq:NNLIF-bc}
    p(V_F,t)=0,\ p(-\infty,t)=0,\ t>0,\quad\quad\quad p(v,0)=p_I(v),\ v\in(-\infty,V_F).
\end{equation}And the firing rate $N(t)$ is determined by \eqref{intro-Nt}
\begin{equation}\label{def-N-eq}
    N(t)=-(a_0+a_1N(t))\p_vp(V_F,t).
\end{equation}
We assume the parameters satisfy
\begin{equation}
    V_L\leq V_R<V_F,\quad \mu_0\in\mathbb{R},\quad b\in\mathbb{R},\quad a_0>0,a_1>0.
\end{equation}
It is worth emphasizing that we work with the case $a_1>0$, which imposes a finite upper bound $-p(V_F,t)< 1/a_1$ to get a finite and non-negative firing rate from \eqref{def-N-eq}. We might interpret the solution of \eqref{def-N-eq} as
\begin{equation}\label{def-N}
\begin{aligned}
        N(t):=\frac{-a_0\p_vp(V_F,t)}{(1+a_1\p_vp(V_F,t))_+}:=\begin{dcases}
    \dfrac{-a_0\p_vp(V_F,t)}{1+a_1\p_vp(V_F,t)},\quad 0\leq -\p_vp(V_F,t)< 1/a_1,\\
    +\infty,\quad  -\p_vp(V_F,t)\geq 1/a_1.
    \end{dcases}
\end{aligned}
\end{equation}Note that the lower bound $-\p_vp(V_F,t)\geq 0$ is expected since the solution is non-negative with $p(V_F,t)=0$ \eqref{eq:NNLIF-bc}. However, in general there is no guarantee of the upper bound $-\p_vp(V_F,t)< 1/a_1$. The possible blow-up of $N(t)$ is the main obstacle towards a global solution of \eqref{eq:NNLIF-original} in the classical sense \cite{carrillo2013CPDEclassical}.

\subsection{Notion of the generalized solution}\label{sec:2.1}
For heuristic purposes, we explain the notion of the generalized solution in an informal way in Section \ref{sec:2.1.1}, and precise definitions and statements follow in Section \ref{sec:212}. 
\subsubsection{Main ideas}\label{sec:2.1.1}
\paragraph{Equation in the dilated timescale $\tau$}
First we rewrite \eqref{eq:NNLIF-original}, substituting \eqref{def-N-eq} for $N(t)$ in the Dirac source term, and get
\begin{equation}\label{eq:NNLIF-original-rw}
\begin{aligned}
    \p_tp+\p_v[(-v+b_0+bN(t))p]=(a_0&+a_1N(t))\p_{vv}p\\& +(a_0+a_1N(t))(-\p_vp(V_F,t))\delta_{v=V_R},\quad v\in(-\infty,V_F),t>0.
\end{aligned}
\end{equation}

The key idea, is to divide both sides of \eqref{eq:NNLIF-original-rw} by $N(t)$. But to avoid the degeneracy from $N(t)=0$, we divide by $N+c$ instead with some $c>0$, whose specific value is unessential, as we shall see in Proposition \ref{prop:gen-uni-c}.

Precisely, we introduce another timescale $\tau$, called \textbf{the dilated timescale}, given by the following change of variable in time 
\begin{equation}\label{change-time-1}
    d\tau=(N(t)+c)dt,\quad
\end{equation}with some fixed constant $c>0$, e.g., $c=1$. 

Then we divide both sides of \eqref{eq:NNLIF-original-rw} by $N(t)+c$ to derive, at least formally, an equation for $\tp(v,\tau(t)):=p(v,t)$ where $\tau(t):=\int_0^{t}(N(s)+c)ds$,
\begin{equation}\label{eq:NNLIF-tau}\begin{aligned}
        \p_{\tau}\tp+\p_v[(-v\tilde{N}+b_c\tilde{N}+b)\tp]=(a_c\tilde{N}&+a_1)\p_{vv}\tp\\& +(a_c\tilde{N}+a_1)(-\p_v\tp(V_F,\tau))\delta_{v=V_R},\quad v\in(-\infty,V_F),\tau>0,
        \end{aligned}
\end{equation} with the same boundary and initial conditions as \eqref{eq:NNLIF-bc}
\begin{equation}\label{eq:NNLIF-bc-tau}
        \tp(V_F,\tau)=0,\ \tp(-\infty,\tau)=0,\ \tau>0,\quad\quad\quad \tp(v,0)=p_I(v),\ v\in(-\infty,V_F).
\end{equation}
In \eqref{eq:NNLIF-tau}, $\tilde{N}$ is a shorthand for $\tilde{N}(\tau)$, which is defined as 
\begin{equation}\label{def-tildeN}
    \tilde{N}(\tau):=\frac{1}{{N}+c}=\frac{(1+a_1\p_v\tp(V_F,\tau))_+}{-a_0\p_v\tp(V_F,\tau)+c(1+a_1\p_v\tp(V_F,\tau))_+}.
\end{equation} 

And the coefficients $b_c,a_c$ are just constants given by \footnote{There is a little inconsistency between the definition of $a_c$ \eqref{def-bcac} and parameters $a_0,a_1$ in \eqref{eq:NNLIF-original-rw}: if we substitute $c=1$ in the definition of $a_c$, it gives $a_0-a_1$ rather than $a_1$. In this paper we either take $c$ to be a fixed but arbitrary positive number or take $c=a_0/a_1$ to set $a_c=0$ (in Section \ref{sec:4e} and \ref{sec:5i}).}	
\begin{equation}\label{def-bcac}
	    b_c:=b_0-cb,\quad\quad\quad a_c:=a_0-ca_1.
\end{equation}
In principle, the dilated timescale $\tau$, and the solution $(\tp,\tilde{N})$ for \eqref{eq:NNLIF-tau} also depend on the parameter $c>0$ in \eqref{change-time-1}. However, as we shall see, the specific choice of $c>0$ won't influence the resulting generalized solution (Proposition \ref{prop:gen-uni-c}). Therefore for convenience, we do not write out the dependence of $c$ for $\tau,\tilde{N}$ and $\tp(v,\tau)$ explicitly unless necessary.

We refer to \eqref{eq:NNLIF-tau}--\eqref{def-tildeN} as \textbf{the equation in the dilated timescale}, in contrast to \eqref{eq:NNLIF-original}, the equation in the original timescale. Several observations indicate that \eqref{eq:NNLIF-tau} is more favorable than  \eqref{eq:NNLIF-original} for the well-posedness. From \eqref{def-tildeN} we know
\begin{equation}\label{up-low-bd-tilde}
    \tilde{N}(\tau)\in[0,1/c],
\end{equation}thanks to $-\p_v\tp(V_F,\tau)\geq 0$ since $p\geq 0$ and $\tp(V_F,\tau)=0$. Therefore in \eqref{eq:NNLIF-tau} the source of nonlinearity $\tilde{N}(\tau)$ is uniformly bounded, in contrast to $N(t)$ in the equation in the original timescale $t$ \eqref{eq:NNLIF-original}. Moreover, \eqref{eq:NNLIF-tau} is uniformly parabolic since the diffusion coefficient \begin{equation}\label{a-sec211-uniform}
    a_c\tilde{N}(\tau)+a_1=\tilde{N}(\tau)a_0+(1-c\tilde{N}(\tau))a_1
\end{equation}
is between $a_1>0$ and $\frac{1}{c}a_0>0$. 

In fact, we can establish the {global well-posedness} for {classical} solutions of the equation in the dilated timescale $\tau$ \eqref{eq:NNLIF-tau}, for any $b\in\mathbb{R}$. Precise statements are postponed to Section \ref{sec:212}. 

In the classical sense, the above derivation from \eqref{eq:NNLIF-original} to \eqref{eq:NNLIF-tau} is rigorous only if $N<+\infty$, since otherwise when $N(t)=+\infty$ the change of variable \eqref{change-time-1} may not be well-defined. Recall that the blow up of the firing rate $N(t)$ is exactly the obstacle for the existence of a global classical solution of \eqref{eq:NNLIF-original} \cite{carrillo2013CPDEclassical}.  However, in the dilated timescale $\tau$, $N=+\infty$ is allowed since it just corresponds to $\tilde{N}=0$ by \eqref{def-tildeN}, which does not hinder the evolution of \eqref{eq:NNLIF-tau}. In view of its global well-posedness, \eqref{eq:NNLIF-tau} may evolve through a period when $N=+\infty$ (which corresponds to $\tilde{N}=0$) and go back to the non-blow up regime {when} $N<+\infty$ (which corresponds to $\tilde{N}>0$). Hence, if we change back from $\tau$ to the original timescale $t$, we may be able to define a generalized solution of \eqref{eq:NNLIF-original} which goes beyond the blow-up time, thus having a longer lifespan. We shall call this ``change-back process'' the inverse time-dilation transform (see the remark below \eqref{def-N-2} or Definition \ref{def:generaized} for a precise definition), {while the process from \eqref{eq:NNLIF-original-rw} to obtain \eqref{eq:NNLIF-tau} by dividing by $N+c$ is called a time-dilation transform.} In the next, we proceed to examine the inverse time-dilation transform in details to define a generalized solution.

\paragraph{(Inverse) time-dilation transform and the generalized solution}Let $(\tp,\tilde{N})$ be a classical solution of the equation in the dilated timescale $\tau$ \eqref{eq:NNLIF-tau}. We aim to define generalized solutions for the equation in the original timescale $t$ \eqref{eq:NNLIF-original}, by ``changing back'' from timescale $\tau$ to timescale $t$. In view of the change of variable \eqref{change-time-1} we define
\begin{equation}\label{change-time-2}
    dt=\frac{1}{N+c}d\tau=\tilde{N}d\tau,\quad\quad\quad t(\tau):=\int_0^{\tau}\tilde{N}(\omega)d\omega,
\end{equation}which gives each $\tau$ a $t=t(\tau)$. We shall assign each $t$ a $\tau=\tau(t)$ by inverting \eqref{change-time-2}. However, regarding that $\tilde{N}\in [0,{1}/{c}]$ as in \eqref{up-low-bd-tilde}, the map $t(\tau)$ in \eqref{change-time-2} is only non-decreasing, not necessarily strictly increasing. Indeed, if $\tilde{N}(\tau)\equiv0$ on an interval $[\tau_1,\tau_2]$, then $t(\tau)\equiv t(\tau_1)$ for all $\tau\in[\tau_1,\tau_2]$, which makes it indefinite to define the inverse $\tau(t)$.

\begin{figure}[htb]
    \centering
    \includegraphics[width=1.0\textwidth]{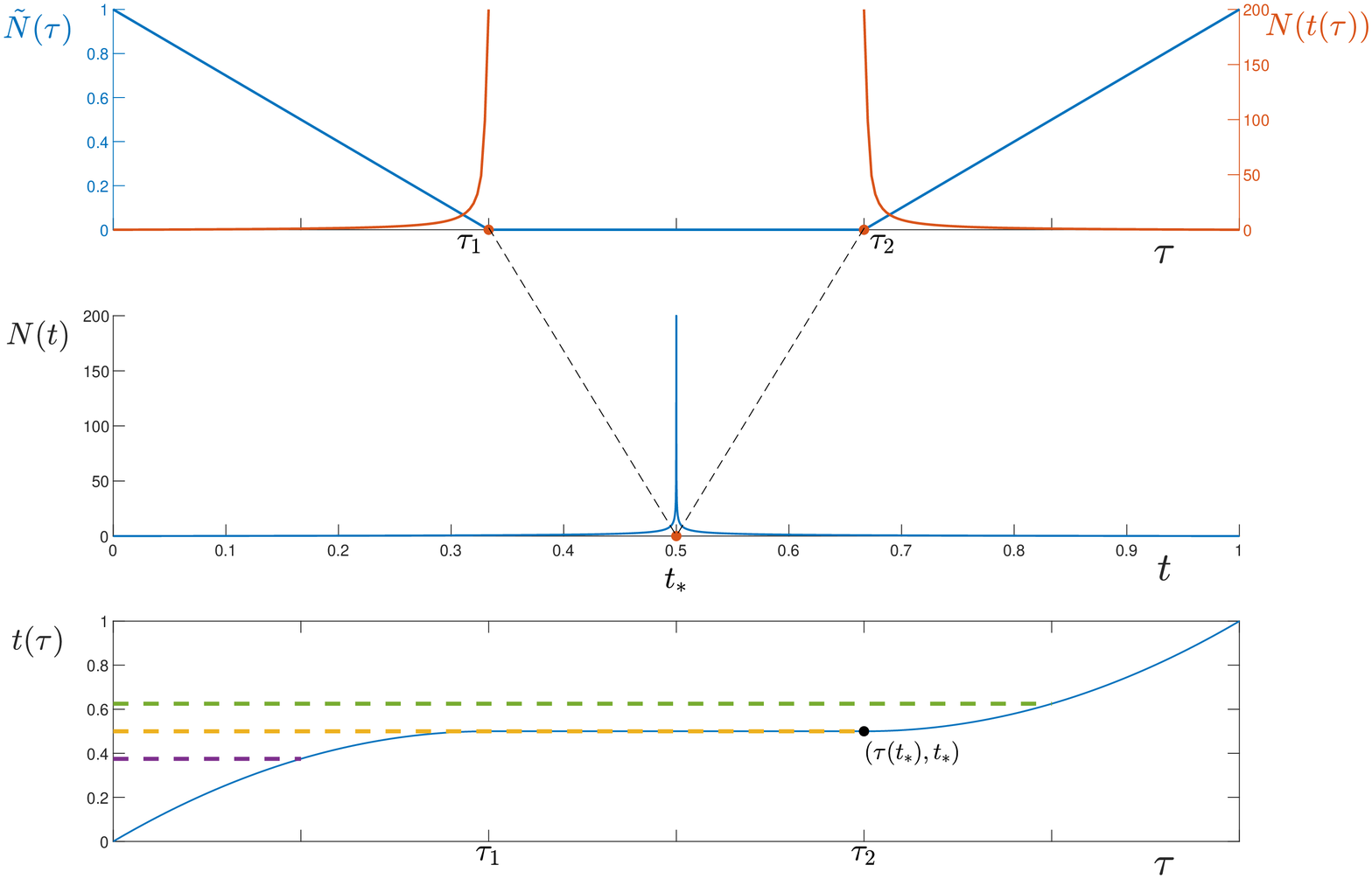}
    \caption{Schematic plots to illustrate the definition of the generalized solution. Upper: $\tilde{N}$ and the corresponding $N$ in the dilated timescale $\tau$. On $[\tau_1,\tau_2]$, $\tilde{N}\equiv0$ and $N\equiv+\infty$. Middle: The firing rate $N(t)$ in the original timescale $t$. The blow-up only happens at the time $t_*=0.5$. The blow-up time $t_*$ is ``dilated'' to the interval $[\tau_1,\tau_2]$ in the original timescale, as indicated by the dash lines connecting the upper plot and the middle plot. Lower: The map $t(\tau)$ \eqref{change-time-2} associated with $\tilde{N}$. For $\tau\in[\tau_1,\tau_2]$, ${t}(\tau)\equiv t(\tau_2)=:t_*$. Dashed lines indicate the change-back map $\tau(t)$ \eqref{change-back-map}.  }
    \label{fig:1-schematic}
\end{figure}

To remove the ambiguity, we consider the following generalized inverse of \eqref{change-time-2}
\begin{equation}\label{change-back-map}
    \tau(t):=\sup\left\{\tau:\  \int_0^{\tau}\tilde{N}(\omega)d\omega=t\right\}.
\end{equation} In principle, picking another $\tau$ with $\int_0^{\tau}\tilde{N}(\omega)d\omega=t$ shall not make an essential difference, but to fix the choice \eqref{change-back-map} is convenient, which makes the map $\tau(t)$ càdlàg. An illustration is given in Figure \ref{fig:1-schematic}. 

We define the generalized solution of \eqref{eq:NNLIF-original} at time $t$ by assigning the value of the classical solution of \eqref{eq:NNLIF-tau} at $\tau=\tau(t)$. Precisely, we define the generalized solution at time $t$ to be
\begin{equation}\label{211-pgen}
    p(v,t):=\tp(v,\tau(t)),
\end{equation}
 and its associated firing rate is given by 
 \begin{equation}\label{211-Ngen}
     N(t):=\left(\frac{1}{\tilde{N}(\tau(t))}-\frac{1}{c}\right)\in[0,+\infty],
 \end{equation} which is from inverting the definition of $\tilde{N}$ \eqref{def-tildeN}. Recalling the definition of $\tilde{N}$ \eqref{def-tildeN}, we see that \eqref{211-pgen} and \eqref{211-Ngen} give the same firing rate as \eqref{def-N}
 \begin{equation}\label{def-N-2}\begin{aligned}
     N(t)&=\frac{-a_0\p_vp(V_F,t)}{(1+a_1\p_vp(V_F,t))_+}=\begin{dcases}
    \frac{-a_0\p_vp(V_F,t)}{1+a_1\p_vp(V_F,t)},\quad 0\leq -\p_vp(V_F,t)< 1/a_1,\\
    +\infty,\quad  -\p_vp(V_F,t)\geq 1/a_1.
    \end{dcases}
\end{aligned}
 \end{equation} In \eqref{211-pgen} and \eqref{211-Ngen}, we shall call $(\tp,\tilde{N})$ the \textbf{time-dilation transform} of $(p,N)$, and call $(p,N)$ the \textbf{inverse time-dilation transform} of $(\tp,\tilde{N})$.
 \vspace{1.5em}

\paragraph{The blow-up of $N$} Let's examine the behavior of the generalized solution when the firing rate $N$ blows up. We shall look at the zeros of $\tilde{N}$, since by \eqref{def-tildeN}
\begin{equation}\label{tildeN-0-N-inf}
    \tilde{N}=0 \iff N=+\infty.
\end{equation}
Suppose the interval $[\tau_1,\tau_2]=\{\tau\geq0:\tilde{N}(\tau)=0\}$ are all zeros of $\tilde{N}$ in the dilated timescale. Then we know from \eqref{change-time-2} $t(\tau_1)=t(\tau_2)=:t_*$ and by the choice \eqref{change-back-map} $$\tau(t_*^+)=\tau(t_*)=\tau_2>\tau_1=\tau(t_*^{-}).$$ Hence in the original timescale, the probability density has a jump at $t=t_*$ as
\begin{equation}\label{jump-sec211}
    p(\cdot,t_*^+)=p(\cdot,t_*)=\tp(\cdot,\tau_2)\neq \tp(\cdot,\tau_1)=p(\cdot,t_*^-).
\end{equation}
 And the firing rate in the timescale $t$ blows up just at $t=t_*$, i.e., $N(t_*)=+\infty$ but for $t\neq t_*$ $N(t)<+\infty$, see Figure \ref{fig:1-schematic}.

The jump of $p$ comes from the fact that the whole interval $[\tau_1,\tau_2]$ is compressed to one point by \eqref{change-back-map}, due to $\tilde{N}=0$. In timescale $\tau$, the blow up time $t_*$ is ``dilated'' to an interval $[\tau_1,\tau_2]$. This is why we call $\tau$ a ``dilated timescale''. In this ``dilated'' interval $[\tau_1,\tau_2]$, while $N\equiv+\infty$, or equivalently $\tilde{N}\equiv0$, the density function $\tp$ evolves continuously. A more detailed discussion on the relationship between $p(\cdot,t_*)$ and $p(\cdot,t_*^-)$ is given in Section \ref{sec:2.3blow}. 

\paragraph{Lifespan}While the solution in the dilated timescale $\tau$ is always global, the generalized solution in the original timescale $t$ may be not. By definition, we can define a generalized solution at time $t$ if and only if the map in \eqref{change-back-map} $\tau(t)<+\infty$, or equivalently
\begin{equation*}
    t<\int_0^{+\infty}\tilde{N}(\tau)d\tau.
\end{equation*}
Therefore, the maximal lifespan of a generalized solution is exactly given by
\begin{equation}\label{lifespan-2.1.1}
    T^*:=\int_0^{+\infty}\tilde{N}(\tau)d\tau.
\end{equation}
{And then the global well-posedness for the generalized solution in the original timescale $t$ is equivalent to $T^*=+\infty$}. Of course, the integral in \eqref{lifespan-2.1.1} does not depend on the choice of $c>0$ either (Proposition \ref{prop:gen-uni-c}). 

\begin{remark}
The reformulation from \eqref{eq:NNLIF-original} to \eqref{eq:NNLIF-original-rw} is indeed essential for the construction of the generalized solution. It turns out that although \eqref{eq:NNLIF-original} and \eqref{eq:NNLIF-original-rw} are equivalent when $N(t)<+\infty$, they give different time dilation transforms when $N=+\infty$. We will discuss more on this in Section \ref{sec:discussion}.
\end{remark}

\subsubsection{Precise definitions and properties}\label{sec:212}

In this section we give the precise definition of the generalized solution of \eqref{eq:NNLIF-original} and its basic properties.

We start with the definition of the classical solution for the equation in the dilated timescale $\tau$ \eqref{eq:NNLIF-tau}, and the statement of its global well-posedness.

Following {definitions} of classical solutions for the NNLIF model in time variable $t$ \eqref{eq:NNLIF-original} \cite{carrillo2013CPDEclassical,roux2021towards}, we give the definition of a classical solution for the NNLIF model in $\tau$ variable \eqref{eq:NNLIF-tau}--\eqref{def-tildeN}.

\begin{definition}[Classical solution in the dilated timescale $\tau$]\label{def:classical-tau}We say the pair $(\tp,\tilde{N})$ is a classical solution for the system in timescale $\tau$ \eqref{eq:NNLIF-tau}--\eqref{def-tildeN} with a given parameter $c>0$ on the time interval $[0,\bar{\tau})$, $0<\bar{\tau}\leq+\infty$ if
\begin{enumerate}
    \item(Regular, integrable, and non-negative $\tp$) $\tp\in C^0\bigl((-\infty,V_F]\times[0,\bar{\tau})\bigr)\cap C^{2,1}\bigl(((-\infty,V_R)\cup(V_R,V_F))\times[0,\bar{\tau})\bigr)\cap L^{\infty}\bigl([0,\bar{\tau}),L^1_+((-\infty,V_F])\bigr)$.
    \item(Boundary flux) For $\tau\in[0,\bar{\tau})$, the one-sided derivatives $\p_v\tp(V_R^{-},\tau),\p_v\tp(V_R^{+},\tau)$ and $\p_v\tp(V_F^{-},\tau)$ are well-defined and finite.
    \item(Decay and Boundary condition) For $\tau\in[0,\bar{\tau})$, $\tp(v,\tau),\p_v\tp(v,\tau)$ goes to zero as $v\rightarrow-\infty$, and $\tp(V_F,\tau)=0$.
    \item(Continuous $\tilde{N}$) $\tilde{N}(\tau)$ is a continuous function from $[0,\bar{\tau})$ to $[0,\frac{1}{c}]$ and the relationship \eqref{def-tildeN} is satisfied for $\tau\in[0,\bar{\tau})$.
    \item $(\tp,\tilde{N})$ satisfies the equation \eqref{eq:NNLIF-tau} in the classical sense on $(-\infty,V_R)\cup(V_R,V_F)$ and in the distributional sense on $(-\infty,V_F)$.
    \item (Initial condition) $\tp(v,0)=p_I(v)$ for $v$ in $(-\infty,V_F]$. 
\end{enumerate}
\begin{remark}\label{rmk:flux-jump}
For $(\tp,\tilde{N})$ with regularities given in Definition \ref{def:classical-tau}, to satisfy \eqref{eq:NNLIF-tau} in the distributional sense on $(-\infty,V_F)$ is equivalent to the requirement that \eqref{eq:NNLIF-tau} is satisfied in the classical sense on $(-\infty,V_R)\cup(V_R,V_F)$ in addition to the following flux jump condition at $V_R$
\begin{equation}
   -(a_c\tilde{N}+a_1)\p_v\tp(V_R^{+},\tau)+(a_c\tilde{N}+a_1)\p_v\tp(V_R^{-},\tau)=-(a_c\tilde{N}+a_1)\p_v\tp(V_F^{-},\tau).
\end{equation}
\end{remark}
\end{definition}
To state the global well-posedness for the classical solution of \eqref{eq:NNLIF-tau}, we need a natural assumption on the initial data as in \cite{carrillo2013CPDEclassical,roux2021towards}. 

\begin{assumption}\label{as:classical-init}
For the initial data $p_I(v)$ we impose the following assumptions.
\begin{enumerate}
    \item $p_I\in L^1((-\infty,V_F))$ is a probability density function, i.e., $p_I(v)\geq 0$ and $\int_{-\infty}^{V_F}p_I(v)dv=1$.
    \item $p_I\in C^0((-\infty,V_F])\cap C^1((-\infty,V_R)\cup(V_R,V_F))$ and is non-negative.
    \item The one-sided derivatives $\frac{d}{dv}p_I(v)$ at $V_R^{-},V_R^{+}$ and $V_F^{-}$ are well defined and finite.
    \item $p_I(v),\p_vp_I(v)$ goes to zero as $v\rightarrow-\infty$, and $p_I(V_F)=0$.
\end{enumerate}
\end{assumption}

With this assumption, we show the global well-posedness for the classical solution of the equation in the dilated timescale $\tau$ \eqref{eq:NNLIF-tau}. 
\begin{theorem}[Global well-posedness in the dilated timescale $\tau$]\label{thm:gwp-tau}
Let the initial data $p_I(v)$ satisfy Assumption \ref{as:classical-init}, and then for any given $c>0$, the system in timescale $\tau$ \eqref{eq:NNLIF-tau}--\eqref{def-tildeN} has a unique classical solution defined in Definition \ref{def:classical-tau}, which is global, i.e., exists on the time interval $[0,+\infty)$.
\end{theorem}
The proof of Theorem \ref{thm:gwp-tau}, given in Section \ref{sc:tau-gwp}, is an adaption of \cite{carrillo2013CPDEclassical}. As discussed in Section 2.1.1, the boundedness of $\tilde{N}$ and the uniform parabolicity thanks to $a_1>0$ play important roles in the proof. 

We define generalized solutions for the equation in the original timescale $t$ \eqref{eq:NNLIF-original}, by applying the transform \eqref{change-back-map} to classical solutions in the dilated timescale $\tau$ \eqref{eq:NNLIF-tau}, as discussed in Section 2.1.1. The precise definition of our generalized solution is given as follows.

\begin{definition}[Generalized solution in timescale $t$]\label{def:generaized}We say the pair $(p,N)$ is a generalized solution for \eqref{eq:NNLIF-original},\eqref{eq:NNLIF-bc} and \eqref{def-N-eq} on time interval $[0,\bar{T})$ if there exists a classical solution $(\tp,\tilde{N})$ for the time-dilated problem \eqref{eq:NNLIF-tau}--\eqref{def-tildeN} on time interval $[0,\bar{\tau})$ with parameter $c>0$ such that
\begin{enumerate}
    \item $\bar{T}=\int_0^{\bar{\tau}}\tilde{N}(\tau)d\tau$.
    \item $p(v,t)=\tp(v,\tau(t))$ for $v\in(-\infty,V_F]$ and $t\in[0,\bar{T})$.
    \item $N(t)=\begin{cases}\frac{1-c\tilde{N}(\tau(t))}{\tilde{N}(\tau(t))},\quad \tilde{N}(\tau(t))\in(0,\frac{1}{c}],\\
    +\infty,\quad \tilde{N}(\tau(t))=0.
    \end{cases}$
\end{enumerate}
Here $\tau(t)$ is the change of variable defined in \eqref{change-back-map}.
\end{definition} 
We refer to $(\tp,\tilde{N})$ as the \textbf{time-dilation transform} of $(p,N)$, and correspondingly, $(p,N)$ as the \textbf{inverse time-dilation transform} of $(\tp,\tilde{N})$. The generalized solution defined in Definition \ref{def:generaized} may be called a \textbf{time-dilated solution} of \eqref{eq:NNLIF-original}.

\begin{remark}Definition \ref{def:generaized} implies that the density function $p$ and the firing rate $N$ satisfies the relation \eqref{def-N-2}.
\end{remark}

By the uniqueness of \eqref{eq:NNLIF-tau} (Theorem \ref{thm:gwp-tau}), we shall show that the generalized solution does not depend on the value of $c>0$ and its uniqueness. Moreover, we can always choose $\bar{\tau}=+\infty$ in the definition thanks to the global well-posedness of \eqref{eq:NNLIF-tau} (Theorem \ref{thm:gwp-tau}), and get the maximal lifespan of the generalized solution
\begin{equation*}
    T^*=\int_0^{+\infty}\tilde{N}(\tau)d\tau,
\end{equation*} as in \eqref{lifespan-2.1.1}. These are summarized in the following proposition.

\begin{proposition}\label{prop:gen-uni-c}
(i)Let the initial data satisfy Assumption \ref{as:classical-init}. Then the generalized solution in Definition \ref{def:generaized} on time interval $[0,\bar{T})$, if exists, is unique and does not depend on the particular choice of $c>0$. (ii) Moreover, its maximal lifespan $T^*$ is given by \eqref{lifespan-2.1.1} which does not depend on the choice of $c>0$ either.
\end{proposition}
\begin{proof}[Proof of Proposition \ref{prop:gen-uni-c}]Suppose $(q_1,\tilde{N}_1)$ and $(q_2,\tilde{N}_2)$ are two solutions for \eqref{eq:NNLIF-tau} with different $c>0$, denoted as $c_1,c_2>0$ respectively, but sharing the same initial data. Let $\tau_1(t)$ and $\tau_2(t)$ be the corresponding maps \eqref{change-back-map}, then for (i) it suffices to show
\begin{equation}\label{tmp-notdependc}
    q_1(\cdot,\tau_1(t))=q_2(\cdot,\tau_2(t)),\quad \forall t\in[0,\bar{T}).
\end{equation}Indeed, this follows from a change of timescale from $\tau_1$ to $\tau_2$, which is well-defined, and the uniqueness of \eqref{eq:NNLIF-tau} for a fixed $c$. To give details, starting from $q_1(\tau_1,\cdot)$, we consider the change of timescale 
\begin{equation}\label{change-c1c2}
    d\tau_2=\frac{N+c_1}{N+c_2}d\tau_1=\frac{1}{1+(c_2-c_1)\tilde{N}(\tau_1)}d\tau_1,
\end{equation}
 which is well-defined since $\frac{N+c_1}{N+c_2}$ is between $1$ and $\frac{c_1}{c_2}>0$. Define the corresponding map
\begin{equation*}
    h(\tau_1):=\int_{0}^{\tau_1}\frac{1}{1+(c_2-c_1)\tilde{N}(\tau_1)}d\tau_1.
\end{equation*} Then the new density function induced from $q_1$ by the change of timescale \eqref{change-c1c2} is given by
\begin{equation}
    \hat{q}_1(\cdot,h(\tau_1)):=q_1(\cdot,\tau_1).
\end{equation} By its construction, $\hat{q}_1$ solves \eqref{eq:NNLIF-tau} with $c=c_2>0$. Hence, by its uniqueness given in Theorem \ref{thm:gwp-tau}, we have
\begin{equation*}
 q_2(\cdot,h(\tau_1))=\hat{q}_1(\cdot,h(\tau_1)):=q_1(\cdot,\tau_1).
\end{equation*}
Finally, noting that $h$ is a strictly increasing function, one directly checks $\tau_2(t)=h(\tau_1(t))$, therefore \eqref{tmp-notdependc} holds. (Here with abuse of notation, we use $\tau_i$ to denote the time variables and $\tau_i(t)$ for the maps \eqref{change-back-map}.) Hence we have shown that the generalized solution does not depend on the choice of $c>0$, its uniqueness also follows from the uniqueness of \eqref{eq:NNLIF-tau} (or taking $c_1=c_2$ in the above proof).

By Definition \ref{def:generaized} and the global well-posedness of \eqref{eq:NNLIF-tau}, the maximal lifespan is given by $T^*=\int_0^{+\infty}\tilde{N}(\tau)d\tau$ as in \eqref{lifespan-2.1.1}. For two different time-dilation transforms $\tilde{N}_1$, $\tilde{N}_2$, we have
\begin{equation*}
    \int_0^{+\infty}\tilde{N}_1(\tau_1)d\tau_1= \int_0^{+\infty}\tilde{N}_1\frac{\tilde{N_2}}{\tilde{N_1}}d\tau_2=\int_0^{+\infty}\tilde{N}_2(\tau_2)d\tau_2,
\end{equation*}noting that $d\tau_1=\frac{\tilde{N}_2}{\tilde{N}_1}d\tau_2$ as in \eqref{change-c1c2}.

\end{proof}
Next, we examine the connection between the generalized solution and the classical solution of \eqref{eq:NNLIF-original}. Here the definition of a classical solution for \eqref{eq:NNLIF-original} is similar to Definition \ref{def:classical-tau} with an additional requirement that the corresponding $N(t)<+\infty$, or equivalently $-a_1\p_vp(V_F,t)<1$.
\begin{proposition}[Relation with the classical solution]\label{prop:classical}
A classical solution for \eqref{eq:NNLIF-original} is always a generalized solution. On the other hand, a generalized solution $(p,N)$ for \eqref{eq:NNLIF-original} on $[0,\bar{T})$ {where} $0<\bar{T}\leq T^*$, is also a classical solution if and only if $N(t)<+\infty$ for all $t\in [0,\bar{T})$.
\end{proposition}
\begin{proof}[Proof of Proposition \ref{prop:classical}]
For a classical solution of \eqref{eq:NNLIF-original} we have $0\leq N(t)<+\infty$ hence the change of variable $d\tau=(N(t)+c)dt$ in \eqref{change-time-1} is well-defined. Then one can transform it to the dilated timescale \eqref{eq:NNLIF-tau} to get its associated time-dilation transform $(\tp,\tilde{N})$ in Definition \ref{def:generaized}.

For a generalized solution $(p,N)$ on $[0,\bar{T})$ to be a classical solution, clearly it has to satisfy that $N(t)<+\infty$ for $t\in[0,\bar{T})$. On the other hand, suppose a generalized solution indeed satisfy $N(t)<+\infty$ for all $t\in [0,\bar{T})$. This implies $\tilde{N}(\tau)>0$ for all $\tau\in[0,\bar{\tau})$ in the dilated timescale. Hence, the map $t(\tau)$ \eqref{change-time-2}, and thus its inverse ``change-back-map'' $\tau(t)$ \eqref{change-back-map} are both $C^1$ strictly increasing function. Therefore, by definition the inverse time-dilation transform preserves all regularities required for $(p,N)$ to be a classical solution.  
\end{proof}

In general, a generalized solution $p$ always has the same spatial regularity as a classical solution, but in time we only expect it to be càdlàg when $N(t)=+\infty$. We summarize the regularity of a generalized solution as follows.
\begin{proposition}[Regularity of the generalized solution]\label{prop:regularity}
For a generalized solution $(p,N)$ of \eqref{eq:NNLIF-original} on $[0,\bar{T})$ {where} {$0<\bar{T}\leq T^*$}, we have
\begin{enumerate}
    \item $p$ has the same regularity in $v$ as classical solutions.
    \item In the original timescale $t$, $p$ is càdlàg, i.e., right-continuous with left limits.
    \item  $N$ is a continuous function\footnote{As usual, the topology of $[0,+\infty]$ is generated by open sets in $[0,+\infty)$ in addition to $(M,+\infty]$, for all $M\in[0,+\infty)$.} from $[0,\bar{T})$ to $[0,+\infty]$. 
    \item If $N(t)<+\infty$ for some $0\leq t<\bar{T}$, then around $t$ locally in time $p$ has the same regularity as a classical solution, in particular $p$ shall be $C^1$ in time in {a neighbourhood of} $t$. 
\end{enumerate}
In particular, if $N(0)<+\infty$, then $p(\cdot,t=0)$ coincides with the initial data $p_I(\cdot)$. Otherwise the initial data $p_I(v)$ shall be interpreted as $p_I(\cdot)=p(\cdot,0^-)$.
\end{proposition}
\begin{proof}[Proof of Proposition \ref{prop:regularity}] These are all direct consequences of Definition \ref{def:generaized}. Note that the càdlàg property of $p$ follows from the càdlàg property of the map $\tau(t)$ defined in \eqref{change-back-map}.
\end{proof}

For a generalized solution in the original timescale, the density function $p$ can has a jump in time when $N(t)=+\infty$, i.e., the firing rate blows up. We postpone the discussion on this scenario to Section \ref{sec:2.3blow}. We shall discuss the global well-posedness first, along which we introduce the {limit equation} at the blow-up, which is an important notion for the description of the blow-up as well as the global well-posedness.

\subsection{Global well-posedness verus non-existence for the generalized solution}\label{sec:2.2gwp}
\subsubsection{Main results}
Now we state the global well-posedness in the original timescale $t$. Depending on the connectivity parameter $b$, we identify three scenarios:
\begin{enumerate}
    \item $b\geq V_F-V_R$, the strongly excitatory case. 
    \item $0<b<V_F-V_R$, the mildly excitatory case.
    \item $b\leq 0$, the inhibitory case \footnote{Strictly speaking, only $b<0$ is the inhibitory case, but here we also include $b=0$ since they can be treated in a unified way.}.
\end{enumerate}
In short, the global well-posedness for the generalized solution holds if and only if $b<V_F-V_R$.

\begin{figure}[htbp]
    \centering
    \includegraphics[width=0.49\linewidth]{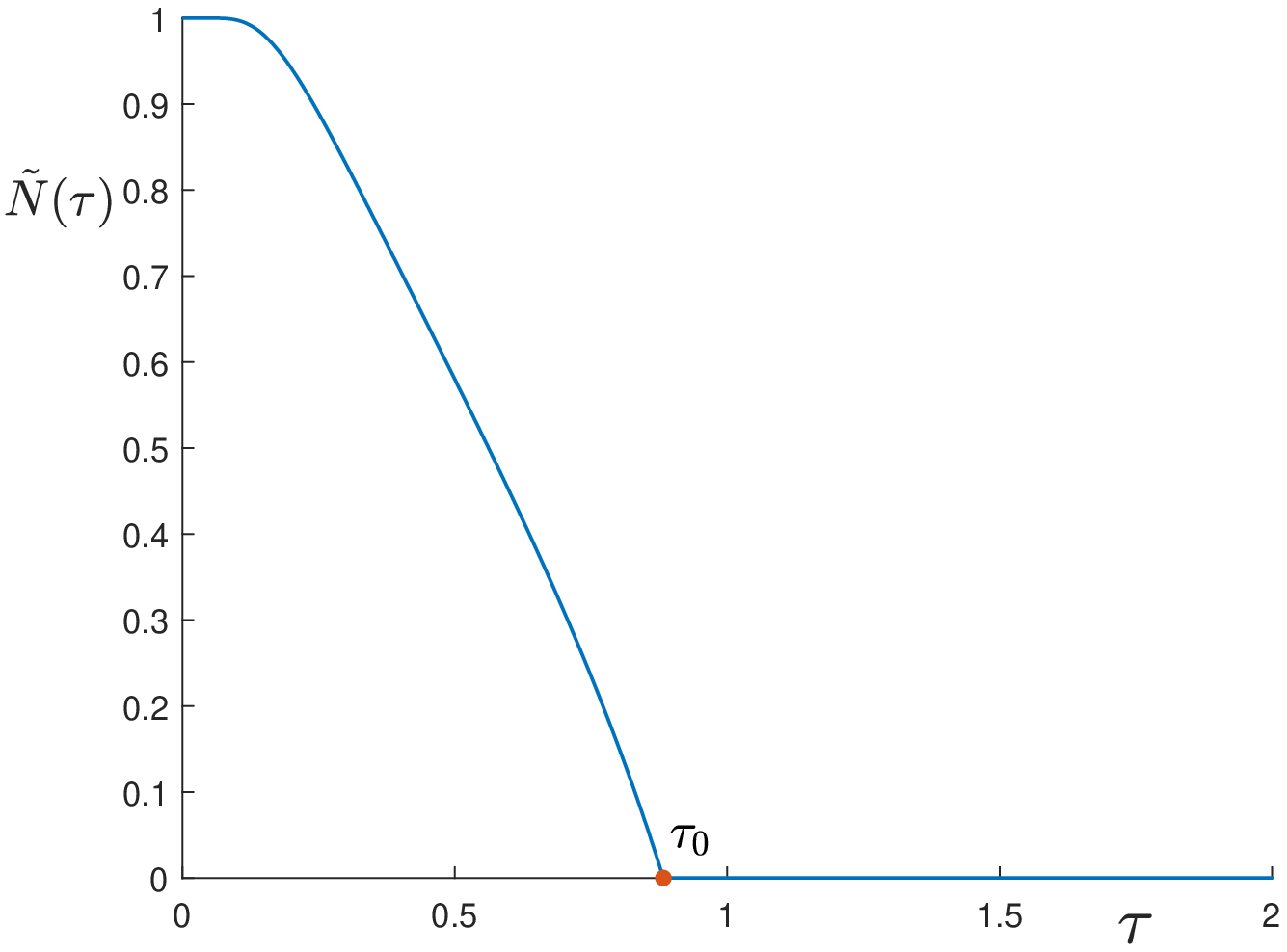}    \includegraphics[width=0.49\linewidth]{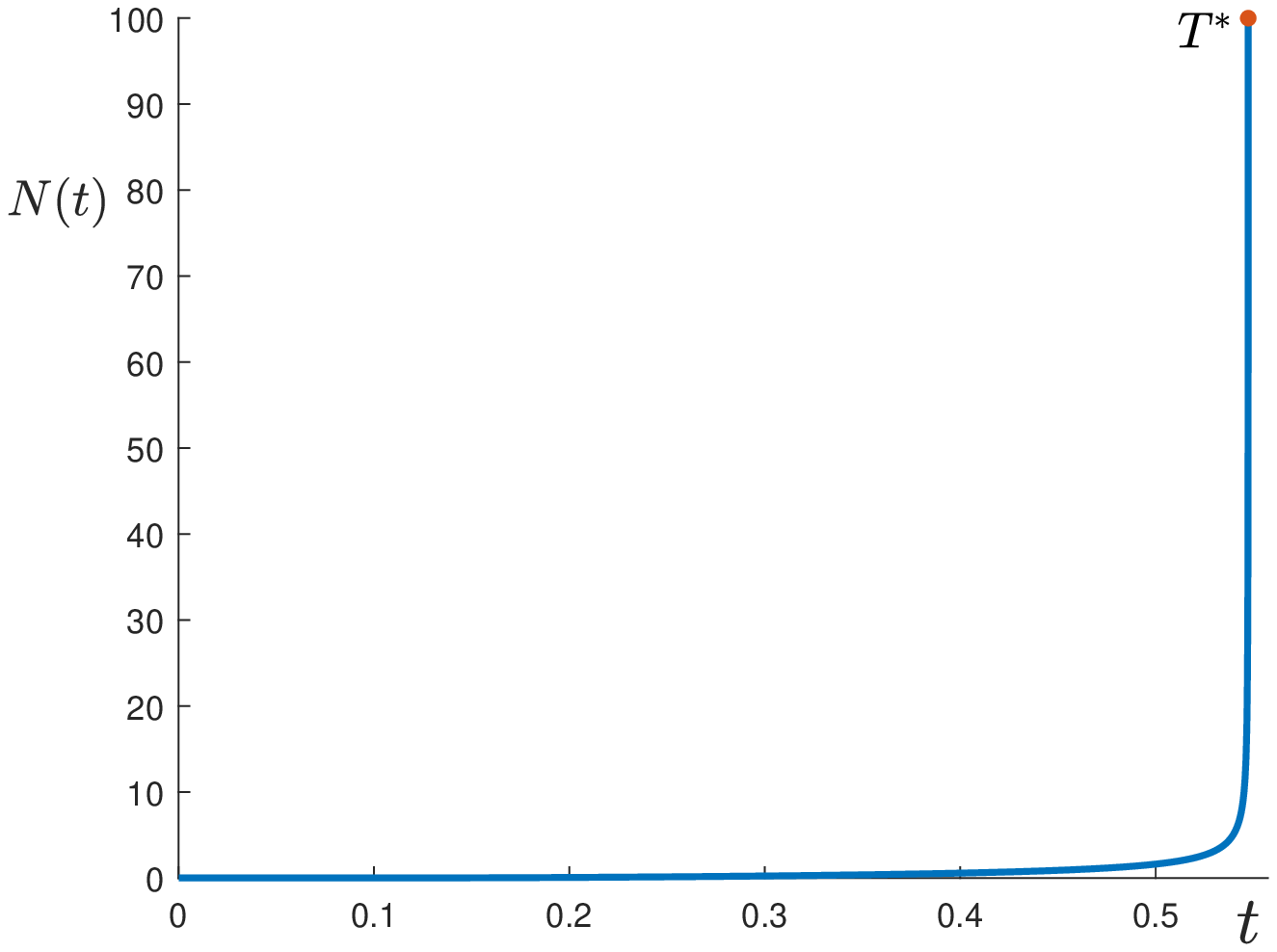}
    \caption{An example of the eternal blow-up (when $b\geq V_F-V_R$) via a numerical simulation. Left: A plot of $\tilde{N}(\tau)$ in the dilated timescale, $\tilde{N}(\tau)\equiv0$ for $\tau\geq \tau_0$, which leads to a finite lifespan $T^*=\int_0^{\tau_0}\tilde{N}(\tau)d\tau$ by \eqref{lifespan-2.1.1}. Right: The corresponding $N(t)$ in the original timescale. The firing rate blows up as the time $t$ approaches the maximal lifespan $T^*$. In this plot, the firing rate is truncated at $N=100$. Parameters: $V_F=1,V_R=0,b=1.5$. The initial data is {a localized Gaussian with variance $0.17$} around $v=-1$. The numerical scheme is a variant version from the one in \cite{hu2021structure}.}
    \label{fig:eternal-blow-up}
\end{figure}

Precisely, we have the following theorem.
\begin{theorem}[Global well-posedness in the original timescale $t$]\label{thm:gwp-t} 
Let the initial data $p_I(v)$ satisfy Assumption \ref{as:classical-init}, there exists a unique generalized solution for the equation in the original timescale $t$ \eqref{eq:NNLIF-original} on the time interval $[0,T^*)$. Here $T^*\in [0,+\infty]$ is the maximal lifespan given by \eqref{lifespan-2.1.1} and we have 
\begin{enumerate}
    \item When $b\geq V_F-V_R$, the strongly excitatory case, there exists an initial data satisfying Assumption \ref{as:classical-init} but giving $T^*=0$, that is, there is no generalized solution on any time interval $[0,t)$ with $t>0$.
    \item When $b<V_F-V_R$ we have $T^*=+\infty$, that is, there exists a unique global generalized solution, provided the following additional assumptions on the drift and the initial data hold.
    \begin{enumerate}
        \item When $0<b<V_F-V_R$, the mildly excitatory case, we additionally assume
        \begin{align}\label{cond-tail-exc}
            &p_I(v)e^{-\frac{b}{a_1}v}\text{ is bounded on $(-\infty,V_R]$},\\
                &b_0-\frac{a_0}{a_1}b=V_L+\mu_0-\frac{a_0}{a_1}b
                \leq V_F.\label{cond-drift-exc}
        \end{align}
        \item When $b\leq 0$, the inhibitory case, we additionally assume
        \begin{equation}\label{cond-drift-inh}
            b_0=V_L+\mu_0\leq V_F.
        \end{equation}
    \end{enumerate}
\end{enumerate}
\end{theorem}
The proof of Theorem \ref{thm:gwp-t} is given in Section \ref{sec:4e} for the excitatory part and in Section \ref{sec:5i} for the inhibitory part.
\begin{remark} Assumptions on the drift \eqref{cond-drift-exc} and \eqref{cond-drift-inh} are not restrictive. A range of $\mu_0$ is allowed since $V_L\leq V_R<V_F$. Besides, we conjecture that they are technical rather than essential for the global well-posedness {(at least based on some preliminary numerics)}.
\end{remark}
Let's explain more on the non-existence of a global solution when $b\geq V_F-V_R$. By Proposition \ref{prop:gen-uni-c}, for a generalized solution, its maximal lifespan $T^*$ is given by \eqref{lifespan-2.1.1}
\begin{equation*}
    T^*=\int_0^{+\infty}\tilde{N}(\tau)d\tau=\int_0^{+\infty}\frac{1}{N+c}d\tau.
\end{equation*} 
In the strongly excitatory case, we can give examples with $T^*<+\infty$ which deny the global well-posedness. Such examples are due to a phenomenon which we refer to as the {\textbf{eternal blow-up}}. It means in the dilated timescale there exists some time $\tau_0\geq 0$ such that
\begin{equation}\label{eternal-blow-up-def}
    \tilde{N}(\tau)\equiv0,\quad\text{for all } \tau\geq \tau_0.
\end{equation} We call it an eternal blow-up since $\tilde{N}=0$ is equivalent to $N=+\infty$ \eqref{tildeN-0-N-inf}. Clearly, \eqref{eternal-blow-up-def} implies that in \eqref{lifespan-2.1.1}
\begin{equation*}
T^*= \int_0^{\tau_0}\tilde{N}(\tau)d\tau\leq \tau_0/c<+\infty,
\end{equation*} since $\tilde{N}\leq 1/c$ by \eqref{def-tildeN}. Thus then a global solution does not exist when there is an eternal blow-up. A numerical illustration is given in Figure \ref{fig:eternal-blow-up}.

The critical threshold $b=V_F-V_R$ is the same as that of the ``physical solution'' in \cite{delarue2015particle} (corresponds to $\alpha=1$ in their notations), which is a generalized solution from the SDE perspective in the regime $a_1=0$. They also give an intuition that when $b\geq V_F-V_R$, a neuron can spike infinite times in an instant. It is interesting to note that our eternal blow-up \eqref{eternal-blow-up-def} has a similar meaning, see Remark \ref{rmk:dtau=infty}.

In the next section we give the proof strategy for Theorem \ref{thm:gwp-t}. There, in particular in Proposition \ref{prop:steady-sec2}, we can see how to identify the threshold $b=V_F-V_R$.

\begin{remark}[When the initial data blows up]
In Theorem \ref{thm:gwp-t} we do not require that the initial data satisfies $-a_1\p_vp_I(V_F)<1$, therefore the firing rate $N(t)$ is allowed to blow up at $t=0$. As described in Proposition \ref{prop:regularity}, when $N(0)=+\infty$ the initial data $p_I$ shall be interpreted as $p_I(\cdot)=p(\cdot,0^-)$. Let's understand Theorem \ref{thm:gwp-tau} in this special case. 

In the dilated timescale, an initial blow-up of the firing rate means $\tilde{N}=0$ for $\tau=0$. Even so, we have the global well-posedness for generalized solutions when $b<V_F-V_R$. In view of the formula for $T^*$ \eqref{lifespan-2.1.1}, this implies that in the dilated timescale, $\tilde{N}(\tau)\not\equiv0$ (or equivalently ${N}\not\equiv+\infty$). In other words, the solution would get out of the blow-up regime eventually when $b<V_F-V_R$. 

But when $b\geq V_F-V_R$ there exist examples such that $\tilde{N}(\tau)\equiv0$ (or equivalently $N\equiv+\infty$) for all $\tau\geq 0$, which is the eternal blow-up defined in \eqref{eternal-blow-up-def} with $\tau_0=0$. It implies that the lifespan in the original time scale $T^*=0$ and hence denies the global well-posedness of the generalized solution in the strongly excitatory case. An investigation into the eternal blow-up is also a motivation towards the proof strategy for Theorem \ref{thm:gwp-t}.
\end{remark}

\subsubsection{Proof strategy: The limit equation at the blow-up}\label{sec:222ProofStrategy}

This section is devoted to the proof strategy of Theorem \ref{thm:gwp-t}, concerning whether a global generalized solution exists in the original timescale $t$.

Theorem \ref{thm:gwp-tau} ensures the global well-posedness of the equation in the dilated timescale $\tau$ \eqref{eq:NNLIF-tau}, which allows us to define the generalized solution. A generalized solution is global in the original timescale $t$ if and only if the lifespan $T^*$ given by \eqref{lifespan-2.1.1} is infinite, i.e.,
\begin{equation}\label{gwp-cond-tN}
    T^*=\int_0^{+\infty}\tilde{N}(\tau)d\tau=+\infty.
\end{equation} Recall the definition of $\tilde{N}$ \eqref{def-tildeN}
    \begin{equation*}
    \tilde{N}(\tau)=\frac{1}{{N}+c}=\frac{(1+a_1\p_v\tp(V_F,\tau))_+}{-a_0\p_v\tp(V_F,\tau)+c(1+a_1\p_v\tp(V_F,\tau))_+}.
    \end{equation*}
    {Condition \eqref{gwp-cond-tN} implies that the firing rate firing rate $N$, or the boundary flux $-a_1\p_v\tp(V_F,\tau)$ can not be always too large}. Hence, the global well-posedness in the original timescale $t$ is a question on the long time behavior of \eqref{eq:NNLIF-tau} in the dilated timescale $\tau$.

 First, we consider the extreme case when $\tilde{N}(\tau)\equiv0$ for all $\tau\geq0$, which simply gives $T^*=0<+\infty$ in \eqref{lifespan-2.1.1}. This is the eternal blow-up defined in \eqref{eternal-blow-up-def} with $\tau_0=0$. In this case, the dynamics of  \eqref{eq:NNLIF-tau} simplifies to 
 \begin{equation}\label{limit-linear}
	\begin{aligned}
	  \p_{\tau}\tp+b\p_v\tp&=a_1\p_{vv}\tp-a_1\p_v\tp(V_F,\tau)\delta_{v=V_R},\quad \tau>0,v\in(-\infty,V_F),\\
	    \tp(V_F,\tau)&=0,\ \tp(-\infty,\tau)=0,\ \tau>0.
	\end{aligned}
	\end{equation}
In fact, \eqref{limit-linear} is a linear equation, with constant coefficients and a flux jump mechanism. We call \eqref{eq:NNLIF-tau} \textbf{the limit equation at the blow-up}. Clearly \eqref{limit-linear} can be solved without the context of the nonlinear equation \eqref{eq:NNLIF-tau}. But a solution $\tp(v,\tau)$ of the linear equation \eqref{limit-linear} also solves the nonlinear equation \eqref{eq:NNLIF-tau} if and only if $\tilde{N}$ defined in \eqref{def-tildeN} is always zero. In view of the definition of $\tilde{N}$ \eqref{def-tildeN}, this is equivalent to the following condition on the boundary flux
	\begin{equation}\label{cond-linear}
	    -a_1\p_v\tp(V_F,\tau)\geq 1,\quad \forall \tau> 0.
	\end{equation}

To find a solution of \eqref{limit-linear} satisfying \eqref{cond-linear}, it is natural to look into the steady state. When $b>0$, i.e., the excitatory case, we have a unique probability density steady state with an explicit formula for its boundary flux as follows.
\begin{proposition}[Part of Proposition \ref{prop:steady-excit}]\label{prop:steady-sec2} When $b>0$, for the limit equation \eqref{limit-linear}, there exists a unique probability density steady state $\tp_{\infty}(v)$. Moreover, its boundary flux is given by
\begin{equation}\label{tmp-formula-flux}
    -a_1\p_v\tp_{\infty}(V_F)=\frac{b}{V_F-V_R}>0.
\end{equation}
\end{proposition} 

The critical threshold $b=V_F-V_R$ emerges in the formula \eqref{tmp-formula-flux}, which implies that $-a_1\p_v\tp_{\infty}(V_F)\geq1$ if and only if $b\geq V_F-V_R$. Therefore, when $b\geq V_F-V_R$, the steady state of \eqref{limit-linear}, which gives $\tilde{N}=0$, is also a steady state of \eqref{eq:NNLIF-tau}. If we set this steady state $\tp_{\infty}(v)$ as the initial data, we shall have $\tilde{N}(\tau)\equiv0$, which implies $T^*=0$ in \eqref{lifespan-2.1.1}. Hence, we prove the first part of Theorem \ref{thm:gwp-t} -- an example of non-existence in the strongly excitatory case. 

Now we describe the proof strategy for the global well-posedness when $b<V_F-V_R$. The key intuition also comes from the limit equation \eqref{limit-linear}. {We first consider the mildly excitatory case, i.e., the case when $0<b<V_F-V_R$.} 

Again we start with examining an extreme case -- the eternal blow-up. Precisely, we ask whether it is possible that $\tilde{N}=0$ holds for all $\tau\geq 0$ for a solution of \eqref{eq:NNLIF-tau} when $0<b<V_F-V_R$. That is equivalent to ask whether there exists a solution of the limit equation \eqref{limit-linear} with $-a_1\p_v\tp(V_F,\tau)\geq 1$ for all $\tau\geq 0$, a question concerning only the linear equation \eqref{limit-linear}. Intuitively the answer is no, because entropy methods \cite{caceres2011analysis} tell us that a solution of \eqref{limit-linear} converges in the long time to the steady state. But the boundary flux of the steady state $-a_1\p_v\tp_{\infty}(V_F)$ is strictly less than $1$ when $0<b<V_F-V_R$ by Proposition \ref{prop:steady-sec2}.

The long time behavior of \eqref{limit-linear} provides the cornerstone for analysis of general cases. We argue by contradiction, if on the contrary
\begin{equation}\label{contra-assump}
        \int_0^{+\infty}\tilde{N}(\tau)d\tau<+\infty,
\end{equation} then $\int_L^{+\infty}\tilde{N}(\tau)d\tau\rightarrow 0$ as $L\rightarrow+\infty$, which implies a smallness of $\tilde{N}(\tau)$. Hence, we can view the nonlinear equation \eqref{eq:NNLIF-tau} as a perturbation of the linear limit equation \eqref{limit-linear}, where the $\tilde{N}(\tau)$ terms are regarded as non-autonomous perturbations. By \eqref{contra-assump} the total influence of such perturbations shall be finite and therefore eventually the dynamics would be dominated by the linear equation \eqref{limit-linear} in the long time. If so, however, by the entropy estimate for the linear equation, the boundary flux of the solution shall be near to that of the steady state, which would contradict \eqref{contra-assump}.

To turn these intuitions into a rigorous proof, we adopt various techniques developed in literature for the NNLIF model \eqref{eq:NNLIF-original}, including the entropy estimate \cite{caceres2011analysis}, the control of the boundary flux \cite{CACERES201481,Antonio_Carrillo_2015}, and the $L^{\infty}$ bounds by constructing super solutions \cite{Antonio_Carrillo_2015}. {The proof for the excitatory case is given in Section \ref{sec:4e}.}

For the inhibitory case, i.e., $b\leq 0$, the general strategy is the same. Certain adjustments are needed, {which are mainly} due to a different behavior of the limit equation \eqref{limit-linear}. {The proof for the inhibitory case is given in Section \ref{sec:5i}.}

\subsection{Characterization of the blow-up scenario}\label{sec:2.3blow}
This section is devoted to characterizing the scenario when the firing rate blows up. We start by revisiting the argument for \eqref{jump-sec211} in a rigorous way as now we have the precise definition of the generalized solution. 

Let $(p,N)$ be a generalized solution with its lifespan $T^*>0$, and $(\tp,\tilde{N})$ be its time-dilation transform as in Definition \ref{def:generaized}. Suppose the firing rate blows up at $t=t_*<T^*$, i.e.,
\begin{equation*}
    N(t_*)=+\infty,\quad 
\end{equation*} Consider $[\tau_1,\tau_2]=\{\tau\geq 0,t(\tau)=t_*\}$, the preimage of $t_*$ under the map $t(\tau)$ in \eqref{change-back-map}. We have $0\leq \tau_1\leq \tau_2<+\infty$, where the last inequality is due to $t_*<T^*$. Then $t(\tau_1)=t(\tau_2)$ implies $\tilde{N}(\tau)=0$ for $\tau\in(\tau_1,\tau_2)$. If $\tau_1<\tau_2$, by continuity we have
\begin{equation}\label{2.3zero}
    \tilde{N}(\tau)=0,\quad \tau\in[\tau_1,\tau_2].
\end{equation}Actually \eqref{2.3zero} also holds in the subtle case $\tau_2=\tau_1$, since $N(t_*)=+\infty$ corresponds to $\tilde{N}(\tau_2)=0$ by definition of the generalized solution.

When $\tau_1<\tau_2$, the blow-up time $t_*$ is dilated to an interval $[\tau_1,\tau_2]$. Thus in the original timescale $p$ is discontinuous in time at $t_*$ since by Definition \ref{def:generaized}
\begin{equation*}
    p(\cdot,t_*^+)=p(\cdot,t_*)=\tp(\cdot,\tau_2)\neq \tp(\cdot,\tau_1)=p(\cdot,t_*^-).
\end{equation*}
We can further characterize the relationship between $p(\cdot,t_*^+)$ and $p(\cdot,t_*^-)$. Thanks to \eqref{2.3zero}, on ``the blow-up interval'' $[\tau_1,\tau_2]$ the dynamics \eqref{eq:NNLIF-tau} is reduced to the linear limit equation \eqref{limit-linear}. This observation leads to {the following proposition.}

\begin{proposition}[Characterization of the jump at a blow-up time]\label{prop:jump}
Let $(p,N)$ be a generalized solution of \eqref{eq:NNLIF-original} on $[0,T^*)$, where $T^*$ is its maximal lifespan. For every $t_*\in[0,{T}^*)$ with $N(t_*)=+\infty$, we have
\begin{equation}\label{semigroup-jump}
    p(\cdot,t_*^+)=p(\cdot,t_*)=S_{\Delta \tau}p(\cdot,t_*^{-}),
\end{equation}where $S_{\tau}(\tau\geq 0)$ is the time-evolution semigroup associated with the linear limit equation \eqref{limit-linear}, and $\Delta\tau\in[0,+\infty)$ can be characterized by \eqref{delta-tau} below.

In other words, suppose we solve the linear limit equation \eqref{limit-linear} with the initial data $p(\cdot,t_*^-)$ and denote the solution as $w(v,\tau)$, i.e.,
\begin{equation}
    \begin{aligned}
    	    \p_{\tau}w+b\p_vw&=a_1\p_{vv}w-a_1\p_vw(V_F,\tau)\delta_{v=V_R},\quad \tau>0,v\in(-\infty,V_F),\\
    	    w(V_F,\tau)&=0,\quad \tau>0,\\
    	    w(v,\tau=0)&=p(v,t_*^{-}),\quad v\in(-\infty,V_F].
    \end{aligned}
\end{equation}Then we have 
\begin{equation}
  p(v,t_*^+)=p(v,t_*)=w(v,\Delta\tau),\quad v\in(-\infty,V_F].
\end{equation}

Here $\Delta {\tau}\in[0,+\infty)$ is given by
\begin{align}\label{delta-tau}
    \Delta \tau:=& \inf_{\tau\geq 0}\{\tau\geq0:-a_1\p_vw(V_F,\tau)< 1\}\\ =&\sup_{\tau\geq 0}\{\tau\geq0:-a_1\p_vw(V_F,s)\geq 1,\quad \forall s\in[0,\tau]\}\label{delta-tau-2}.
\end{align}
\end{proposition}
\begin{proof}[Proof of Proposition \ref{prop:jump}]
Following the discussion above the statement of this proposition, it remains to show that $\tau_2-\tau_1$ coincides with $\Delta\tau$ defined in \eqref{delta-tau}, or equivalently \eqref{delta-tau-2}. Note that by \eqref{def-tildeN}, the boundary flux condition in \eqref{delta-tau-2} is equivalent to that the corresponding $\tilde{N}=0$. Hence, by \eqref{2.3zero} we have $\tau_2-\tau_1\leq \Delta\tau$. 

Suppose on the contrary $\Delta\tau>\tau_2-\tau_1$. We can construct another solution of \eqref{eq:NNLIF-tau} on $[\tau_1,\tau_1+\Delta \tau]$ with $\tilde{N}\equiv0$, by copying $w$ from the time interval $[0,\Delta\tau]$. Note that the definition of $\tau_2$ ensures that $\tilde{N}\not\equiv0$ on any (right) neighborhood of its. Hence, we have two different solutions of \eqref{eq:NNLIF-tau} with a same initial data at $\tau=\tau_1$. This contradicts the uniqueness of \eqref{eq:NNLIF-tau}, which is ensured by Theorem \ref{thm:gwp-tau}. 

Therefore $0\leq \Delta \tau=\tau_2-\tau_1$, in particular, $\Delta\tau<+\infty$. For the last inequality, we recall $\tau_2<+\infty$ due to $t_*< T^*$, which implies $\int_{\tau_1}^{+\infty}\tilde{N}(w)dw>0$. Therefore $\tilde{N}$ can not vanish for all $\tau\geq\tau_1$.
\end{proof}

\begin{figure}[htbp]
    \centering
    \includegraphics[width=1.0\linewidth]{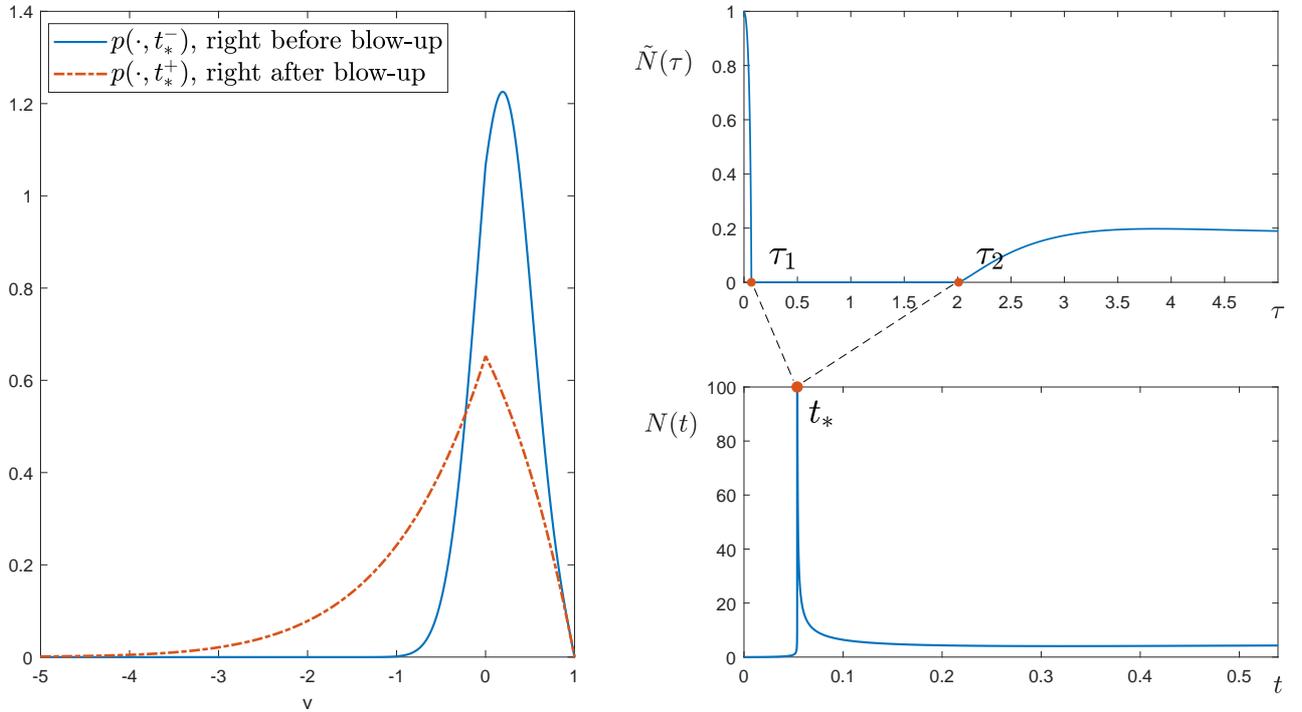}
    \caption{A numerical example of the generalized solution to illustrate what happens when the firing rate blows up. {Right top: $\tilde{N}(\tau)$ in the dilated timescale.} $\tilde{N}\equiv0$ on $[\tau_1,\tau_2]$. Right bottom: Firing rate $N(t)$ in the original timescale, truncated at $N=100$. The firing rate blows up at $t=t_*$, which corresponds to the interval $[\tau_1,\tau_2]$ in the dilated timescale. Left: Profiles of $p(\cdot,t_*^-)$, the solution just before the blow-up, and $p(\cdot,t_*^+)$, the solution just after the blow-up. Parameters: $V_F=1,V_R=0,b=0.9$ with $a_0=0.5,a_1=1,b_0=0$. The initial data is localized around $v=0.2$. The numerical scheme is a variant version from the one in \cite{hu2021structure}. }
    \label{fig:23-jump}
\end{figure}

\begin{remark}
Proposition \ref{prop:jump} (trivially) extends to the regular case when $N(t_*)<+\infty$ with $\Delta\tau=0$. Then $S_{\Delta\tau}=S_0=id$ is the identity operator, and \eqref{semigroup-jump} says that the density function is continuous in time when $N<+\infty$. 
\end{remark}
\begin{remark}\label{rmk:dtau=infty}
Formally, the eternal blow-up \eqref{eternal-blow-up-def} when $b\geq V_F-V_R$ can be viewed as $\Delta\tau=+\infty$ (or even $\Delta\tau \geq +\infty$). In this case a solution can never return to the non-blowup regime ($N<+\infty$) and can not continue to evolve in the original timescale. The neurons keep actively spiking with $N=+\infty$ in the dilated timescale, but this happens within an instant in the original timescale. 
\end{remark}

To summarize, when the firing rate blows up, the density function evolves according to the limit equation \eqref{limit-linear} for a time $\Delta \tau$ (in the dilated timescale), which is by \eqref{delta-tau} the {minimal time} for the boundary flux to become below the critical value, or equivalently to get out of the blow-up regime. All this happens in an infinitesimal time in the original timescale but typically $\Delta\tau>0$, which leads to a jump of the density function. A numerical illustration is given in Figure \ref{fig:23-jump}.

The following two remarks discuss the dynamics of the limit equation \eqref{limit-linear}.

\begin{remark}We identify two groups of terms in equation \eqref{eq:NNLIF-original-rw}. The first group is
\begin{equation}\label{first group}
    -\p_v[(-v+b_0)p]+a_0\p_{vv}p-a_0\p_vp(V_F,t)\delta_{v=V_R},
\end{equation}which does not depend on $N(t)$. The second group is, on the other hand, the terms containing a factor $N(t)$
\begin{equation}\label{second group}
    \bigr(-b\p_vp+a_1\p_{vv}p-a_1\p_vp(V_F,t)\delta_{v=V_R}\bigl)N(t).
\end{equation} 
From a modeling point of view, \eqref{first group} describes the dynamics of a single and independent neuron, while \eqref{second group} is the contribution of the interaction within a neuronal network, which is proportion to $N(t)$. If the firing rate $N(t)$ is very small, then \eqref{second group} is small and \eqref{first group} dominates. While in the (nearly) blow-up regime when $N(t)$ is very large, \eqref{second group} {is also} very large, which motivates us to study the dynamics in the dilated timescale $\tau$.

In the limit equation \eqref{limit-linear}, \eqref{first group} disappears since it is divided by $N=+\infty$. In particular, the parameters $b_0,a_0$ for the single-neuron dynamics -- the leaky mechanism and the external input, does not influence \eqref{limit-linear}. What matters then is the network connectivity parameter $b\in\mathbb{R}$ and {the scaling parameter of the intrinsic noise} $a_1>0$.
\end{remark}
\begin{remark}We might view \eqref{limit-linear} as a drift-diffusion equation with a unique fire-reset mechanism. The direction of the drift is determined by $b\in\mathbb{R}$: rightwards for $b>0$ (the excitatory case) and leftwards for $b<0$ (the inhibitory case). And we see again $a_1>0$ is crucial to make \eqref{limit-linear} a well-posed parabolic problem. 
\end{remark}

\section{Global well-posedness in the dilated timescale $\tau$}\label{sc:tau-gwp} 
In this section, we prove Theorem \ref{thm:gwp-tau}, the global well-posedness for classical solutions of the equation in the dilated timescale $\tau$ \eqref{eq:NNLIF-tau}. For notation simplicity, we drop the tilde in $\tp$ since in this section we work exclusively in the dilated timescale. First we recall \eqref{eq:NNLIF-tau}
\begin{equation*}
\begin{aligned}
        \p_{\tau}p+\p_v[(-v\tilde{N}+b_c\tilde{N}+b)p]=(a_c\tilde{N}&+a_1)\p_{vv}p\\& +(a_c\tilde{N}+a_1)(-\p_vp(V_F,\tau))\delta_{v=V_R},\quad v\in(-\infty,V_F),\tau>0,
        \end{aligned}
\end{equation*}with $p(V_F,\tau)=0,\tau\geq 0$ where $\tilde{N}$ is given by \eqref{def-tildeN}
\begin{equation*}
    \tilde{N}(\tau):=\frac{1}{{N}+c}=\frac{(1+a_1\p_vp(V_F,\tau))_+}{-a_0\p_vp(V_F,\tau)+c(1+a_1\p_vp(V_F,\tau))_+}\in [0,1/c].
\end{equation*} Here $c>0$ is the parameter to avoid the degeneracy from $N=0$, and $b_c,a_c$ are constants defined in \eqref{def-bcac}.

Classical solutions for the equation in the original timescale $t$ \eqref{eq:NNLIF-original} have been analysed in \cite{carrillo2013CPDEclassical}, where they have established the local well-posedness for any $b\in\mathbb{R}$ and the global well-posedness for $b\leq 0$, for \eqref{eq:NNLIF-original} with $a_1=0$. Our treatment for \eqref{eq:NNLIF-tau} is an adaption of \cite{carrillo2013CPDEclassical}, since \eqref{eq:NNLIF-tau}, which is formally derived from \eqref{eq:NNLIF-original}, shares similar structures with it. However, there are two differences. We shall prove the global well-posedness for classical solutions of \eqref{eq:NNLIF-tau}, regardless of the connectivity parameter $b\in\mathbb{R}$, i.e., for both the excitatory and the inhibitory case. Moreover, we work in the scenario $a_1>0$ rather than $a_1=0$. 

An important feature of \eqref{eq:NNLIF-tau}, the equation in the dilated timescale $\tau$, is that the source of nonlinearity $\tilde{N}\in[0,1/c]$ is always bounded as in \eqref{up-low-bd-tilde}. This is in contrast to $N(t)$ in \eqref{eq:NNLIF-original} for the original timescale $t$, which is allowed to take any value among $[0,+\infty)$ in principle. Our proof shall show that the boundedness of $\tilde{N}$ indeed plays a key role for the existence of a global solution.

Besides, when $a_1>0$, together with $a_0>0$, the equation possesses uniform parabolicity in the dilated timescale $\tau$, in the sense that the diffusion coefficient in \eqref{eq:NNLIF-tau} satisfies
\begin{equation*}
    \max(a_0/c,a_1)\geq a_c\tilde{N}+a_1=a_0\tilde{N}+a_1(1-c\tilde{N})\geq \min(a_0/c,a_1)>0,
\end{equation*} thanks to the definition of $a_c$ \eqref{def-bcac} and that $\tilde{N}\in[0,1/c]$ by \eqref{def-tildeN}. Hence, by providing enough diffusion, $a_1>0$ gives crucial benefits in the dilated timescale $\tau$, while it is an obstacle for the well-posedness of classical solutions in the original timescale $t$ \cite[Appendix A]{Antonio_Carrillo_2015}. 

Although the global well-posedness of \eqref{eq:NNLIF-tau} is the foundation for our generalized solution, its proof is independent of the reasoning in Section 4 and 5, and thus may be skipped at a first reading.

The proof strategy of Theorem \ref{thm:gwp-tau} is to first transform \eqref{eq:NNLIF-tau} to a Stefan-like free boundary problem, which is given in Section 3.1. Then we can reduce the PDE problem to an integral equation for $M(s)$, the boundary flux, which is an analogy of the firing rate in the free boundary formulation. In this way, with further analysis, the global well-posedness for the equivalent free boundary problem is established in Section 3.2. {This proof strategy closely follows \cite{carrillo2013CPDEclassical}.}

\subsection{Transform to a Stefan free boundary problem}
 By a translation, without loss of generality we set $V_F=0$ and therefore $V_R<0$. We introduce the following notations for the drift and diffusion coefficients, recalling definitions of $b_c,a_c$ \eqref{def-bcac} 
\begin{align}\label{def-tilde-mu}
    \tilde{\mu}(\tau)&:=b+b_c\tilde{N}(\tau)=\tilde{N}b_0+(1-c\tilde{N})b,\\\tilde{a}(\tau)&:=a_c\tilde{N}(\tau)+a_1=\tilde{N}a_0+(1-c\tilde{N})a_1.\label{def-tilde-a}
\end{align}Then we rewrite \eqref{eq:NNLIF-tau} as
\begin{equation}\label{NNLIF-3.1}
\p_\tau p+\p_v((-v\tilde{N}+\tilde{\mu})p)=\tilde{a}\p_{vv}p-\tilde{a}\p_vp(0,\tau)\delta_{v=V_R},\quad v<0,\, \tau>0.
\end{equation}

In this section, we transform \eqref{NNLIF-3.1} into an equivalent free boundary problem, which resembles the classical Stefan problem. Then in the next section, we do a further reduction to an integral equation of the boundary flux and prove the global well-posedness.

\subsubsection{Changes of variables}
First, we utilize a transform which is sightly modified from a classical change of variable \cite{carrillo2000asymptotic}. In our case, the reformulation involves $\tilde{N}(\tau)$, due to the dependence of the coefficients in \eqref{NNLIF-3.1} on it. 

Introduce a new space variable $y=\beta(\tau)v$ and a new time variable $s=S(\tau)$, where
\begin{align}\label{def-beta}
    \beta(\tau)&:=\exp\left(\int_0^{\tau}\tilde{N}({\tau}')d{\tau}'\right),\\
    \label{def-S}
    S(\tau)&:=\int_0^{\tau}\beta^2(\tau')\tilde{a}(\tau')d\tau'.
\end{align} Denote the solution in the new variables as $q(y,s)$, given by
\begin{equation}\label{change-p-q}
    \beta(\tau)q(\beta(\tau)v,S(\tau))=p(v,\tau).
\end{equation}
Then one derives from \eqref{NNLIF-3.1} via lengthy but direct calculations, whose details are presented in Appendix \ref{app:pf-sec3-cal},
\begin{equation}\label{reduce-q}
    \p_sq=\p_{yy}q-D(s)\p_yq+M(s)\delta_{y=V_R\beta(T(s))},\quad y<0,\,s>0,
\end{equation}where $M(s)$ is the boundary flux
\begin{equation}\label{3.1-M}
    M(s):=-\p_yq(0,s),
\end{equation} and $T(s)$ is the inverse map of $S(\tau)$ \eqref{def-S}, satisfying
\begin{equation}\label{def-T'}
    T'(s)=\frac{1}{\beta^2(T(s))\tilde{a}(T(s))},\quad T(0)=0.
\end{equation}
Finally, the drift $D(s)$ in \eqref{reduce-q} is given by
\begin{equation}\label{3.1-D}
    D(s):=\frac{1}{\beta(T(s))\tilde{a}(T(s))}\tilde{\mu}(T(s)).
\end{equation}

To justify that \eqref{reduce-q}--\eqref{3.1-D} indeed form a closed system, we need to represent $\beta(T(s))$ and $D(s)$ in terms of $M(s)$, using only the information in the new variables. This additional work comes from the fact that our change of variables \eqref{change-p-q} depends on the solution through $\tilde{N}$. From \eqref{def-tildeN} we represent $\tilde{N}$ in the new variables
\begin{equation}\label{sec311-tildeN}
\tilde{N}(T(s))=\frac{\bigl(1-a_1\beta^2(T(s))M(s)\bigr)_+}{a_0\bigl(\beta^2(T(s))M(s)\bigr)+c\bigl(1-a_1\beta^2(T(s))M(s)\bigr)_+}\in[0,1/c],
\end{equation}

Actually, we only need to solve for $\beta(T(s))$, thanks to that $\tilde{\mu}(T(s))$ and $\tilde{a}(T(s))$ are determined by $\tilde{N}(T(s))$ in \eqref{def-tilde-mu}-\eqref{def-tilde-a}, and that $\tilde{N}(T(s))$ is a function of $\beta(T(s))$ in \eqref{sec311-tildeN} when $M(s)$ is given. Applying the inverse change of variable $\tau=T(s)$ \eqref{def-T'} to \eqref{def-beta} we get
\begin{align}\notag
     \beta(T(s))&=\exp(\int_0^{s}\tilde{N}(T(\bar{s}))T'(\bar{s})d\bar{s})\\
     &=\exp(\int_0^{s}\tilde{N}(T(\bar{s}))\frac{1}{\beta^2(T(\bar{s}))\tilde{a}(T(\bar{s}))}d\bar{s})\label{sec311-beta-0}
\end{align} 
This integral equation \eqref{sec311-beta-0} indeed fully characterizes $\beta(T(s))$ under a given $M$, which allows us to denote it as $\beta(s;M)$. 

With $\beta(T(s))$ represented by $\beta(s;M)$, we can also represent $D(s)$ as $D(s;M)$  and rewrite \eqref{reduce-q} as a closed system
\begin{equation}\label{reduce-q-close}
    \p_sq=\p_{yy}q-D(s;M)\p_yq+M(s)\delta_{y=V_R\beta(s;M)},\quad y<0,\,s>0,
\end{equation} where $M(s)=-\p_yq(0,s)$. The complete justifications, including the precise definitions of $\beta(s;M)$ and $D(s;M)$, are postponed to Section \ref{sec-Dbeta-dependence}. By a further change of variable along the characteristics of \eqref{reduce-q-close}, we get the free boundary problem in Section \ref{sec:31fb}.

\subsubsection{The free boundary problem}\label{sec:31fb}
In view of the characteristics of \eqref{reduce-q-close}, we consider a further change of variable $u(x,s)=q(y,s)$, where
\begin{equation}\label{final-change-of-variable}
    x=y-\int_0^sD(\bar{s};M)d\bar{s}.
\end{equation} 
Then we get the following free boundary problem, which is equivalent to \eqref{NNLIF-3.1}.

\begin{lemma}[Corresponds to {\cite[Lemma 2.1]{carrillo2013CPDEclassical}}]\label{cpde-2.1}The system in the dilated timescale $\tau$ \eqref{eq:NNLIF-tau}--\eqref{def-tildeN} is equivalent to the free boundary problem
\begin{spacing}{1.2}
\begin{equation}\label{eq:freeboundary}
    \begin{cases}
    \p_su=\p_{xx}u+M(s)\delta_{x=\ell_R(s)},\quad x<\ell(s),\,s>0,\\
    \ell(s)=\ell_I-\int_0^sD(\bar{s};M)d\bar{s},\quad s>0,\\
    \ell_R(s)=\ell(s)+V_R\beta(s;M),\quad s>0,\\
    M(s)=-\p_xu|_{x=\ell(s)},\quad s>0,\\
    u(-\infty,s)=0,\ u(\ell(s),s)=0,\quad s>0,\\
    u(x,0)=u_I(x),\quad x<\ell_I.
    \end{cases}
\end{equation}
\end{spacing}
Here $\beta(s;M)$ and $D(s;M)$ are given by Lemma \ref{lemma:existence} in Section \ref{sec-Dbeta-dependence}. 
\end{lemma}

Hence, to prove the global well-posedness of \eqref{eq:NNLIF-tau}, it is equivalent to study \eqref{eq:freeboundary}. We shall establish the global well-posedness for classical solutions of \eqref{eq:freeboundary}. The definition of a classical solution of \eqref{eq:freeboundary} and the assumption on its initial data can be directly deduced from Definition \ref{def:classical-tau} and Assumption \ref{as:classical-init} via the change of variable discussed above. Thus they are omitted here.

The nonlinearity of \eqref{eq:freeboundary} comes from the fact that $\ell(s)$ and $\ell_R(s)$ depend on the solution through its boundary flux $M(s)$. Compared to the free boundary problem in \cite{carrillo2013CPDEclassical}, the only but important difference is how $\ell(s)$ and $\ell_R(s)$ move. Here we have a uniform a priori control on the moving speeds of the free boundaries $\ell(s),\ell_R(s)$, which is crucial for getting a global solution.

\begin{lemma}\label{lemma:Lip-s}
For a classical solution of the free boundary problem \eqref{eq:freeboundary} on a time interval $[0,S)$, we have the following uniform controls on the moving speeds of $\ell(s),\ell_R(s)$
\begin{equation*}
    |\ell'(s)|\leq L,\quad |\ell_R'(s)|\leq L,\quad s\in[0,S),
\end{equation*} which implies that $\ell(s)$ and $\ell_R(s)$ are Lipschitz continuous with a global Lipschitz constant
\begin{equation*}
    |\ell(s_1)-\ell(s_2)|+|\ell_R(s_1)-\ell_R(s_2)|\leq 2L|s_1-s_2|,\quad\quad \forall s_1,s_2\in [0,S).
\end{equation*} Here $L>0$ is a (universal) constant independent of the solution, or the time $S$.
\end{lemma}

Such a uniform bound, which is the key for a global well-poseness, comes from the fact that $\tilde{N}(\tau)\in[0,1/c]$ by \eqref{def-tildeN} in the Fokker-Planck formulation \eqref{eq:NNLIF-tau}. The rigorous proof of Lemma \ref{lemma:Lip-s} needs the precise definitions of $\beta(s;M)$ and $D(s;M)$, and is postponed to the end of Section \ref{sec-Dbeta-dependence}.

\subsubsection{Definitions and properties of $\beta(s;M)$ and $D(s;M)$}\label{sec-Dbeta-dependence} 

In this section, we give and justify the definitions of $\beta(s;M)$ and $D(s;M)$, therefore completing the transform from \eqref{NNLIF-3.1} to \eqref{reduce-q-close} and the free boundary problem \eqref{eq:freeboundary}. Besides, properties of $\beta(s;M)$ and $D(s;M)$ are given as preparations for the well-posedness proof in the next section.

Given a non-negative function $M(s)\in C[0,S_0]$ with some $S_0>0$, we want to find $\beta(s;M)\in C[0,S_0]$ which represents $\beta(T(s))$. In view of \eqref{sec311-beta-0}, $\beta(s;M)$ shall be characterized by  
\begin{equation}\label{312-beta}
    \beta(s;M)=\exp\left(\int_0^{s}\frac{1}{\beta^2(\bar{s};M)}\frac{\tilde{N}(\bar{s};M)}{\tilde{a}(\bar{s};M)}d\bar{s}\right),
\end{equation} where $\tilde{N}(s;M)$ is given by a similar equation to \eqref{sec311-tildeN}
\begin{equation}
    \tilde{N}(s;M)=\frac{\bigl(1-a_1\beta^2(s;M)M(s)\bigr)_+}{a_0\bigl(\beta^2(s;M)M(s)\bigr)+c\bigl(1-a_1\beta^2(s;M)M(s)\bigr)_+},\label{312-N}
\end{equation}  and $\tilde{a}(s;M)$ is defined as in \eqref{def-tilde-a}
\begin{equation}
    \tilde{a}(s;M):=\tilde{N}(s;M)a_0+(1-c\tilde{N}(s;M))a_1.\label{312-a-s-M}
\end{equation}

Taking the derivative with respect to $s$ in \eqref{312-beta}, we find 
\begin{align}\notag
    \beta'(s;M)&=\frac{1}{\beta^2({s;M})}\frac{\tilde{N}({s};M)}{\tilde{a}(s;M)}\exp\left(\int_0^{s}\frac{1}{\beta^2(\bar{s};M)}\frac{\tilde{N}(\bar{s};M)}{\tilde{a}(\bar{s};M)}d\bar{s}\right)\\&=\frac{1}{\beta({s;M})}\frac{\tilde{N}({s};M)}{\tilde{a}(s;M)}.\label{ode-beta-312}
\end{align} Thanks to \eqref{312-a-s-M} $\tilde{a}(s;M)$ is determined by $\tilde{N}(s;M)$ and thanks to \eqref{312-N} $\tilde{N}(s;M)$ is determined by $\beta(s;M)$ and $M(s)$. Therefore \eqref{ode-beta-312} gives an ODE about $\beta$, which can be solved as long as we know $M(s)$. Solving \eqref{ode-beta-312} we can obtain $\beta$ and subsequently $\tilde{N}$ and $\tilde{a}$ from \eqref{312-N} and \eqref{312-a-s-M}. Then we define $D(s;M)$ as in \eqref{3.1-D}
\begin{equation}\label{312-D-s-M}
        D(s;M)=\frac{1}{\beta(s;M)\tilde{a}(s;M)}\tilde{\mu}(s;M),
\end{equation} where, corresponding to \eqref{def-tilde-mu}, we denote
\begin{align}\label{312-mu-s-M}
        \tilde{\mu}(s;M)&:=\tilde{N}(s;M)b_0+(1-c\tilde{N}(s;M))b.
\end{align}

In this way we can show that $\beta(s;M)$ and $D(s;M)$ are well-defined, formulated as the following lemma.
\begin{lemma}\label{lemma:existence}
For a non-negative function $M(s)\in C[0,S_0]$ with some $S_0>0$, there exists unique $\beta(s;M),D(s;M)\in C[0,S_0]$ such that \eqref{312-beta}-\eqref{312-a-s-M} and \eqref{312-D-s-M}-\eqref{312-mu-s-M} hold.
\end{lemma} The proof of Lemma \ref{lemma:existence} just involves some calculations to check the well-posedness of an ODE. And we postpone it to Appendix \ref{app:pf-sec3}.

Lemma \ref{lemma:existence} completes the descriptions of \eqref{reduce-q-close} and \eqref{eq:freeboundary}. In the light of the inverse map \eqref{def-T'}, we define $$T(s)=\int_0^s\frac{1}{\beta^2(s';M)\tilde{a}(s';M)}ds'.$$ Then we can directly verify that the change of variable $$\frac{1}{\beta(s;M)}p\left(\frac{y}{\beta(s;M)},T(s)\right)=q(y,s)$$ sends a solution of \eqref{reduce-q-close} back to a solution of \eqref{NNLIF-3.1}.

For the well-posedness of \eqref{eq:freeboundary}, more knowledge on $D(s;M),\beta(s;M)$ is needed, summarized as the following lemma.

\begin{lemma}\label{lemma:key-bound}
(i) For a non-negative function $M(s)\in C[0,S_0]$ with some $S_0>0$, let $\beta(s;M)$ and $D(s;M)$ be as in Lemma \ref{lemma:existence}. Then we have the following bounds for $s\in[0,S_0]$
\begin{equation}\label{bound-beta}
    1\leq \beta(s;M)\leq \sqrt{Cs+1},
\end{equation} 
\begin{equation}\label{bound-D-key}
    |D(s;M)|\leq C.
\end{equation}
where $C>0$ is a constant independent of $S_0$ or $M(s)$. 

(ii) For two non-negative functions $M_1,M_2\in C[0,S_0]$ with some $S_0\leq S<+\infty$, we have
\begin{equation}\label{D-pertub}
    \|D(\cdot;M_1)-D(\cdot;M_2)\|_{C[0,S_0]}\leq C_S\|M_1-M_2\|_{C[0,S_0]}.
\end{equation}
Moreover, if for some $0\leq\sigma \leq S_0$, $M_1(t)\equiv M_2(t)$ for all $t\in[0,S_0-\sigma)$, i.e., $M_1$ and $M_2$ only start to differ from the time $S_0-\sigma$, then
\begin{equation}\label{beta-pertub}
    \|\beta(\cdot;M_1)-\beta(\cdot;M_2)\|_{C[S_0-\sigma,S_0]}\leq \sigma C_S\|M_1-M_2\|_{C[S_0-\sigma,S_0]}
\end{equation}
Here the constant $C_S$ depends only on $S$, an upper bound of $S_0$. In other words, the Lipschitz constant is uniformly bounded for a bounded family of $S_0$.
\end{lemma}
In Lemma \ref{lemma:key-bound}-(ii) we show that $\beta$ and $D$ depend on $M$ in a Lipschitz continuous way, which is important for the fixed point argument towards the local well-posedness. The time interval $[S_0-\sigma,S_0]$ is introduced because in the well-posedness proof we need to construct a solution starting at some time $s>0$ other than $s=0$.

The bound \eqref{bound-D-key} is the main ingredient for the uniform control Lemma \ref{lemma:Lip-s}, which is the key for the global well-posedness. The full proof of Lemma \ref{lemma:key-bound}, which just involves some elementary calculations and ODE analysis, is postponed to Appendix \ref{app:pf-sec3}. 

We finish this section with a self-contained proof of Lemma \ref{lemma:Lip-s}. 

\begin{proof}[Proof of Lemma \ref{lemma:Lip-s}]
We shall recover the bounds through the definitions of $\beta(s;M)$ and $\tilde{N}(s;M)$
\begin{align*}
    \beta(s;M)&\geq 1,\quad \tilde{N}(s;M)\in[0,1/c],\\
    \tilde{a}(s;M)&\geq \min(a_0/c,a_1)>0,\\
    |\tilde{\mu}(s;M)|&\leq \max(|b_0|/c,|b|).
\end{align*}
Indeed, by \eqref{312-beta} we get $\beta\geq 0$ therefore by \eqref{312-N} we have $\tilde{N}\in[0,1/c]$. Therefore by \eqref{312-a-s-M} we obtain $\tilde{a}(s;M)\geq \min(a_0/c,a_1)>0$. Then using \eqref{312-beta} again we get $\beta\geq 1$. Also by \eqref{312-mu-s-M} and $\tilde{N}\in[0,1/c]$ we have $|\tilde{\mu}(s;M)|\leq \max(|b_0|/c,|b|)$.

Hence, for $\ell(s)$ we get
\begin{align*}
    |\ell'(s)|=|D(s;M)|&=\frac{1}{\beta(s;M)\tilde{a}(t;M)}|\tilde{\mu}(t;M)|\\&\leq \frac{\max(|b_0|/c,|b|)}{\min(a_0/c,a_1)}<+\infty.
\end{align*}
For $\ell_R(s)$, it remains to deal with $\beta'(s;M)$ since $|\ell_R'(s)|\leq |\ell'(s)|+|\beta'(s;M)|$. By \eqref{ode-beta-312} we have
\begin{align*}
    |\beta'(s;M)|=|\frac{1}{\beta(s;M)}\frac{\tilde{N}(s;M)}{\tilde{a}(s;M)}|\leq \frac{1}{\min(a_0,a_1c)}<+\infty.
\end{align*}
Therefore we conclude the proof and the Lipschitz constant $L$ can be taken as $\frac{\max(|b_0|,|b|c)}{\min(a_0,a_1c)}+\frac{1}{\min(a_0,a_1c)}$.
\end{proof}

\subsection{Well-posedness of the equivalent free boundary problem}
In this section we establish the global well-posedness of the free boundary problem \eqref{eq:freeboundary}, which is equivalent to \eqref{eq:NNLIF-tau}, therefore proving Theorem \ref{thm:gwp-tau}. {This step is also a modification from the approach in \cite{carrillo2013CPDEclassical}} (see also classical treatments on the Stefan problem \cite{friedman1959free,MR0181836}). 

We first establish the local well-posedness of \eqref{eq:freeboundary}, with an existence time controlled by an upper bound of the derivative of the solution. Then we shall show that the derivative of a solution can not blow up in finite time, which leads to the global well-posedness. Here the uniform bounds on the boundary moving speeds $\ell'(s),\ell_R'(t)$, given in Lemma \ref{lemma:Lip-s}, play a key role for the proof.

\subsubsection{Local well-posedness}\label{sec:local}

This section is devoted to proving the following local well-posedness result for \eqref{eq:freeboundary}.
\begin{theorem}[Corresponds to {\cite[Theorem 3.1]{carrillo2013CPDEclassical}}]\label{cpde-3.1}The free boundary problem \eqref{eq:freeboundary} has a unique classical solution on some time interval $[0,\sigma]$, provided that the initial data $u_I(x)$ satisfies {the assumption corresponding to Assumption \ref{as:classical-init}}. Here $\sigma>0$ only depends on an upper bound of $$\sup_{-\infty<x\leq \ell_I}\left|\frac{\p u_I}{\p x}\right|.$$
Moreover, suppose we have a classical solution on the time interval $[0,s_0]$, then we could extend it uniquely to the time interval $[0,s_0+\sigma]$. Here $\sigma>0$ depends on an upper bound of $s_0$ and an upper bound of $$\sup_{-\infty<x\leq \ell(s_0)}\left|\frac{\p u(x,s_0)}{\p x}\right|.$$
\end{theorem}
The proof strategy is to reduce \eqref{eq:freeboundary} to an integral equation for the boundary flux $M(s)$ {on which a fixed point argument can be applied}. 

We start with deriving integral equations from \eqref{eq:freeboundary} via the heat kernel
\begin{equation}\label{heat-kernel}
G(x,s,\xi,\tau)=\frac{1}{[4\pi(s-\tau)]^{1/2}}\exp\left(-\frac{|x-\xi|^2}{4(s-\tau)}\right),\quad\quad s>\tau>0,\quad x,\xi\in\mathbb{R}.
\end{equation} The derivations are identical to those in \cite[Section 3.1]{carrillo2013CPDEclassical}, which do not depend on the specific forms of $\ell(s)$ and $\ell_R(s)$. First one can represent a solution $u(x,s)$ of \eqref{eq:freeboundary} as follows
\begin{equation}\label{u-3.4}
    \begin{aligned}
        u(x,s)=\int_{-\infty}^{\ell(0)}&G(x,s,\xi,0)u_I(\xi)d\xi-\int_0^sM(\tau)G(x,s,\ell(\tau),\tau)d\tau\\
        &+\int_0^sM(\tau)G(x,s,\ell_R(\tau),\tau)d\tau.
    \end{aligned}
\end{equation}
In \eqref{u-3.4}, the first term comes from the initial data, the second term reflects the loss due to Dirichlet zero boundary condition at $x=\ell(s)$, and the third term is from the Dirac source at $x=\ell_R(s)$. Then one can derive an integral equation for $M(s)=-\p_xu|_{x=\ell(s)}$
\begin{equation}\label{M-3.6}
    \begin{aligned}
        M(s)=-2\int_{-\infty}^{\ell(0)}&G(\ell(s),s,\xi,0)u_I'(\xi)d\xi+2\int_0^sM(\tau)G_x(\ell(s),s,\ell(\tau),\tau)d\tau\\&-2\int_0^sM(\tau)G_x(\ell(s),s,\ell_R(\tau),\tau)d\tau.
    \end{aligned}
\end{equation}

In this way, the free boundary problem \eqref{eq:freeboundary} is reduced to an integral equation for $M(s)$ \eqref{M-3.6}. This integral formulation provides a platform to perform a fixed point argument for local-wellposedness.

\begin{theorem}[Corresponds to {\cite[Theorem 3.2]{carrillo2013CPDEclassical}} ]\label{cpde-3.2}The integral equation \eqref{M-3.6} has a unique solution $M(s)\in C[0,\sigma]$, provided that $u_I(x)$ satisfies the assumption corresponding to Assumption \ref{as:classical-init}. Here $\sigma>0$ only depends on an upper bound of $$\sup_{-\infty<x\leq \ell_I}\left|\frac{\p u_I}{\p x}\right|.$$
Moreover, suppose we have a classical solution of \eqref{eq:freeboundary} on time interval $[0,t_0]$. Then the following equation
\begin{equation}\label{M-3.6-t0}
    \begin{aligned}
        M(s)=-2\int_{-\infty}^{\ell(s_0)}&G(\ell(s),s,\xi,t_0)u_x(\xi,s_0)d\xi+2\int_{s_0}^{s}M(\tau)G_x(\ell(s),s,\ell(\tau),\tau)d\tau\\&-2\int_{s_0}^{s}M(\tau)G_x(\ell(s),s,\ell_R(\tau),\tau)d\tau,
    \end{aligned}
\end{equation}
which is \eqref{M-3.6} with the starting time changed to $s_0>0$, has a unique solution $M(s)\in C[s_0,s_0+\sigma]$. Here $\sigma>0$ only depends on an upper bound of $s_0$ and an upper bound of $$\sup_{-\infty<x\leq \ell(s_0)}\left|\frac{\p u(x,s_0)}{\p x}\right|.$$
\end{theorem}
The proof is based on the Banach fixed point Theorem. We shall use the Lipschitz dependence of $D(s;M),\beta(s;M)$ on $M$ in Lemma \ref{lemma:key-bound} to give estimates. By choosing a quantitatively small enough $\sigma$, we gain smallness to construct a contraction map. The detailed proof is omitted here, since it follows the same lines as \cite[Proof of Theorem 3.2]{carrillo2013CPDEclassical}.

With a $M(s)\in C[0,\sigma]$ given by Theorem \ref{cpde-3.2}, we get (a mild solution) $u$ via \eqref{u-3.4}. One can show such an integral formulation \eqref{u-3.4} indeed gives a classical solution, by an identical augment as \cite[Corollary 3.3]{carrillo2013CPDEclassical}. Therefore, Theorem \ref{cpde-3.1} is proved.

\subsubsection{Global well-posedness}\label{sec:global}
In view of the local result Theorem \ref{cpde-3.1}, towards a global well-posedness we need to control $\sup_{x\leq \ell(s)}|u_x(x,s)|$. 

First, we show that a uniform bound for $\sup_{x\leq \ell(s)}|u_x(x,s)|$ can be implied by a {uniform estimate} on $M(s)$.
\begin{proposition}[Corresponds to{ \cite[Proposition 4.1]{carrillo2013CPDEclassical}}]\label{cpde-4.1}
For a classical solution of \eqref{eq:freeboundary} on some time interval $[0,s_0)$ with $0<s_0 <+\infty$. Suppose for some $0<\eps<s_0$, the following holds
\begin{equation}
    M^*=\sup_{s\in(s_0-\eps,s_0)}M(s)<\infty.
\end{equation}
Then we have
\begin{equation}
    \sup\{|u_x(x,s)|:x\in(-\infty,\ell(s)],s\in[s_0-\eps,s_0)\}<\infty,
\end{equation}and the bound depends only on $M^*,s_0$ and $U_0$. Here $U_0$ is defined by
\begin{equation*}
    U_0:=\sup_{x\in(-\infty,\ell(s_0-\eps)]}|u_x(x,s_0-\eps)|<\infty.
\end{equation*} $U_0<\infty$ is ensured by the definition of a classical solution. 
\end{proposition}

Such a uniform estimate on $u_x$ is obtained through the integral formula \eqref{u-3.4}. The proof is omitted here since it is identical to the proof of \cite[Proposition 4.1]{carrillo2013CPDEclassical}, thanks to the Lipschitz continuity for $\ell,\ell_R$ in Lemma \ref{lemma:Lip-s}. 

Then, one can extend the lifespan of a solution as long as $M(s)$ is bounded.
\begin{proposition}[Corresponds to {\cite[Theorem 4.2]{carrillo2013CPDEclassical}}]\label{cpde-4.2}
For a classical solution of \eqref{eq:freeboundary} on the time interval $[0,S)$ with $0<S<+\infty$, if 
\begin{equation}\label{3.2-bound-M}
\sup_{0\leq s<S}M(s)<+\infty,
\end{equation} then the solution can be extended to $[0,S+\sigma)$ for some $\sigma>0$.
\end{proposition}
Indeed, if \eqref{3.2-bound-M} holds by Proposition \ref{cpde-4.1} we get bounds on $\sup_{x\leq \ell(s)}|u_x(x,s)|$. Then applying Theorem \ref{cpde-3.1} for $s\in[0,S)$, we can extend the solution to $[s,s+\sigma]$ with a uniform $\sigma>0$. For a detailed proof, one can see \cite[Theorem 4.2]{carrillo2013CPDEclassical}

Then, the following proposition rules out the possibility that $M(s)$ blows up in a finite time, which allows the solution to be extended step by step towards infinity.
\begin{proposition}[Corresponds to {\cite[Proposition 4.3]{carrillo2013CPDEclassical}}]\label{cpde-4.3}
For any classical solution on a time interval $[0,s_0)$ with $0<s_0<+\infty$, we have
\begin{equation}\label{cpde-4.3-1}
    \sup_{0\leq s<s_0}M(s)<+\infty.
\end{equation}
Precisely, there exists an $\eps>0$ small enough such that for any classical solution one can control 
\begin{equation}\label{cpde-4.3-2}
    \sup_{s_0-\eps<s<s_0}M(s)\leq C_0<\infty,
\end{equation} with the bound $C_0$ only depends on
\begin{equation}
    \mathscr{M}_0:=\sup_{x\in(-\infty,\ell(s_0-\eps)]}|u_x(x,s_0-\eps)|<\infty.
\end{equation} Here $\mathscr{M}_0<+\infty$ is implied by the definition of a classical solution. And we assume $\eps<s_0$, otherwise we can choose a smaller $\eps$ which does not matter for the analysis of global well-posedness.
\end{proposition}
We postpone the proof of Proposition \ref{cpde-4.3} to Appendix \ref{app:pf-sec3}. In this case, the key to the proof is the uniform Lipschitz continuity of $\ell(s)$ (Lemma \ref{lemma:Lip-s}) due to the boundedness of $\tilde{N}$ \eqref{up-low-bd-tilde}.

Now we prove Theorem \ref{thm:gwp-tau}, the global well-posedness of classical solutions for the equation in the dilated timescale $\tau$ \eqref{eq:NNLIF-tau}.
\begin{proof}[Proof of Theorem \ref{thm:gwp-tau}]
By Lemma \ref{cpde-2.1}, it is equivalent to show the global well-posedness of the free boundary problem \eqref{eq:freeboundary}. Suppose the maximal time interval for the solution is $[0,S^*)$, by the local result Theorem \ref{cpde-3.1}, $S^*>0$. If $S^*<+\infty$, then by Proposition \ref{cpde-4.3}
\begin{equation*}
    \sup_{0\leq s<S^*}M(s)<+\infty.
\end{equation*} Hence by Proposition \ref{cpde-4.2} the solution can be extended to $[0,S^*+\sigma)$ with some $\sigma>0$, which leads to a contradiction against the definition of $S^*$. Hence, $S^*=+\infty$, which means that the solution is global.
\end{proof}

\section{Excitatory case: Global well-posedness in $t$ versus eternal blow-up}\label{sec:4e}

With global classical solutions in the dilated timescale $\tau$ \eqref{eq:NNLIF-tau} at hand, in the original timescale $t$ generalized solutions of \eqref{eq:NNLIF-original} are well-defined as in Definition \ref{def:generaized} and Proposition \ref{prop:gen-uni-c}. Moreover, the maximal lifespan of a generalized solution is exactly given by \eqref{lifespan-2.1.1} 
\begin{equation*}
    T^*=\int_0^{+\infty}\tilde{N}(\tau)d\tau,
\end{equation*} by Proposition \ref{prop:gen-uni-c}. Owing to the definition of $\tilde{N}$ \eqref{def-tildeN}
\begin{equation*}
    \tilde{N}(\tau)=\frac{1}{{N}+c}=\frac{(1+a_1\p_v\tp(V_F,\tau))_+}{-a_0\p_v\tp(V_F,\tau)+c(1+a_1\p_v\tp(V_F,\tau))_+},
\end{equation*}we see a longer lifespan is closely related with upper controls on the boundary flux $-a_1\p_v\tp(V_F,\tau)$. Therefore, the global well-posedness for generalized solutions of \eqref{eq:NNLIF-original}, which is equivalent to $T^*=+\infty$, is a problem concerning the long time behavior of \eqref{eq:NNLIF-tau} in the dilated timescale $\tau$. 

In this section we prove Theorem \ref{thm:gwp-t} for the excitatory case, i.e., the case $b>0$. We shall show $T^*=+\infty$ for $0<b<V_F-V_R$ and give examples with $T^*=0$ when $b\geq V_F-V_R$. The inhibitory case when $b\leq 0$ is treated in the next section. The general strategy, as sketched in Section \ref{sec:222ProofStrategy}, is the same.

First in Section \ref{sec:41} we study the limit equation \eqref{limit-linear}. In the strongly excitatory case $b\geq V_F-V_R$, the steady state of \eqref{limit-linear} provides an example of the eternal blow-up defined in \eqref{eternal-blow-up-def} therefore denying the global well-posedness in that case. Then in Section \ref{sec:42}, a contradiction argument allows us to treat \eqref{eq:NNLIF-tau} as a perturbation of \eqref{limit-linear}. And we prove the global well-posedness of the generalized solution of \eqref{eq:NNLIF-original} in the mildly excitatory case, i.e., the case $0<b<V_F-V_R$.

For notation convenience, we shall drop the tilde in $\tp$ in this section and the next section, as we will work in timescale $\tau$ only to prove $T^*=+\infty$.

\subsection{The limit equation and the eternal blow up}\label{sec:41}

We first study the limit equation \eqref{limit-linear}
\begin{equation*}
	\begin{aligned}
	  \p_{\tau}p+b\p_vp&=a_1\p_{vv}p-a_1\p_vp(V_F,\tau)\delta_{v=V_R},\quad \tau>0,v\in(-\infty,V_F),\\
	    p(V_F,\tau)&=0,\ p(-\infty,\tau)=0,\ \tau>0.
	\end{aligned}
\end{equation*} A solution of \eqref{limit-linear} is also a solution of \eqref{eq:NNLIF-tau} if and only if the corresponding $\tilde{N}\equiv0$, or equivalently
\begin{equation*}
    -a_1\p_vp(V_F,\tau)\geq 1,\quad \tau>0.
\end{equation*} This is the so-called an eternal blow-up as defined in \eqref{eternal-blow-up-def}.

We start with the steady state of \eqref{limit-linear}, which gives the counter-example against the global well-posedness in the case $b\geq V_F-V_R$. Then the long time behavior of \eqref{limit-linear} is analyzed via entropy estimates, which serves as the foundation for the global well-posedness proof for the case $0<b<V_F-V_R$.

	\subsubsection{Steady state}
	The steady state of \eqref{limit-linear} satisfies
	\begin{equation}\label{steady-limit}
	\begin{aligned}
	\p_v(bp)&=a_1\p_{vv}p-a_1\p_vp(V_F,\tau)\delta_{v=V_R},\quad v\in(-\infty,V_F),\\
	p(V_F)&=0,\ p(-\infty)=0.
	\end{aligned}
	\end{equation}
	Solving \eqref{steady-limit} for $b>0$, we obtain the following proposition.
	\begin{proposition}\label{prop:steady-excit}
	When $b>0$, the steady state of \eqref{limit-linear} is unique up a multiplying constant. Moreover, we can pick a constant such that the steady state is a probability density on $(-\infty,V_F)$, given by
		\begin{equation}\label{formula-steady-state}
	    p_{\infty}(v)=\frac{1}{Z}\begin{cases}
	    (e^{\frac{b}{a_1}(V_F-V_R)}-1)e^{\frac{b}{a_1}(v-V_R)},\quad v\leq V_R,\\
	    e^{\frac{b}{a_1}(V_F-V_R)}-e^{\frac{b}{a_1}(v-V_R)},\quad v\in(V_R,V_F),
	    \end{cases}
	\end{equation} where the normalization constant $Z=(V_F-V_R)e^{\frac{b}{a_1}(V_F-V_R)}$. In particular, the boundary flux at $v=V_F$ for \eqref{formula-steady-state} is given by
	\begin{equation}\label{def-Minfty}
	   M_{\infty}:=-a_1\p_v{p_{\infty}}(V_F)=\frac{b}{V_F-V_R}>0.
	\end{equation} Thus $M_{\infty}\geq 1$ when $b\geq V_F-V_R$ and $0<M_{\infty}<1$ when $0<b<V_F-V_R$.
	\end{proposition}
	\begin{proof} By a translation, we set $V_R=0$ WLOG to simplify calculations.
	On $(-\infty,0)\cup(0,V_F)$, the solution satisfies
	\begin{equation*}
	    \p_{vv}p-\frac{b}{a_1}\p_vp=0.
	\end{equation*}
	On $(-\infty,0)$ by the  decay at infinity, $p$ is given by $c_1e^{\frac{b}{a_1}v}$. On $(0,V_F)$ by the boundary condition at $V_F$, $p$ is given by $c_2(e^{\frac{b}{a_1}V_F}-e^{\frac{b}{a_1}v})$. At $v=0$, thanks to the flux jump condition
	\begin{equation*}
	    a_1\p_vp(0^+)-a_1\p_vp(0^-)=a_1\p_vp(V_F),
	\end{equation*} we get that $c_1=(e^{\frac{b}{a_1}V_F}-1)c_2$. Therefore the steady state is unique up to a multiplying constant. We can choose an appropriate $Z$ such that
	\begin{equation*}
	    p_{\infty}(v)=\frac{1}{Z}\begin{cases}
	    (e^{\frac{b}{a_1}V_F}-1)e^{\frac{b}{a_1}v},\quad v\leq 0,\\
	    e^{\frac{b}{a_1}V_F}-e^{\frac{b}{a_1}v},\quad v\in(0,V_F),
	    \end{cases}
	\end{equation*}is a probability density. Such a normalization constant $Z$ is given by
	\begin{equation*}
	    Z=\int_{-\infty}^0(e^{\frac{b}{a_1}V_F}-1)e^{\frac{b}{a_1}v}dv+\int_0^{V_F}(e^{\frac{b}{a_1}V_F}-e^{\frac{b}{a_1}v})dv=V_Fe^{\frac{b}{a_1}V_F}.
	\end{equation*}
Then, it is direct to compute the diffusion flux at $v=V_F$ 
	\begin{equation*}
	   M_{\infty}:=-a_1\p_v{p_{\infty}}(V_F)=\frac{1}{Z}be^{\frac{b}{a_1}V_F}=\frac{b}{V_F}.
	\end{equation*}
	\end{proof}

In the strongly excitatory case, i.e., $b\geq V_F-V_R$, by Proposition \ref{prop:steady-excit} the steady state of \eqref{limit-linear} satisfies $-a_1\p_v p_{\infty}(V_F)\geq 1$. Therefore, $p_{\infty}$ is also a steady state of the nonlinear equation in the dilated variable \eqref{eq:NNLIF-tau}, with $\tilde{N}=0$. Taking this steady state as the initial data, and then we clearly get a solution to the time-evolution problem with $\tilde{N}(\tau)\equiv0$ for all $\tau\geq 0$. Hence, the eternal blow-up \eqref{eternal-blow-up-def} happens and the corresponding lifespan of the generalized solution in the original time scale is zero as given in \eqref{lifespan-2.1.1}. 
\begin{corollary}\label{cor:eternal-bl}When $b\geq V_F-V_R$, there exists initial data satisfying Assumption \ref{as:classical-init}, such that the solution to \eqref{eq:NNLIF-tau} in the dilated timescale satisfies \begin{equation*}
     \tilde{N}(\tau)\equiv0,\quad \forall\tau\geq 0.
\end{equation*} Therefore, the lifespan of the corresponding generalized solution of \eqref{eq:NNLIF-original} $T^*$ is zero, as given in \eqref{lifespan-2.1.1}.
\end{corollary}
Using the steady state of \eqref{limit-linear} as a counter-example, Corollary \ref{cor:eternal-bl} denies the global well-posedness of the generalized solution when $b\geq V_F-V_R$. The first part of Theorem \ref{thm:gwp-t} is then proved. 

	\subsubsection{Relative entropy estimate}
	
	 The steady state \eqref{formula-steady-state} allows us to study the long time behavior of \eqref{limit-linear}{, which is a linear equation}, using the relative entropy method {as in \cite[Section 4]{caceres2011analysis} for the NNLIF equation}. We aim to derive a control on the boundary flux, {by modifying the treatments in \cite[Appendix A]{CACERES201481} and \cite[Theorem 2.1]{Antonio_Carrillo_2015}}. Estimates in this section serve as a basis for the global well-posedness proof in Section \ref{sec:42}.

	Let $p(v,\tau)$ be a solution of the limit equation \eqref{limit-linear}. We denote its boundary flux as
	\begin{equation}\label{def-M}
	    M(\tau):=-a_1\p_vp(v,\tau)\geq 0.
	\end{equation}
	Let $p_{\infty}(v)$ be the steady state \eqref{formula-steady-state} and $M_{\infty}$ be its boundary flux \eqref{def-Minfty}. We use the following notations for the ratios {as in \cite{Antonio_Carrillo_2015}}
	\begin{equation}\label{def-h-nu}
	    h(v,\tau):=\frac{p(v,\tau)}{p_{\infty}(v)},\quad \nu(\tau):=\frac{M(\tau)}{M_{\infty}}.
	\end{equation}
	For a $C^2$ convex function $G:\mathbb{R}  \rightarrow \mathbb{R}$, we denote the relative entropy
	\begin{equation}\label{def-GR}
		   S_G[p](\tau):=\int_{-\infty}^{V_F}p_{\infty}G(h(v,\tau))dv.
	\end{equation}
	
	Calculating $\frac{d}{d\tau}S_G[p](\tau)$, we have the entropy dissipation as summarized in the following. 
	\begin{proposition}
	\label{prop:entropy-linear} 
	Let $p(v,\tau)$ with its boundary flux $M(\tau)$ be a solution of the limit equation \eqref{limit-linear} and let $p_{\infty}$ be the steady state \eqref{formula-steady-state} with its boundary flux $M_{\infty}$. For a $C^2$ convex function $G:\mathbb{R}  \rightarrow \mathbb{R}$, let $S_G[p](t)$ be the relative entropy as defined in \eqref{def-GR}, we have \begin{equation*}
		       \frac{d}{d\tau}S_G[p](\tau)=-D_G[p](\tau)\leq 0,
		    \end{equation*} where the entropy dissipation $D_G[p](\tau)$ is given by
		    \begin{equation}\label{entropy-dis}
		        \begin{aligned}
		       D_G[p](\tau)&:=\int_{-\infty}^{V_F}a_1G''(h)|\p_vh|^2p_{\infty}dv\\&+M_{\infty}\bigl[G(\nu(\tau))-G(h(V_R,\tau))-G'(h(V_R,\tau))\bigl(\nu(\tau)-h(V_R,\tau)\bigr)\bigr]\geq 0.
		        \end{aligned}
		    \end{equation}
	\end{proposition}
	\begin{proof}[Proof of Proposition \ref{prop:entropy-linear}] We start from checking the boundary conditions of $h$.
	At $v=V_F$ since $0=p(V_F,\tau)=p_{\infty}(V_F,\tau)$, by L.H\^{o}pital's rule we get
	\begin{equation*}
	    h(V_F,\tau)=\frac{\p_vp(V_F,\tau)}{\p_vp_{\infty}(V_F)}=\frac{M}{M_{\infty}}=\nu(\tau).
	\end{equation*}
At $v=V_R$ we rewrite the flux jump condition for $p$ as
\begin{equation*}
a_1\p_v(hp_{\infty})|_{V_R^+}-a_1\p_v(hp_{\infty})|_{V_R^-}=a_1\p_v(hp_{\infty})|_{V_F},
\end{equation*} and obtain
\begin{equation*}
    a_1\p_vh(V_R^+,\tau)-a_1\p_vh(V_R^-,\tau)=-\frac{M_{\infty}}{p_{\infty}(V_R)}(\nu(\tau)-h(V_R,\tau)).
\end{equation*}
Then writing $p=hp_{\infty}$ in \eqref{limit-linear}, we derive an equation for $h$
\begin{equation}\label{eq:h}
\begin{aligned}
     \p_{\tau}h=-b\p_vh+&a_1\p_{vv}h+2a_1\p_vp_{\infty}/p_{\infty}\p_vh\\&+\frac{M_{\infty}}{p_{\infty}(V_R)}(\nu(\tau)-h(V_R,\tau))\delta_{v=V_R}.
\end{aligned}
\end{equation}
Multiplying \eqref{eq:h} by $G'(h)$, we have
\begin{equation}\label{eq:Gh}
    \begin{aligned}
     \p_{\tau}G(h)=-b\p_vG(h)+&a_1\p_{vv}G(h)+2a_1\p_vp_{\infty}/p_{\infty}\p_vG(h)-a_1|\p_vh|^2G''(h)\\&+\frac{M_{\infty}}{p_{\infty}(V_R)}(\nu(\tau)-h(V_R,\tau))G'(h(V_R,\tau))\delta_{v=V_R}.
\end{aligned}
\end{equation}
Then we multiply \eqref{eq:Gh} by $p_{\infty}$, and integrate by parts. Noting that $\p_{vv}$ is understood in the distributional sense, we deduce
\begin{equation}\label{tmp-IH}
    \begin{aligned}
         \frac{d}{d\tau}S_G[p](\tau)&=-\int_{-\infty}^{V_F}a_1|\p_vh|^2G''(h)dv\\&+\int_{-\infty}^{V_F}dv\left(a_1\p_{vv}(G(h)p_{\infty})+M_{\infty}(\nu(\tau)-h(V_R,\tau))G'(h)\delta_{v=V_R}+M_{\infty}G(h)\delta_{v=V_R}\right).
    \end{aligned}
\end{equation}
Finally, by the boundary conditions for $h$ and $p_{\infty}$, we obtain
\begin{align*}
        \int_{-\infty}^{V_F}a_1\p_{vv}(G(h)p_{\infty})dv&=a_1\p_v(G(h)p_{\infty})|_{v=V_F}\\&=a_1\p_vp_{\infty}(V_F)G(h(V_F,\tau))=-M_{\infty}G(\nu(\tau)).
\end{align*}
Plugging this to \eqref{tmp-IH}, we get \eqref{entropy-dis}. The two terms in \eqref{entropy-dis} are non-negative because $G$ is convex.
	\end{proof}
	Taking $G(x)=(x-1)^2$ in Proposition \ref{prop:entropy-linear}, we get
	\begin{equation}\label{entropy-quadratic}
	\begin{aligned}
		   \frac{d}{d\tau}\left(\int_{-\infty}^{V_F}p_{\infty}(h(v,\tau)-1)^2dv\right)=\frac{d}{d\tau}S_G[p](\tau)&= -D_G[p](\tau)\\&=-2\int_{-\infty}^{V_F}a_1|\p_vh|^2p_{\infty}dv-M_{\infty}(\nu(\tau)-h(V_R,\tau))^2.
	\end{aligned}
	\end{equation}
	The first term in this dissipation can be estimated via the following Poinc\'{a}re inequality. 
	\begin{proposition}\label{prop:poincare}When $b>0$, let $p_{\infty}$ be the steady state \eqref{formula-steady-state} of \eqref{limit-linear}. For $h(v)\in H^1((-\infty,V_F);p_{\infty}(v)dv)$ with $\int_{-\infty}^{V_F}(h-1)p_{\infty}dv=0$, we have the following Poinc\'{a}re inequality
	\begin{equation*}
	\alpha \int_{-\infty}^{V_F}p_{\infty}(h-1)^2dv\leq \int_{-\infty}^{V_F}p_{\infty}|\p_vh|^2dv,
	\end{equation*} where $\alpha>0$ is a constant independent of $h$.  
	\end{proposition}
	The proof of this Poinc\'{a}re inequality is similar to that of \cite[Proposition 4.3]{caceres2011analysis} and is postponed to Appendix \ref{app:poincare}.
	
	Together with the entropy dissipation \eqref{entropy-quadratic}, Proposition \ref{prop:poincare} implies the exponential convergence to the steady state, in the norm induced by the entropy $S_G[p](\tau)$ with $G=(x-1)^2$. 
	\begin{corollary}\label{cor:exp-convergence}
	When $b>0$, for a solution $p(v,\tau)$ of \eqref{limit-linear} with the initial data satisfying 
	\begin{equation}\label{init-exp-linear}
	    \int_{-\infty}^{V_F}p(v,0)dv=1,\qquad \int_{-\infty}^{V_F}p(v,0)^2\frac{dv}{p_{\infty}(v)}<\infty,
	\end{equation} Then $p(\cdot,\tau)$ converges to the steady state $p_{\infty}$ exponentially in the following sense 
	\begin{equation*}
	    \int_{-\infty}^{V_F}(p(v,\tau)-p_{\infty}(v))^2\frac{dv}{p_{\infty}(v)}\leq e^{- 2\alpha \tau} \left(\int_{-\infty}^{V_F}(p(v,0)-p_{\infty}(v))^2\frac{dv}{p_{\infty}(v)}\right),\quad \tau\geq 0.
	\end{equation*}
	\end{corollary}
	\begin{proof}Recall the notation $h(v,\tau):=\frac{p(v,\tau)}{p_{\infty}(v)}$ in \eqref{def-h-nu}. Thanks to the conservation of mass of \eqref{limit-linear} and the first condition in \eqref{init-exp-linear}, we have $\int_{-\infty}^{V_F}(h(v,\tau)-1)p_{\infty}(v)dv=\int_{-\infty}^{V_F}(p(v,\tau)-p_{\infty}(v))dv=0$. Note that \eqref{init-exp-linear} ensures $S_G[p](0)<+\infty$ with $G=(x-1)^2$. Applying Proposition \ref{prop:poincare} to $h(\cdot,\tau)$ in the entropy dissipation \eqref{entropy-quadratic}, which gives $\frac{d}{d\tau}S_G[p](\tau)\leq-2\alpha S_G[p](\tau)$, we conclude the result by Gronwall's inequality. 
	\end{proof} 
	
	This convergence result suggests that in the long time the boundary flux $M(\tau)=-a_1\p_vp(V_F,\tau)$ as defined in \eqref{def-M} shall also be close to that of the steady state $M_{\infty}$. However, the convergence in $L^2((-\infty,V_F),\frac{dv}{p_{\infty}(v)})$ given by Corollary \ref{cor:exp-convergence} is not strong enough to imply that $M(\tau)$ also converges to $M_{\infty}$. Instead, using the treatments in \cite{CACERES201481,Antonio_Carrillo_2015}, we control $\nu(\tau)=M(\tau)/M_{\infty}$ using entropy dissipation $D_G[p](\tau)$ as follows.
	\begin{lemma}\label{lemma:control-nu}
	When $b>0$, for $G=(x-1)^2$, there exists some constant $\eps>0$ such that
	\begin{equation}
	    D_G[p](\tau)\geq \int_{-\infty}^{V_F}a_1|\p_vh|^2p_{\infty}dv+\eps M_{\infty}(\nu(\tau)-1)^2,
	\end{equation} for any $p(v,\tau)$ with $\int_{-\infty}^{V_F}p(v,\tau)dv=1$.
	\end{lemma}
	\begin{proof}[Proof of Lemma \ref{lemma:control-nu}]
	For some $0<\eps<\frac{1}{2}$ to be determined, by the elementary inequality $a^2+2\eps b^2\geq \eps (a-b)^2$, we get
	\begin{equation}
	    M_{\infty}(\nu(\tau)-h(V_R,\tau))^2\geq \eps M_{\infty}(\nu(\tau)-1)^2-2\eps M_{\infty}(h(V_R,\tau)-1)^2.
	\end{equation}
	By Proposition \ref{prop:poincare} we have
	\begin{equation}\label{tmp-lm41}
	    \int_{-\infty}^{V_F}a_1|\p_vh|^2p_{\infty}dv\geq  \frac{1}{2}\int_{-\infty}^{V_F}a_1|\p_vh|^2p_{\infty}dv+\frac{\alpha}{2}\int_{-\infty}^{V_F}a_1(h-1)^2p_{\infty}dv.
	\end{equation}
	 Next, we apply the Sobolev injection $H^1(I)\subset L^{\infty}(I)$ to $h(v)-1$ on a small neighbourhood $I$ of $V_R=0$. The interval $I$ is chosen such that $p_{\infty}(v)$ as in \eqref{formula-steady-state} is bounded from below on $I$. Then noting $\p_v(h-1)=\p_vh$, from \eqref{tmp-lm41} we obtain
	\begin{equation}
	    C\int_{-\infty}^{V_F}|\p_vh|^2p_{\infty}dv\geq (h(V_R,\tau)-1)^2,
	\end{equation} with some $C>0$.
	Thanks to the expression of $D_G[p](t)$ in \eqref{entropy-quadratic}, choosing $\eps$ such that $2\eps M_{\infty}C<a_1$, and $\eps<1/2$, we get the desired results.
	\end{proof}
	Then we can show $M(\tau)$ is close to $M_{\infty}$ in the long term as follows.
	\begin{corollary}\label{cor:limit-linear-integral}
	When $b>0$, for a solution $p(v,\tau)$ of the limit equation \eqref{limit-linear} with the initial data satisfying \eqref{init-exp-linear}, we have 
	\begin{equation}
	    \int_0^{+\infty}(M(\tau)-M_{\infty})^2d\tau<\infty,
	\end{equation} where $M(\tau):=-a_1\p_vp(V_F,\tau)$ is the boundary flux as in \eqref{def-M}.
	\end{corollary}
	\begin{proof}
	By Lemma \ref{lemma:control-nu}, we get $\eps M_{\infty} (\nu(\tau)-1)^2\leq D_G[p](\tau)$ with $G=(x-1)^2$. Then by Proposition \ref{prop:entropy-linear} we obtain
	\begin{align*}
	    	    \eps \int_0^TM_{\infty}(\nu(\tau)-1)^2d\tau\leq \int_0^T D_G[p](\tau)d\tau&= \int_0^T-\frac{d}{d\tau}S_G[p](\tau)d\tau\\&=S_G[p](0)-S_G[p](T)\leq S_G[p](0),
	\end{align*} since $S_G[p]$ is non-negative. Taking the $T\rightarrow+\infty$ limit, we deduce
	\begin{equation*}
	    \int_0^{+\infty}M_{\infty}(\nu(\tau)-1)^2d\tau\leq\frac{1}{\eps}S_G[p](0)<+\infty,
	\end{equation*} since condition \eqref{init-exp-linear} ensures $S_G[p](0)<+\infty$. Note that $(M-M_{\infty})^2=M_{\infty}^2(\nu(\tau)-1)^2$. The proof is completed since $M_{\infty}>0$ as in \eqref{def-Minfty}.
	\end{proof}

\subsection{Global well-posedness in $t$ in {the} mildly excitatory case}\label{sec:42}

With the results on the linear limit equation \eqref{limit-linear} as preparations, this section is devoted to the global well-posedness of the generalized solution in the original timescale $t$ \eqref{eq:NNLIF-original} when $0<b<V_F-V_R$. In other words, we aim to prove the mildly excitatory case of Theorem \ref{thm:gwp-t}.

As shown in Proposition \ref{prop:gen-uni-c}, the lifespan of a generalized solution is independent of the parameter $c>0$ in the change of variable $d\tau=(N+c)dt$. Hence, without loss of generality we take $c=\frac{a_0}{a_1}$ to set $a_c=0$ in view of \eqref{def-bcac}, therefore simplifying the diffusion coefficient. With this choice $c=\frac{a_0}{a_1}$, the equation in $\tau$ \eqref{eq:NNLIF-tau} becomes 
\begin{equation}\label{eq:NNLIF-tau-simple}
    \begin{aligned}
        \p_{\tau}p+\p_v[(-v\tilde{N}+b^*\tilde{N}+b)p]=&a_1\p_{vv}p\\&+a_1(-\p_vp(V_F,\tau))\delta_{v=V_R},\quad v\in(-\infty,V_F),\tau>0,
        \end{aligned}
\end{equation} with the boundary condition $p(V_F,\tau)=0,\,\tau>0$. Here in place of $b_c$ in \eqref{def-bcac}, $b^*$ is given by 
\begin{equation}\label{def-b^*}
    b^*=b_0-\frac{a_0}{a_1}b.
\end{equation}
Then the drift condition \eqref{cond-drift-exc} is equivalent to 
\begin{equation}
    b^*\leq V_F.
\end{equation}
As as in \eqref{def-M}, we denote 
\begin{equation}
    M(\tau)=-a_1\p_{v}p(V_F,\tau).
\end{equation}
Then together with $c=\frac{a_0}{a_1}$, from \eqref{def-tildeN} we have
\begin{equation}\label{tildeN-M}
\tilde{N}(\tau)=\frac{(1-M(\tau))_+}{\frac{a_0}{a_1}M(\tau)+\frac{a_0}{a_1}(1-M(\tau))_+}=\frac{a_1}{a_0}\frac{(1-M(\tau))_+}{M(\tau)+(1-M(\tau))_+}\in [0,\frac{a_1}{a_0}].
\end{equation}

We argue by contradiction. Suppose the solution does not globally exist, then in view of the expression of its lifespan \eqref{lifespan-2.1.1}, we have
\begin{equation}\label{contra-gwpt-cond}
    \int_0^{+\infty}\tilde{N}(\tau)d\tau<\infty,
\end{equation} which implies $\int_T^{+\infty}\tilde{N}(\tau)d\tau\rightarrow0$ as $T\rightarrow+\infty$. This allows us to see \eqref{eq:NNLIF-tau-simple} as a non-autonomous perturbation from the linear limit equation \eqref{limit-linear}, with the additional term $\p_v[(-v+b^*)\tilde{N}p]$. {By \eqref{contra-gwpt-cond}, this perturbation shall diminish as time evolves.}

Hence, under the assumption \eqref{contra-gwpt-cond}, the long time behavior of \eqref{eq:NNLIF-tau-simple} is expected to be similar to that of \eqref{limit-linear}. In particular in the long term the boundary flux $M(\tau)$ should also be close to $M_{\infty}$ \eqref{def-Minfty}, which is less than one when $0<b<V_F-V_R$. However, this would imply some positive lower bound of $\tilde{N}$ thanks to \eqref{tildeN-M}, which contradicts the assumption \eqref{contra-gwpt-cond}. Our proof is to fulfill this intuitive plan in a rigorous way.

First we note that it is sufficient to get an estimate similar to Corollary \ref{cor:limit-linear-integral}.

\begin{lemma}\label{lemma:contra-exc} When $0<b<V_F-V_R$, if a solution of \eqref{eq:NNLIF-tau-simple} satisfies 
\begin{equation}\label{variance-cond}
    \int_0^{+\infty}(M(\tau)-M_{\infty})^2d\tau<\infty,
\end{equation} then \eqref{contra-gwpt-cond} can not hold. Here $M(\tau)=-a_1\p_vp(V_F,\tau)$ as in \eqref{def-M}.
\end{lemma}
\begin{proof}[Proof of Lemma \ref{lemma:contra-exc}]Let $\delta:=1-M_{\infty}$. By Proposition \ref{prop:steady-excit}, $\delta>0$ when $0<b<V_F-V_R$. Therefore \eqref{variance-cond} implies that the Lebesgue measure of the set 
$$A:=\{\tau\geq 0:M(\tau)\geq 1-\delta/2 \}\subseteq\{\tau\geq 0:|M(\tau)-M_{\infty}|\geq \delta/2 \}$$ is finite. However when $\tau\notin A$ and $\tau\geq 0$, we have $0\leq M(\tau)\leq 1-\frac{\delta}{2}$ and thus $\tilde{N}(\tau)\geq \frac{a_1}{a_0}\frac{\delta}{2}$ by \eqref{tildeN-M}, which implies that
\begin{equation*}
    \int_0^{+\infty}\tilde{N}(\tau)d\tau\geq \int_{\mathbb{R}^+\backslash A} \frac{a_1}{a_0}\frac{\delta}{2}dt=+\infty,
\end{equation*} Therefore \eqref{contra-gwpt-cond} does not hold.
\end{proof}

By this lemma, we thus aim to prove \eqref{variance-cond} by carrying out the relative entropy estimate to \eqref{eq:NNLIF-tau-simple}. We need to control the extra terms involving $\tilde{N}(\tau)$. To this aim, we first give the $L^{\infty}$ bounds of $h=p/p_{\infty}$, using the method of super-solutions developed in \cite{Antonio_Carrillo_2015}. Then, we do the entropy estimate. Under the assumption \eqref{contra-gwpt-cond} we prove \eqref{variance-cond}, which in turn contradicts \eqref{contra-gwpt-cond} and therefore the proof is completed.

\subsubsection{$L^{\infty}$ bound}

We establish the $L^{\infty}$ bounds in this section as preparations for the entropy estimate.

 First, we construct super-solutions of \eqref{eq:NNLIF-tau-simple} similar to \cite[Defintion 4.2 and Lemma 4.4]{Antonio_Carrillo_2015}. {The super-solution consists of two ingredients: the steady state \eqref{formula-steady-state} of the limit equation \eqref{limit-linear} and a time factor, where the former is to tackle the terms without $\tilde{N}$ and the latter is to deal with the terms involving $\tilde{N}$.}
\begin{proposition}\label{prop:super-solution} When $b>0$ and $b^*\leq V_F$, there exists a $\gamma>0$ such that for any given solution $p$ and $\tilde{N}$ of \eqref{eq:NNLIF-tau-simple},  
\begin{equation}\label{eq:formula-super}
    q(v,\tau):=\exp\left(\gamma\int_0^{\tau}\tilde{N}(s)ds\right)p_{\infty}(v)\geq 0
\end{equation} is a super solution of \eqref{eq:NNLIF-tau-simple} in the following sense: $q(V_F,\tau)=0$ and 
\begin{equation}\label{def-super}
        	    \p_{\tau}q+\p_v((-v\tilde{N}+b^*\tilde{N}+b)q)-a_1\p_{vv}q\geq -a_1\p_vq(V_F,\tau)\delta_{v=V_R},
\end{equation} on $(-\infty,V_F)\times(0,+\infty)$ in the sense of distributions. Here $p_{\infty}$ is the steady state of the limit equation \eqref{limit-linear} given in \eqref{formula-steady-state}.
\end{proposition}
\begin{proof}[Proof of Proposition \ref{prop:super-solution}]
We have $q(V_F,\tau)=0$ since $p_{\infty}(V_F)=0$, i.e., the steady state \eqref{formula-steady-state} satisfies the boundary condition.

To check \eqref{def-super} in the sense of distribution, as in Remark \ref{rmk:flux-jump} it suffices to check that $q$ satisfies \eqref{def-super} in the classical sense on $((-\infty,V_R)\cap(V_R,V_F))\times(0,+\infty)$ and a flux jump inequality at $V_R$ 
\begin{equation}
    -a_1\p_vq(V_R^{+},\tau)+a_1\p_vq(V_R^{-},\tau)\geq -a_1\p_vq(V_F^{-},\tau).
\end{equation}
Since $p_{\infty}(v)$ is the steady state of \eqref{limit-linear}, the flux jump condition is satisfied (it becomes an equality) and the terms without $\tilde{N}$ in \eqref{def-super} cancel. Therefore, we only need to check the extra terms related with $\tilde{N}$
\begin{equation*}
    \gamma \tilde{N} q-\tilde{N}q+\tilde{N}(-v+b^*)\p_vq\geq 0,
\end{equation*} which is equivalent to 
\begin{equation}\label{tmp-super}
        \gamma \tilde{N} p_{\infty}-\tilde{N}p_{\infty}+\tilde{N}(-v+b^*)\p_vp_{\infty}\geq 0,
\end{equation} on $(-\infty,V_R)\cup(V_R,V_F)$.

On $(-\infty,V_R)$, from \eqref{formula-steady-state} $\p_vp_{\infty}=\frac{b}{a_1}p_{\infty}\geq 0$. Note that $-v\geq -V_R$ when $v\leq V_R$. Thus in this case it is sufficient to choose $\gamma>1+\frac{b}{a_1}|b^*-V_R|$. 

On $(V_R,V_F)$, by \eqref{formula-steady-state} $p_{\infty}=\frac{1}{Z}(e^{\frac{b}{a_1}V_F}-e^{\frac{b}{a_1}v})e^{-\frac{b}{a_1}V_R}$, which is zero at $V_F$. Since in this case $\p_vp_{\infty}(v)\leq 0$, the last term in \eqref{tmp-super} is positive if $v\geq \max(V_R,b^*)$. Therefore it suffices to choose
\begin{equation*}
    \gamma>1+\sup_{v\in(V_R,b^*)}\frac{(b^*-v)\frac{b}{a_1}e^{\frac{b}{a_1}v}}{e^{\frac{b}{a_1}V_F}-e^{\frac{b}{a_1}v}}.
\end{equation*} Noting that the $\sup$ term above is finite when $b^*\leq V_F$, we conclude the proof.
\end{proof}We remark that $\gamma>0$ in Proposition \ref{prop:super-solution} does not depend on a particular solution, while the super-solution $q(v,\tau)$ depends on a specific $\tilde{N}$.

By comparing a solution with a super-solution, we can derive $L^{\infty}$ bounds as in \cite[Lemma 4.3]{Antonio_Carrillo_2015}.
\begin{proposition}\label{prop:linfty}When $b>0$ and $b^*\leq V_F$, suppose the initial data $p_I(v)$ of \eqref{eq:NNLIF-tau-simple} satisfy \begin{equation}\label{cond-init}
    p_I(v)\leq C_Ip_{\infty}(v)
\end{equation} for some $C_I>0$. Then the following $L^{\infty}$ bound holds
\begin{equation}
    p(v,\tau)\leq C_I\exp\left(\gamma\int_0^{\tau}\tilde{N}(s)ds\right)p_{\infty}(v),\quad \tau\geq 0,v\in(-\infty,V_F].
\end{equation} Here $\gamma$ is the constant for the super-solution in Proposition \ref{prop:super-solution}.
\end{proposition}
\begin{proof}[Proof of Proposition \ref{prop:linfty}]
 Let $\bar{q}(v,\tau)=C_Iq(v,\tau)$, where $q$ is the super-solution constructed in Proposition \ref{prop:super-solution}. Therefore $w(v,\tau):=\bar{q}(v,\tau)-p(v,\tau)$ satisfies $w(v,0)\geq 0$ and is also a super-solution of \eqref{eq:NNLIF-tau-simple} in the sense of Proposition \ref{prop:super-solution}.

Note that with a given $\tilde{N}$, \eqref{eq:NNLIF-tau-simple} is a linear equation which keeps non-negativity thanks to standard convex entropy estimates {(see for example \cite[Section 4]{caceres2011analysis} in a similar context)}. We deduce that $w(v,\tau)\geq 0$ for $\tau\geq0,v\in(-\infty,V_F]$, which gives the desired estimate. 
\end{proof}

Under the contradiction assumption \eqref{contra-gwpt-cond}, Proposition \ref{prop:linfty} directly implies a uniform-in-time $L^{\infty}$ bound as follows.
\begin{corollary}\label{cor:Linfinity}
When $b>0$ and $b^*\leq V_F$, suppose the initial data of \eqref{eq:NNLIF-tau-simple} satisfy $p_I(v)\leq C_Ip_{\infty}(v)$ for some $C_I>0$. If additionally, \eqref{contra-gwpt-cond} holds, then {for the solution $p$ of \eqref{eq:NNLIF-tau-simple}} we have 
\begin{equation}
    p(v,\tau)\leq Cp_{\infty}(v),\quad \tau\geq0,v\in(-\infty,V_F],
\end{equation} for some constant $C$ independent of time.
\end{corollary}
\begin{proof}
Due to \eqref{contra-gwpt-cond}, for any $\tau\geq 0$, $\int_0^{\tau}\tilde{N}(s)ds\leq \int_0^{+\infty}\tilde{N}(s)ds<+\infty$. Thus the result follows from Proposition \ref{prop:linfty}.
\end{proof}
\subsubsection{Relative entropy estimate. Proof of Theorem \ref{thm:gwp-t}-(2)-(a)}

In this part, we extend the relative entropy estimate of the limit equation \eqref{limit-linear} to \eqref{eq:NNLIF-tau-simple}. Then we fulfill the aforementioned contradiction argument to establish the global well-posedness for generalized solutions of \eqref{eq:NNLIF-original} when $0<b<V_F-V_R$, proving Theorem \ref{thm:gwp-t}-(2)-(a), 

For a solution $p(v,\tau)$ of \eqref{eq:NNLIF-tau-simple} with its boundary flux $M(\tau)$ \eqref{def-M}, as in \eqref{def-h-nu} we use the notations
\begin{equation*}
    h(v,\tau)=\frac{p(v,\tau)}{p_{\infty}(v)},\quad \nu(\tau)=\frac{M(\tau)}{M_{\infty}}.
\end{equation*} Here $p_{\infty}$ is the steady state of the limit equation \eqref{limit-linear} and $M_{\infty}$ is its boundary flux as in Proposition \ref{prop:steady-excit}. And as in \eqref{def-GR} we consider the relative entropy for $p$ against $p_{\infty}$ 
\begin{equation*}
    S_G[p](\tau):=\int_{-\infty}^{V_F}p_{\infty}G(h(v,\tau))dv.
\end{equation*}Similar to Proposition \ref{prop:entropy-linear} we calculate $\frac{d}{d\tau}S_G[p](\tau)$. The difference is that now $p(v,\tau)$ is governed by \eqref{eq:NNLIF-tau-simple} instead of \eqref{limit-linear}, thus terms involving $\tilde{N}$ appear. 
	\begin{proposition}\label{prop:entropy-NNLIF} Let $p(v,\tau)$ with its boundary flux $M(\tau)$ be a solution of the nonlinear equation \eqref{eq:NNLIF-tau-simple} and let $p_{\infty}$ be the steady state \eqref{formula-steady-state} of the linear limit equation \eqref{limit-linear} with its boundary flux $M_{\infty}$. For a $C^2$ convex function $G:\mathbb{R}  \rightarrow \mathbb{R}$, let $S_G[p](\tau)$ be the relative entropy as defined in \eqref{def-GR}. Then we have
		    \begin{equation} \frac{d}{d\tau}S_G[p](\tau)=-D_G[p](\tau)+\mathcal{E}_G[p](\tau),
		    \end{equation} where the dissipation $D_G[p](\tau)$ is still given by  \eqref{entropy-dis}
		    \begin{equation*}
		        \begin{aligned}
		   D_G[p](\tau)&:=\int_{-\infty}^{V_F}a_1G''(h)|\p_vh|^2p_{\infty}dv\\&+M_{\infty}[G(\nu(\tau))-G(h(V_R,\tau)-G'(h(V_R,\tau))(\nu(\tau)-h(V_R,\tau))].
		        \end{aligned}
		    \end{equation*}
		    And the ``perturbation'' $\mathcal{E}_G[p](\tau)$ is given by 
		    \begin{align}\label{def-mathcalE}
    \mathcal{E}_G[p](\tau)&=-\tilde{N}(\tau)\int_{-\infty}^{V_F}\bigl(G'(h)h-G(h)\bigr)\p_v\bigl((-v+b^*)p_{\infty}\bigr)dv\\&=-\tilde{N}(\tau)\int_{-\infty}^{V_F}\bigl(G'(h)h-G(h)\bigr)\bigl((-v+b^*)\p_vp_{\infty}-p_{\infty}\bigr)dv\label{def-mathcalE-2}
		    \end{align}
\end{proposition}
\begin{proof}[Proof of Proposition \ref{prop:entropy-NNLIF}]
We rewrite \eqref{eq:NNLIF-tau-simple} as a perturbation from \eqref{limit-linear}.
\begin{equation}\label{NNLIF-fast-a1-Perturb}
\begin{aligned}
         	    \p_{\tau}p=-b\p_vp&+a_1\p_{vv}p-a_1\p_vp(V_F,\tau)\delta_{v=0}\\&-\tilde{N}\p_v((-v+b^*)p),\quad \tau>0,v\in(-\infty,V_F).
\end{aligned}
\end{equation}
Note that $p$ for \eqref{NNLIF-fast-a1-Perturb} has the same boundary condition and the flux jump condition as \eqref{limit-linear}. Thus here $h(v,\tau)$ has the same boundary conditions as in Proposition \ref{prop:entropy-linear}. Therefore terms without $\tilde{N}$ can be manipulated in the exact same way. As in Proposition \ref{prop:entropy-linear}, we multiply \eqref{NNLIF-fast-a1-Perturb} by $G'(h)$ and integrate to obtain
\begin{equation}
    \frac{d}{dt}S_G[p](\tau)=-D_G[p](\tau)+\mathcal{E}_G[p](\tau),
\end{equation}where the dissipation $D_G[p]$ is the same as given in \eqref{entropy-dis}. And the ``extra'' $\tilde{N}$ terms give rise to to $\mathcal{E}_G[p](\tau)$, which is given by 
\begin{align*}
    \mathcal{E}_G[p](\tau)=-\tilde{N}(\tau)\int_{-\infty}^{V_F}\p_v((-v+b^*)p)G'(h)dv.
\end{align*}Thanks to the zero boundary condition at $V_F$, we integrate by parts
\begin{align*}
    \mathcal{E}_G[p](\tau)&=-\tilde{N}(\tau)\int_{-\infty}^{V_F}-(-v+b^*)pG''(h)\p_vhdv\\&=-\tilde{N}(\tau)\int_{-\infty}^{V_F}-(-v+b^*)p_{\infty}hG''(h)\p_vhdv
\end{align*}Noting that $hG''(h)\p_vh=\p_v(G'(h)h-G(h))$, we do another integration by parts to transfer the derivative to $p_{\infty}$
\begin{align*}
    \mathcal{E}_G[p](\tau)&=-\tilde{N}(\tau)\int_{-\infty}^{V_F}\bigl(G'(h)h-G(h)\bigr)\p_v\bigl((-v+b^*)p_{\infty}\bigr)dv\\&=-\tilde{N}(\tau)\int_{-\infty}^{V_F}\bigl(G'(h)h-G(h)\bigr)\bigl((-v+b^*)\p_vp_{\infty}-p_{\infty}\bigr)dv\notag
\end{align*} where the boundary term vanishes due to that $p_{\infty}(V_F)=0$.  

\end{proof}

Now we give the proof of Theorem \ref{thm:gwp-t} in the mildly excitatory case, i.e., when $0<b<V_F-V_R$.
\begin{proof}[Proof of Theorem \ref{thm:gwp-t}. mildly excitatory case]
First we check that for initial data $p_I(v)$, there exists a constant $C_I$ such that
\begin{equation}\label{proof-inital-ex}
    p_I(v)\leq C_Ip_{\infty}(v),\quad v\in(-\infty,V_F].
\end{equation} Due to the continuity of $p_I$ and expression of $p_{\infty}$ \eqref{formula-steady-state}. It suffices to check $v=V_F$ and $v\rightarrow-\infty$. For $v=V_F$, it follows from that $p_I(V_F)=0$, $\frac{d}{dv}p_I(V_F^{-})$ is finite by Assumption \ref{as:classical-init} and $-\p_vp_{\infty}(V_F^{-})>0$ by Proposition \ref{prop:steady-excit}. For $v\rightarrow-\infty$, it is ensured by \eqref{cond-tail-exc}.

Now suppose the contrary the solution is not global, by Proposition \ref{prop:gen-uni-c} in the dilated timescale $\int_0^{+\infty}\tilde{N}(\tau)dt<\infty$. Then thanks to \eqref{proof-inital-ex} by Corollary \ref{cor:Linfinity} we have the uniform $L^{\infty}$ bound
\begin{equation}
    h(v,\tau)\leq C<\infty,\quad v\in(-\infty,V_F),\tau>0.
\end{equation} 

Then we control the ``perturbation'' $\mathcal{E}_G[p](\tau)$  \eqref{def-mathcalE-2} with $G=(x-1)^2$. Thanks to the uniform bound of $h$, $hG'(h),G(h)$ are also bounded uniformly. Hence we derive
\begin{align*}
        |\mathcal{E}_G[p](\tau)|&\leq \tilde{N}(\tau)\left|\int_{-\infty}^{V_F}\bigl(G'(h)h-G(h)\bigr)\bigl((-v+b^*)\p_vp_{\infty}-p_{\infty}\bigr)dv\right|,\\&\leq  C\tilde{N}(\tau)\left|\int_{-\infty}^{V_F}|(-v+b^*)\p_vp_{\infty}|+p_{\infty}dv\right|\\ &\leq C \tilde{N}(\tau).
\end{align*} In the last line, to conclude the integrability of these two terms, we use the formula for $p_{\infty}$ \eqref{formula-steady-state}, in particular its exponential decay at infinity.

Therefore by the relative entropy calculation in Proposition \ref{prop:entropy-NNLIF}, we have
\begin{equation}
   0\leq D_G[p](\tau)\leq -\frac{d}{d\tau}S_G[p](\tau)+C\tilde{N}(\tau).
\end{equation}
Thus we get 
\begin{align*}
        \int_0^{+\infty}D_G[p](\tau)d\tau&\leq S_G[p](0)-\liminf_{t\rightarrow +\infty}S_G[p](\tau)+C\int_0^{+\infty}\tilde{N}(\tau)d\tau\\&\leq  S_G[p](0)+C\int_0^{+\infty}\tilde{N}(\tau)d\tau<\infty,
\end{align*} where $S_G[p](0)<\infty$ thanks to \eqref{proof-inital-ex} and $\int_0^{+\infty}\tilde{N}(\tau)d\tau<\infty$ is the contradiction assumption.

Then, using the entropy dissipation $D_G[p]$ to control the boundary flux by Lemma \ref{lemma:control-nu}, we deduce
\begin{equation}
  \int_0^{+\infty}(M(\tau)-M_{\infty})^2d\tau=\frac{1}{M_{\infty}}\int_0^{+\infty}M_{\infty}(\nu(\tau)-1)^2d\tau\leq  C\int_0^{+\infty}D_G[p](\tau)d\tau<\infty.
\end{equation}
Finally, by Lemma \ref{lemma:contra-exc} this leads to a contradiction with the assumption $\int_0^{+\infty}\tilde{N}(\tau)d\tau<\infty$. Then the proof is completed.
\end{proof}

\section{Inhibitory case: Global well-posedness in $t$}\label{sec:5i}

In this section we prove the global well-posedness for the generalized solution of \eqref{eq:NNLIF-original} in the inhibitory case, i.e., when $b\leq 0$, which is the last part of Theorem \ref{thm:gwp-t}. The proof strategy is similar to the mildly excitatory case in Section 4 but some adjustments are needed. A main difference lies in the long time behavior of the limit equation \eqref{limit-linear}.

As in Section 4, we first study the limit equation \eqref{limit-linear}. Then we are able to treat the equation in timescale $\tau$ as a perturbation of \eqref{limit-linear}.

\subsection{Limit equation}
We start with the limit equation \eqref{limit-linear} as in Section 4.1
\begin{equation*}
	\begin{aligned}
	  \p_{\tau}p+b\p_vp&=a_1\p_{vv}p-a_1\p_vp(V_F,\tau)\delta_{v=V_R},\quad \tau>0,v\in(-\infty,V_F),\\
	    p(V_F,\tau)&=0,\ p(-\infty,\tau)=0,\ \tau>0.
	\end{aligned}
\end{equation*} 
\subsubsection{Non-$L^1$ positive steady state}
	As in the excitatory case, the steady state of \eqref{limit-linear} satisfies \eqref{steady-limit}
	\begin{equation*}
	\begin{aligned}
	\p_v(bp)&=a_1\p_{vv}p-a_1\p_vp(V_F,\tau)\delta_{v=V_R},\quad v\in(-\infty,V_F),\\
	p(V_F)&=0,\ p(-\infty)=0.
	\end{aligned}
	\end{equation*}
	
	However, in the inhibitory case when $b\leq 0$, the steady equation \eqref{steady-limit} does not have a solution, unless we drop the requirement on the decay when $v\rightarrow-\infty$. Nevertheless, if we relax this condition, we can still find positive steady states, summarized as the following. 
	
\begin{proposition}\label{prop:steady-inhibitory} When $b\leq 0$, the limit equation \eqref{limit-linear} does not have steady states in $L^1(-\infty,V_F)$. Nevertheless, if we do not require the decay at infinity, we have the following ``steady state'' which is positive for $v<V_F$. When $b<0$ it is given by 
\begin{equation}\label{formula-steady-b<0}
    p_{\infty}(v)=\begin{cases}
    \bigl(1-e^{\frac{b}{a_1}(V_F-V_R)}\bigr)e^{\frac{b}{a_1}(v-V_R)},\quad v\leq V_R,\\
    e^{\frac{b}{a_1}(v-V_R)}-e^{\frac{b}{a_1}(V_F-V_R)},\quad v\in(V_R,V_F].
    \end{cases}
\end{equation}
While when $b=0$ the steady state is given by
\begin{equation}\label{formula-steady-b=0}
    p_{\infty}(v)=\begin{cases}
    V_F-V_R,\quad v\leq V_R,\\
    V_F-v,\quad v\in(V_R,V_F].
    \end{cases}
\end{equation}
\end{proposition}
\begin{proof}The calculation is similar to Proposition \ref{prop:steady-excit}. WLOG we can set $V_R=0$ by a translation argument.
	On $(-\infty,0)\cup(0,V_F)$, the solution satisfies
	\begin{equation}
	    \p_{vv}p-\frac{b}{a_1}\p_vp=0.
	\end{equation}
	
	We consider the case $b<0$ first. On $(-\infty,0)$, $p$ is given by $c_1e^{\frac{b}{a_1}v}+c_0$. On $(0,V_F)$ by the boundary condition at $V_F$, $p$ is given by $c_2(e^{\frac{b}{a_1}V_F}-e^{\frac{b}{a_1}v})$. At $v=0$, thanks to the flux jump condition 
	\begin{equation*}
	    a_1\p_vp(0^+)-a_1\p_vp(0^-)=a_1\p_vp(V_F),
	\end{equation*} we get that $c_1=(e^{\frac{b}{a_1}V_F}-1)c_2$, which in turn gives $c_0=0$ due to the continuity of $p$ at $v=0$. Since $b<0$ if we want an integrable steady state, we have to set $c_1=0$, which also make $c_2=c_1/(e^{\frac{b}{a_1}V_F}-1)=0$. Then we only get the zero solution. On the other hand, if we choose $c_2=-1$, we obtain the positive non-$L^1$ steady state \eqref{formula-steady-b<0}.
	
	For $b=0$, similarly we have $p(v)=c_1v+c_0$ on $(-\infty,0)$ and $p(v)=c_2(v-V_F)$ on $(0,V_F)$. The flux jump condition gives $c_1=0$. The continuity of $p$ at $v=0$ gives $c_0=-V_Fc_2$. Thus a non-zero steady state can not be integrable. And choosing $c_2=-1,c_0=V_F$ we get \eqref{formula-steady-b=0}.
	
\end{proof}

We can not normalize these ``steady states'' $p_{\infty}(v)$ when $b\leq0$ since they are not in $L^1$. Nevertheless, for later references we fix the expressions of $p_{\infty}$ as in Proposition \ref{prop:steady-inhibitory} and denote its boundary flux as
\begin{equation*}
    M_{\infty}=-a_1\p_vp_{\infty}(v)>0,
\end{equation*} using the same notation as in \eqref{def-Minfty}.
\subsubsection{Relative entropy estimate}
Though the ``steady state'' $p_{\infty}$ in Proposition \ref{prop:steady-inhibitory} is not in $L^1$, we can still use it to derive the relative entropy estimate. By the relative entropy estimate we can see the long time behavior of \eqref{limit-linear} and control the boundary flux, as in Section 4.1.2 and \cite{caceres2011analysis,CACERES201481,Antonio_Carrillo_2015}.

The entropy dissipation calculation in Proposition \ref{prop:entropy-linear}, which does not depend on the sign of $b$ or any specific property of $p_{\infty}$, still holds when $b\leq 0$ with $p_{\infty}$ defined in Proposition \ref{prop:steady-inhibitory}. We recall the result here as follows. 

Let $p(v,\tau)$ be a solution of the limit equation \eqref{limit-linear} and $p_{\infty}$ be as in Proposition \ref{prop:steady-inhibitory}, with their boundary fluxes $M(\tau)$ and $M_{\infty}$. We still use the notations in \eqref{def-h-nu}
    \begin{equation*}
	    h(v,\tau):=\frac{p(v,\tau)}{p_{\infty}(v)},\quad \nu(\tau):=\frac{M(\tau)}{M_{\infty}}.
	\end{equation*}	For a $C^2$ convex function $G:\mathbb{R}  \rightarrow \mathbb{R}$, let the relative entropy associated with $G$ be as in \eqref{def-GR} \begin{equation*}
		   S_G[p](\tau)=\int_{-\infty}^{V_F}p_{\infty}G(h(v,\tau))dv.
		 \end{equation*}
		    Then as in Proposition \ref{prop:entropy-linear} \begin{equation}\label{decrease-Sg}
		   \frac{d}{dt}S_G[p](\tau)=-D_G[p](\tau)\leq 0,
		    \end{equation} where the dissipation $D_G[p](\tau)$ is given by \eqref{entropy-dis}
		    \begin{equation*}
		  \begin{aligned}D_G[p](\tau)&:=\int_{-\infty}^{V_F}a_1G''(h)|\p_vh|^2p_{\infty}dv\\&+M_{\infty}[G(\nu(\tau))-G(h(V_R,\tau))-G'(h(V_R,\tau))(\nu(\tau)-h(V_R,\tau))].
		        \end{aligned}
		    \end{equation*}
		    
The difference is in the choice of the convex function $G$. Now $p_{\infty}$ is not in $L^1$, if we still choose $G(x)=(x-1)^2$, then the entropy $ S_G[p](\tau)$ will not be a finite number in general. Instead, we shall choose $G(x)=x^2$ which gives
\begin{align}\label{x2entropy-2}
    S_G[p](\tau)=\int_{-\infty}^{V_F}p_{\infty}h^2dv&=\int_{-\infty}^{V_F}p^2(v,\tau)\frac{dv}{p_{\infty}(v)}\\&=\int_{-\infty}^{V_F}p(v,\tau)h(v,\tau)dv,\label{x2entropy}
\end{align}which is finite, provided that, e.g., $h(v,\tau)$ is bounded. 

From \eqref{decrease-Sg} we know $S_G[p](\tau)$ is decreasing in time. Let's examine the indication of this on the long time behavior. Plugging in the expressions of $p_{\infty}$ into \eqref{x2entropy-2}, when $b<0$ we obtain
\begin{equation}\label{weighted-b<0}
\begin{aligned}
             S_G[p](\tau)=\frac{1}{1-\exp(\frac{b}{a_1}(V_F-V_R))}&\int_{-\infty}^{V_R}p^2(v,\tau)\frac{dv}{\exp(\frac{b}{a_1}(v-V_R))}\\+&\int_{V_R}^{V_F}p^2(v,\tau)\frac{dv}{\exp(\frac{b}{a_1}(v-V_R))-\exp(\frac{b}{a_1}(V_F-V_R))}.
\end{aligned}
\end{equation} And when $b=0$ we have
\begin{equation}\label{weighted-b=0}
   S_G[p](\tau)=\int_{-\infty}^{V_R}p^2(v,\tau)\frac{dv}{V_F-V_R}+\int_{V_R}^{V_F}p^2(v,\tau)\frac{dv}{V_F-v}.
\end{equation}We can see $S_G[p]$ as a weighed integral of $p^2$ with the weight ${1}/{p_{\infty}(v)}$ by \eqref{x2entropy-2}. Intuitively, the decreasing of $S_G[p](\tau)$ may indicate that the mass of $p$ tends to move from where the weight is large to where the weight is small. Note that when $b<0$, $p_{\infty}$ goes to infinity as $v\to-\infty$ and thus ${1}/{p_{\infty}(v)}$ goes to zero when $v$ goes to $-\infty$. Hence the dissipation of this entropy may tell that the mass shall move to $-\infty$ as time evolves. For the case $b=0$ in \eqref{weighted-b=0}, the weight at $v=V_F$ is also the highest. Therefore when $b\leq0$, we shall expect that the mass near $V_F$ becomes small in the long time, which intuitively implies the boundary flux $M(\tau)$ shall be close to zero.  

Indeed, using the entropy dissipation $D_G[p](\tau)$, we can control the boundary flux $M(\tau)$ via an inhibitory counterpart of Lemma \ref{lemma:control-nu}. We can control the distance between $\nu(\tau)=M(\tau)/M_{\infty}$ and a $\delta_K(\tau)$ which can be made small. Precisely, we have the following lemma.
	\begin{lemma}\label{lemma:control-nu-inh}
	When $b\leq 0$, for $G=x^2$, for any $K>0$ there exists $\eps=\eps(K)>0$ such that
	\begin{equation}
	    D_G[p](\tau)\geq \int_{-\infty}^{V_F}a_1|\p_vh|^2p_{\infty}dv+\eps M_{\infty}\left(\nu(\tau)-\delta_K(\tau)\right)^2,
	\end{equation} for any $p(v,\tau)$ with $\int_{-\infty}^{V_F}p(v,\tau)dv=1$. Here $\delta_K(\tau)$ is the mean of $h$ on $[V_R-K,V_R]$ which satisfies
	\begin{equation}\label{def-deltaK}
	    \delta_K(\tau):=\frac{1}{K}\int_{V_R-K}^{V_R}h(v,\tau)dv\leq \frac{C}{K},
	\end{equation} with a constant $C$ independent of $p(v,\tau)$ or $K$.
	\end{lemma}
	\begin{proof}[Proof of Lemma \ref{lemma:control-nu-inh}]
	By formulas in Proposition \ref{prop:steady-inhibitory}, $1/p_{\infty}(v)$ is bounded for $v\leq V_R=0$, therefore 
	\begin{equation*}
	\begin{aligned}
	    \delta_K(\tau):=\frac{1}{K}\int_{V_R-K}^{V_R}h(v,\tau)dv&\leq \frac{C}{K}\int_{V_R-K}^{V_R}p(v,\tau)dv\\&\leq \frac{C}{K}\int_{-\infty}^{V_F}p(v,\tau)dv\leq \frac{C}{K}.
	\end{aligned}
	\end{equation*}
	For some $0<\eps<\frac{1}{2}$ to be determined, by the elementary inequality $a^2+2\eps b^2\geq \eps (a-b)^2$, we deduce in the same way as for Lemma \ref{lemma:control-nu}
	\begin{equation}
	    M_{\infty}(\nu(\tau)-h(V_R,\tau))^2\geq \eps M_{\infty}(\nu(\tau)-\delta_K(\tau))^2-2\eps M_{\infty}(h(V_R,\tau)-\delta_K(\tau))^2.
	\end{equation}
	By definition $\delta_K(\tau)$ is the average of $h$ on $[V_R-K,V_R]$. Thus, by Poincar\'{e}-Wirtinger inequality on $[V_R-K,V_R]$ and noting that $p_{\infty}(v)$ is bounded from below on $[V_R-K,V_R]$, we have
		\begin{equation}
	    \int_{-\infty}^{V_F}a_1|\p_vh|^2p_{\infty}dv\geq C\int_{V_R-K}^{V_R}(h-\delta_K(\tau))^2dv.
	\end{equation}
	Together with the Sobolev injection from $H^1(I)$ to $L^{\infty}(I)$ applying to $h(\cdot,\tau)-\delta_K(\tau)$ on the interval $I=[V_R-K,V_R]$, we get for some $C_K>0$
	\begin{equation}
	    C_K\int_{-\infty}^{V_F}|\p_vh|^2p_{\infty}dv\geq (h(V_R,\tau)-\delta_K(\tau))^2.
	\end{equation}	Finally, note that the entropy dissipation $D_G[p]$ for $G=x^2$ is the same as that for $G=(x-1)^2$, as in \eqref{entropy-quadratic}
	\begin{equation*}
	    D_G[p](\tau)=2\int_{-\infty}^{V_F}a_1|\p_vh|^2p_{\infty}dv+M_{\infty}(\nu(\tau)-h(V_R,\tau))^2.
	\end{equation*}
	Thus by choosing $\eps>0$ small enough such that $2\eps M_{\infty}C_K<a_1$, and $\eps<1/2$ we get the desired results.
	\end{proof}
	By Lemma \ref{lemma:control-nu-inh} we can obtain a global control on $M(\tau)$ similar to Corollary \ref{cor:limit-linear-integral}.
	\begin{corollary}\label{cor:limit-linear-integral-inh} When $b\leq 0$, suppose $p(v,\tau)$ is a solution to the limit equation \eqref{limit-linear} with the finite initial entropy \eqref{x2entropy}, i.e., $S_G[p](0)<+\infty$ with $G(x)=x^2$. Then for every $\eps>0$ there exists a function in time $\delta(\tau)$ satisfying $0\leq \delta(\tau)<\eps$ for all $\tau\geq 0$ such that
	\begin{equation}
	    \int_0^{+\infty}(M(\tau)-\delta(\tau))^2d\tau<+\infty,
	\end{equation} where $M(\tau)=-a_1\p_vp(V_F,\tau)$ is the boundary flux as in \eqref{def-M}.
	\end{corollary}
\begin{proof}
By Lemma \ref{lemma:control-nu-inh} if we choose $K=K(\eps)$ large enough, we have $0\leq M_{\infty}\delta_K(\tau)<\eps$ for all $\tau\geq 0$. And in an identical way as Lemma \ref{lemma:control-nu}, we get
\begin{align*}
    \int_0^{+\infty}(M(\tau)-M_{\infty}\delta_K(\tau))^2d\tau \leq C \int_0^{+\infty}D_G[p](\tau)d\tau \leq C S_G[p](0)<+\infty.
\end{align*}
Thus we obtain the desired result by taking $\delta(\tau)=M_{\infty}\delta_K(\tau)$.
\end{proof}

\subsection{Global well-posedness in $t$ in the inhibitory case}

With preparations on the limit equation \eqref{limit-linear}, now we start to analyze the equation in the dilated timescale \eqref{eq:NNLIF-tau}. As in Section 4.2, we set $c=\frac{a_0}{a_1}$ WLOG and thus work with
\begin{equation*}
        \begin{aligned}
        \p_{\tau}p+\p_v[(-v\tilde{N}+b^*\tilde{N}+b)p]=&a_1\p_{vv}p\\&+a_1(-\p_vp(V_F,\tau))\delta_{v=V_R},\quad v\in(-\infty,V_F),\tau>0,
        \end{aligned}
\end{equation*} {with} the boundary condition $p(V_F,\tau)=0,\,\tau>0$. We still use the notation $M(\tau)=-a_1\p_{v}p(V_F,\tau)$ \eqref{def-M} for the boundary flux, and recall in this case $\tilde{N}$ is determined by $M$ in \eqref{tildeN-M}
\begin{equation*}
\tilde{N}(\tau)=\frac{(1-M(\tau))_+}{\frac{a_0}{a_1}M(\tau)+\frac{a_0}{a_1}(1-M(\tau))_+}=\frac{a_1}{a_0}\frac{(1-M(\tau))_+}{M(\tau)+(1-M(\tau))_+}\in [0,\frac{a_1}{a_0}].
\end{equation*}

A difference lies in the assumption on the drift $b^*$, which is given by \eqref{def-b^*}\begin{equation*}
    b^*=b_0-\frac{a_0}{a_1}b.
\end{equation*} In this inhibitory case, instead of $b^*\leq V_F$ \eqref{cond-drift-exc}, we assume \eqref{cond-drift-inh}
\begin{equation}
    b_0\leq V_F,
\end{equation}which implies in \eqref{eq:NNLIF-tau-simple}
\begin{equation*}
    b^*\tilde{N}+b=b_0\tilde{N}+b(1-\frac{a_0}{a_1}\tilde{N})\leq b_0\tilde{N}\leq V_F\tilde{N},
\end{equation*}thanks to \eqref{def-b^*} and $\tilde{N}\in[0,\frac{a_1}{a_0}]$ \eqref{tildeN-M}.

We argue by contradiction as in the excitatory case. Suppose the solution does not globally exists, in view of the lifespan formula \eqref{lifespan-2.1.1}, we have \eqref{contra-gwpt-cond}
\begin{equation*}
        \int_0^{+\infty}\tilde{N}(\tau)d\tau<\infty,
\end{equation*} which gives a smallness for $\tilde{N}(\tau)$. Therefore we may see \eqref{eq:NNLIF-tau-simple} as a perturbation of the limit equation \eqref{limit-linear}. However, as indicated in Section 5.1, when $b\leq0$ a solution of \eqref{limit-linear} moves away from $V_F$ in the long time, which would lead to small boundary flux $M(\tau)$. In view of \eqref{tildeN-M}, this would result in some positive lower bound of $\tilde{N}$ \eqref{tildeN-M}, which contradicts \eqref{contra-gwpt-cond}. As in the mildly excitatory case, our proof is a justification of this intuition.

First, we note that it suffices to get an estimate similar to Corollary \ref{cor:limit-linear-integral-inh}.
\begin{lemma}\label{lemma:contra-inh} If a solution of \eqref{eq:NNLIF-tau-simple} satisfies 
\begin{equation}\label{variance-cond-inh}
    \int_0^{+\infty}(M(\tau)-\delta(\tau))^2dt<\infty,
\end{equation}where $\delta(\tau)$ is a function on $[0,+\infty)$ satisfying $$\delta(\tau)\leq \delta_0<1,\quad \tau\geq0,$$ then \eqref{contra-gwpt-cond} can not holds. Here $M(\tau)=-a_1\p_vp(V_F,\tau)$ as in \eqref{def-M}.
\end{lemma}
Lemma \ref{lemma:contra-inh} is the inhibitory counterpart of Lemma \ref{lemma:contra-exc}, and its proof, which we omit here, is essentially the same as Lemma \ref{lemma:contra-exc}.

To derive \eqref{variance-cond-inh}, we first establish $L^{\infty}$ bounds and then use the relative entropy estimate.

\subsubsection{$L^{\infty}$ bound}

To derive the $L^{\infty}$ bounds, first we construct super-solutions similar to \cite{Antonio_Carrillo_2015} and Proposition \ref{prop:super-solution}.

\begin{proposition}\label{prop:super-solution-inh} When $b\leq0$ and $b_0\leq V_F$, there exists a $\gamma>0$ such that for a given solution $p$ and $\tilde{N}$ of \eqref{eq:NNLIF-tau-simple},  
\begin{equation}\label{eq:formula-super-inh}
    0\leq q(v,\tau):=\begin{cases}
    \exp\left(\gamma\int_0^t\tilde{N}(s)ds\right),\quad v\leq V_R,\\
    \exp\left(\gamma\int_0^t\tilde{N}(s)ds\right)\frac{V_F-v}{V_F-V_R},\quad v\in(V_R,V_F].
    \end{cases}
\end{equation} is a super solution of \eqref{eq:NNLIF-tau-simple} in the following sense: $q(V_F,\tau)=0$ and 
\begin{equation}\label{def-super-inh}
        	    \p_{\tau}q+\p_v((-v\tilde{N}+b^*\tilde{N}+b)q)-a_1\p_{vv}q\geq -a_1\p_vq(V_F,\tau)\delta_{v=V_R},
\end{equation} on $(-\infty,V_F)\times(0,+\infty)$ in the sense of distributions.
\end{proposition}
\begin{proof}[Proof of Proposition \ref{prop:super-solution-inh}]
By \eqref{eq:formula-super-inh} we have $q(V_F,\tau)=0$, it remains to show \eqref{def-super-inh}

Similar to Proposition \ref{prop:super-solution}, it suffices to check that $q$ satisfies \eqref{def-super-inh} in the classical sense on $((-\infty,V_R)\cap(V_R,V_F))\times(0,+\infty)$ and a flux jump inequality:
\begin{equation*}
    -a_1\p_vq(V_R^+,\tau)+a_1\p_vq(V_R^-,\tau)\geq a_1\p_vq(V_F^-,\tau).
\end{equation*}It is direct to check that the flux jump condition is satisfied. Indeed for a fixed $\tau$, $q(\cdot,\tau)$ is a steady state of \eqref{limit-linear} with $b=0$. Therefore, we only need to check the following inequality
\begin{equation}\label{tmp-super-2}
    \gamma \tilde{N} q-\tilde{N}q+(-v\tilde{N}+b^*\tilde{N}+b)\p_vq\geq 0,
\end{equation} on the two intervals $(-\infty,V_R)\cup(V_R,V_F)$.

On $(-\infty,V_R)$, $\p_vq=0$ thus \eqref{tmp-super-2} reduces to 
\begin{equation*}
    \gamma \tilde{N} q -\tilde{N}q\geq 0.
\end{equation*}And it suffices to choose $\gamma\geq1$.

On $(V_R,V_F)$, in this case $\p_vq<0$. We rewrite the drift coefficient by plugging in the definition of $b^*$ \eqref{def-b^*}
\begin{equation*}
    -v\tilde{N}+b_0\tilde{N}+b(1-\frac{a_0}{a_1}\tilde{N})\leq -v\tilde{N}+b_0\tilde{N},
\end{equation*} since $b\leq 0$ and $\tilde{N}\in[0,\frac{a_1}{a_0}]$ by \eqref{tildeN-M}. Therefore, in this case owing to $\p_vq\leq 0$, it suffices to ensure
\begin{equation*}
        \gamma \tilde{N} q-\tilde{N}q+(-v\tilde{N}+b_0\tilde{N})\p_vq\geq 0,
\end{equation*}
The last term is positive if $v\geq \max(V_R,b_0)$ thanks to $\p_vq\leq 0$. Therefore it is sufficient to choose
\begin{equation*}
    \gamma>1+\sup_{v\in(V_R,b_0)}\frac{(b_0-v)(-\p_vq)}{q}=1+\sup_{v\in(V_R,b_0)}\frac{b_0-v}{V_F-v}.
\end{equation*} Note that the $\sup$ term above is finite thanks to the assumption $b_0\leq V_F$. We conclude the proof.
\end{proof}
\begin{remark}[Comparison to the excitatory case]
In Proposition \ref{prop:super-solution}, we construct super-solutions for $b>0$ via the steady state of the associated limit equation \eqref{limit-linear}. In this inhibitory case, we construct super-solutions for all $b\leq 0$ by the steady state of \eqref{limit-linear} when $b=0$, which is closer to the setting of ``universal super-solution'' in \cite[Lemma 4.4]{Antonio_Carrillo_2015}. This choice in some sense uses the benefits of $b\leq 0$.
\end{remark}

With super-solutions, we give $L^{\infty}$ bounds as in Proposition \ref{prop:linfty} and \cite[Lemma 4.3]{Antonio_Carrillo_2015}.
\begin{proposition}\label{prop:linfty-inh}When $b\leq 0$ and $b_0\leq V_F$, if the initial data $p_I(v)$ of \eqref{eq:NNLIF-tau-simple} is bounded and satisfies \begin{equation}\label{cond-init-inh}
\limsup_{v\rightarrow V_F^-}\frac{p_I(v)}{V_F-v}<\infty.
\end{equation}Then the following $L^{\infty}$ bound holds
\begin{equation}
    p(v,\tau)\leq \begin{cases}
    C_I\exp\left(\gamma\int_0^t\tilde{N}(s)ds\right),\quad v\leq V_R,\\
    C_I\exp\left(\gamma\int_0^t\tilde{N}(s)ds\right)\frac{V_F-v}{V_F-V_R},\quad v\in(V_R,V_F],
    \end{cases}\quad t\geq 0,v\in(-\infty,V_F].
\end{equation} Here $\gamma$ is the constant for the super-solution in Proposition \ref{prop:super-solution-inh} and $C_I>0$ is a constant depending on the initial data.
\end{proposition}
\begin{proof}[Proof of Proposition \ref{prop:linfty-inh}]
Noting that by the condition \eqref{cond-init-inh} and boundedness of $p_{I}(v)$, we can choose $C_I$ large enough, such that
\begin{equation*}
    p_I(v)\leq C_Iq(v,0)=\begin{cases}
    C_I,\quad v\leq V_R,\\
    C_I\frac{V_F-v}{V_F-V_R},\quad v\in(V_R,V_F],
    \end{cases},\quad v\in(-\infty,V_F].
\end{equation*}
Then the rest of the proof is identical to that of Proposition \ref{prop:linfty}.
\end{proof}

Under the contradiction assumption \eqref{contra-gwpt-cond}, Proposition \ref{prop:linfty-inh} directly implies uniform-in-time $L^{\infty}$ bounds similar to Corollary \ref{cor:Linfinity} in the excitatory case. For later purposes, we give bounds on both $p$ and $h=p/p_{\infty}$. 
\begin{corollary}\label{cor:Linfinity-inh}
When $b\leq 0$ and $b_0\leq V_F$, for \eqref{eq:NNLIF-tau-simple} with bounded initial data satisfying \eqref{cond-init-inh}. If additionally, \eqref{contra-gwpt-cond} holds then we have the uniform bounds
\begin{equation}
    p(v,\tau)\leq C,\quad \tau\geq0,v\in(-\infty,V_F),
\end{equation} and
\begin{equation}
    p(v,\tau)\leq Cp_{\infty}(v),\quad \tau\geq0,v\in(-\infty,V_F),
\end{equation} 
for some constant $C$ independent of time, where $p_{\infty}$ is the ``steady state'' of \eqref{limit-linear} given in Proposition \ref{prop:steady-inhibitory}.
\end{corollary}
\begin{proof}
Due to \eqref{contra-gwpt-cond}, for any $\tau\geq 0$, $\int_0^{\tau}\tilde{N}(s)ds\leq \int_0^{+\infty}\tilde{N}(s)ds<+\infty$. Then the result follows from Proposition \ref{prop:super-solution-inh} and that $q(v,0)\leq C\min(1,p_{\infty}(v))$ by the formula in Proposition \ref{prop:steady-inhibitory}, where $q(v,0)=\begin{cases}
    1,\quad v\leq V_R,\\
    \frac{V_F-v}{V_F-V_R},\quad v\in(V_R,V_F],
    \end{cases}$.
\end{proof}

\subsubsection{Relative entropy estimate. Proof of Theorem \ref{thm:gwp-t}-(2)-(b)}

Now we prove the global well-poseness for generalized solutions of \eqref{eq:NNLIF-original} in the inhibitory case $b\leq 0$, completing the proof of Theorem \ref{thm:gwp-t}-(2)-(b).

We start with recalling the relative entropy calculations. For a solution $p(v,\tau)$ of \eqref{eq:NNLIF-tau-simple} with its boundary flux $M(\tau)$ \eqref{def-M}, as before we use the notations \eqref{def-h-nu}
\begin{equation*}
    h(v,\tau)=\frac{p(v,\tau)}{p_{\infty}(v)},\quad \nu(\tau)=\frac{M(\tau)}{M_{\infty}}.
\end{equation*} Here $p_{\infty}$ is the ``steady state'' of the linear equation \eqref{limit-linear}, whose expression is given in Proposition \ref{prop:steady-inhibitory}, and $M_{\infty}$ is its boundary flux. 

In the inhibitory case, we choose the convex function $G(x)=x^2$ and consider the associated relative entropy for $p$ against $p_{\infty}$ \eqref{def-GR}
\begin{equation}\label{inh-GR}
    S_G[p](\tau):=\int_{-\infty}^{V_F}p_{\infty}G(h(v,\tau))dv=\int_{-\infty}^{V_F}p^2(v,\tau)/p_{\infty}dv
\end{equation} Note that the calculation in Proposition \ref{prop:entropy-NNLIF}, which does not depend on the sign of $b$, also holds for $b\leq 0$, and here we recall
 \begin{equation}\label{dS/dt-inh} \frac{d}{d\tau}S_G[p](\tau)=-D_G[p](\tau)+\mathcal{E}_G[p](\tau),
\end{equation} where the dissipation $D_G[p](\tau)$ is given by  \eqref{entropy-dis} with $G(x)=x^2$
 \begin{equation}\label{tmp-pf-inh-dis}
		        \begin{aligned}
	    D_G[p](\tau)=2\int_{-\infty}^{V_F}a_1|\p_vh|^2p_{\infty}dv+M_{\infty}(\nu(\tau)-h(V_R,\tau))^2.
		        \end{aligned}
		    \end{equation}
		    And the ``perturbation'' $\mathcal{E}_G[p](\tau)$ is given by \eqref{def-mathcalE}
		  \begin{align}\notag
    \mathcal{E}_G[p](\tau)&=-\tilde{N}(\tau)\int_{-\infty}^{V_F}\bigl(G'(h)h-G(h)\bigr)\bigl((-v+b^*)\p_vp_{\infty}-p_{\infty}\bigr)dv\\&=-\tilde{N}(\tau)\int_{-\infty}^{V_F}h^2\bigl((-v+b^*)\p_vp_{\infty}-p_{\infty}\bigr)dv.\label{tmp-pf-mathcalE}
	\end{align} 

Now we start the proof, whose general strategy is identical to the mildly excitatory case. 
\begin{proof}[Proof of Theorem \ref{thm:gwp-t}. Inhibitory case]
Suppose the contrary the solution is not global, by the lifespan formula \eqref{lifespan-2.1.1} in Proposition \ref{prop:gen-uni-c}, in the dilated timescale \eqref{contra-gwpt-cond} holds, i.e., $\int_0^{+\infty}\tilde{N}(\tau)dt<\infty$..

First we claim that to get a contradiction, it suffices to show
\begin{equation}\label{pf-tmp-inh-DG}
    \int_0^{+\infty}D_G[p](\tau)d\tau<+\infty,
\end{equation} where $D_G[p]$ given by \eqref{tmp-pf-inh-dis} is the entropy dissipation with $G=x^2$. Indeed, by Lemma \ref{lemma:control-nu-inh}, we can choose $K$ large enough, such that $M_{\infty}\delta_K\leq \frac{1}{2}<1$ and that there exists $\eps>0$ satisfying
\begin{equation*}
    \frac{\eps}{M_{\infty}}(M(\tau)-M_{\infty}\delta_K(\tau))^2={\eps}{M_{\infty}}(\nu(\tau)-\delta_K(\tau))^2\leq D_G[p](\tau).
\end{equation*} Integrating in $\tau$, if \eqref{pf-tmp-inh-DG} holds we get
\begin{equation*}
    \int_0^{+\infty}(M(\tau)-M_{\infty}\delta_K(\tau))^2d\tau\leq C    \int_0^{+\infty}D_G[p](\tau)d\tau<+\infty,
\end{equation*} which contradicts \eqref{contra-gwpt-cond} thanks to Lemma \ref{lemma:contra-inh}, with the choice $\delta(\tau)=M_{\infty}\delta_K(\tau)$.

Next, by the entropy calculation \eqref{dS/dt-inh}, we have
\begin{align}\notag
   0\leq \int_0^{+\infty}D_G[p](\tau)&\leq S_G[p](0)-\limsup_{\tau\rightarrow+\infty}S_G[p](\tau)+\int_0^{+\infty} \bigl|\mathcal{E}_G[p](\tau)\bigr|d\tau\\ &\leq S_G[p](0)+\int_0^{+\infty} \bigl|\mathcal{E}_G[p](\tau)\bigr|d\tau,\label{pf-tmp-inh-two-term}
\end{align} since $S_G[p](\tau)\geq 0$. Then to show \eqref{pf-tmp-inh-DG}, it reduces to bound the two terms in \eqref{pf-tmp-inh-two-term}.

We claim that for the initial data $p_I(v)$, thanks to Assumption \ref{as:classical-init} there exists a constant $C_I$ such that
\begin{equation}\label{proof-inital-inh}
    p_I(v)\leq C_Ip_{\infty}(v),\quad v\in(-\infty,V_F].
\end{equation}Indeed, by the expression in Proposition \ref{prop:steady-inhibitory}, $p_{\infty}$ has positive lower bound on $(-\infty,V_F-\eps]$ for all $\eps>0$. Then thanks to the boundedness of $p_I$ by Assumption \ref{as:classical-init}, it suffices to check for $v$ in a left neighbourhood of $V_F$. This can be ensured by $p_I(V_F)=p_{\infty}(V_F)=0$ together with $\frac{d}{dv}p_I(V_F^-)<\infty,\frac{d}{dv}p_{\infty}(V_F^-)>0$ thanks to Assumption \ref{as:classical-init} and Proposition \ref{prop:steady-inhibitory}.

For the first term in \eqref{pf-tmp-inh-two-term}, \eqref{proof-inital-inh} implies 
\begin{equation*}
        S_G[p](0)=\int_{-\infty}^{V_F}p^2_I(v)/p_{\infty}(v)dv\leq C_I \int_{-\infty}^{V_F}p_I(v)dv=C_I<+\infty,
\end{equation*}since $p_I(v)$ is a probability density by Assumption \ref{as:classical-init}.

It remains to show the second term 
\begin{equation*}
    \int_0^{+\infty} \bigl|\mathcal{E}_G[p](\tau)\bigr|d\tau<+\infty.
\end{equation*} Using the contradiction assumption $\int_0^{+\infty}\tilde{N}(\tau)d\tau<+\infty$ \eqref{contra-gwpt-cond}, it suffices to derive a uniform-in-$\tau$ bound
\begin{equation*}
    \bigl|\mathcal{E}_G[p](\tau)\bigr|=\left|\int_{-\infty}^{V_F}h^2(v,\tau)\bigl((-v+b^*)\p_vp_{\infty}-p_{\infty}\bigr)dv\right|\leq C<+\infty,\quad \tau\geq 0.
\end{equation*}
By the triangle inequality, we have
\begin{align*}
   \left|\int_{-\infty}^{V_F}h^2\bigl((-v+b^*)\p_vp_{\infty}-p_{\infty}\bigr)dv\right|&\leq   \left|\int_{-\infty}^{V_F}h^2(-v+b^*)\p_vp_{\infty}dv\right|+ \left|\int_{-\infty}^{V_F}h^2p_{\infty}dv\right|.
\end{align*}
Thanks to \eqref{proof-inital-inh}, by Corollary \ref{cor:Linfinity-inh} both $p$ and $h=p/p_{\infty}$ are uniformly bounded. Therefore for the second term, we have
\begin{equation*}
    \left|\int_{-\infty}^{V_F}h^2p_{\infty}dv\right|=  \left|\int_{-\infty}^{V_F}hp dv\right|\leq C \left|\int_{-\infty}^{V_F}p dv\right|=C.
\end{equation*} And for the first term, we split the integral into two parts as follows
\begin{align*}
    \left|\int_{-\infty}^{V_F}h^2(-v+b^*)\p_vp_{\infty}dv\right|&\leq \left|\int_{-\infty}^{V_R}h^2(-v+b^*)\p_vp_{\infty}dv\right|+\left|\int_{V_R}^{V_F}h^2(-v+b^*)\p_vp_{\infty}dv\right|\\& \leq \left|\int_{-\infty}^{V_R}h^2(-v+b^*)\p_vp_{\infty}dv\right|+C.
\end{align*} Finally, note that when $v\leq V_R$, $|(-v+b^*)\p_vp_{\infty}|\leq Cp_{\infty}^2$ by the expression in Proposition \ref{prop:steady-inhibitory}. Thus we obtain
\begin{align*}
   \left|\int_{-\infty}^{V_R}h^2(-v+b^*)\p_vp_{\infty}dv\right| &\leq C\left|\int_{-\infty}^{V_R}h^2p_{\infty}^2dv\right|\\ &= C\int_{-\infty}^{V_R}p^2dv\leq C \int_{-\infty}^{V_R}p(v,\tau)dv<+\infty,
\end{align*} thanks to the uniform bound of $p$. Then the proof is complete.

\end{proof}

\section{Discussion}\label{sec:discussion}

\subsection{On possible definitions of generalized solutions}

In this work, we have introduced a generalized solution for the NNLIF model \eqref{eq:NNLIF-intro}. When the classical solution ceases to exist, there appears uncertainty in broadening the notion of solutions. Moreover, it is possible to incorporate certain additional mechanism in establishing a generalized solution without having conflicts with the classical one.   

We start with clarifying the importance of the reformulation in defining the generalized solution in this paper. Our generalized solution is based on introducing the dilated timescale $\tau$ via dividing by $N(t)+c$ with $c>0$, or the time dilation transform. Before we do the time dilation transform, we rewrite \eqref{eq:NNLIF-original} into \eqref{eq:NNLIF-original-rw}. This simple reformulation is actually crucial. Let's recall \eqref{eq:NNLIF-original}
\begin{equation*}
    \p_tp+\p_v[(-v+b_0+bN(t))p]=(a_0+a_1N(t))\p_{vv}p+N(t)\delta_{v=V_R},\quad v\in(-\infty,V_F),t>0,
\end{equation*} and the reformulation \eqref{eq:NNLIF-original-rw}
\begin{equation*}
\begin{aligned}
    \p_tp+\p_v[(-v+b_0+bN(t))p]=(a_0&+a_1N(t))\p_{vv}p\\& +(a_0+a_1N(t))(-\p_vp(V_F,t))\delta_{v=V_R},\quad v\in(-\infty,V_F),t>0.
\end{aligned}
\end{equation*}
The only difference is the coefficient of the Dirac source term at $v=V_R$. In \eqref{eq:NNLIF-original} it is $N(t)$ while in \eqref{eq:NNLIF-original-rw} it is $(a_0+a_1N(t))(-\p_vp(V_F,t)).$ These two coefficients are the equal by \eqref{def-N-eq} for classical solutions and thus the two formulations are equivalent in this case. 

However, \eqref{eq:NNLIF-original} and \eqref{eq:NNLIF-original-rw} give different time dilation transforms when $N=+\infty$. If we carry out the time dilation transform to \eqref{eq:NNLIF-original}, we formally obtain \begin{equation}\label{eq:flux-limited}
    \begin{aligned}
        \p_{\tau}\tp+\p_v[(-v\tilde{N}+b_c\tilde{N}+b)\tp]=(a_c\tilde{N}+a_1)\p_{vv}\tp +(1-c\tilde{N})\delta_{v=V_R},\quad v\in(-\infty,V_F),\tau>0,
        \end{aligned}
\end{equation}
where the coefficient of the Dirac source term is $N/(N+c)=1-c\tilde{N}$. While applying time dilation transform to \eqref{eq:NNLIF-original-rw}, we obtain \eqref{eq:NNLIF-tau}, where the coefficient of the Dirac source term is $(a_c\tilde{N}+a_1)(-\p_v\tp(V_F,\tau))$. The former flux $1-c\tilde{N}$ is always bounded while the latter can be arbitrarily large  if $-\p_v\tp(V_F,\tau)$ is large. In view of this, we may call \eqref{eq:flux-limited} the flux-limited version and \eqref{eq:NNLIF-tau} the flux-unlimited version.

The flux-unlimited version \eqref{eq:NNLIF-tau}, which is the topic of this paper, is mathematically more {akin} to the classical solution of the NNLIF model in the original timescale \eqref{eq:NNLIF-intro}. This allows us to develop its well-posed theory using tools {from the analysis of classical solutions for the NNLIF model} \cite{caceres2011analysis,carrillo2013CPDEclassical,CACERES201481,Antonio_Carrillo_2015}. 

 The flux-limited version \eqref{eq:flux-limited} also deserves further understanding. If we naively integrate in space, then formally the reset at $V_R$ can be less than the loss at $V_F$, which gives an apparent violation of mass conservation. Therefore the precise formulation of \eqref{eq:flux-limited} might be more sophisticated, but it might be a worthy challenge.

Another issue on defining a generalized solution is the refractory period. Physically, after a neuron spikes it enters a short refractory period, during which it does not respond to any stimulus. From a modeling perspective, one can impose a finite  or an infinitesimal refractory period. The latter case can be viewed as a limit of the former. If we add a finite refractory period, the model as well as its solution structure is substantially changed (see e.g. \cite{CACERES201481}). But an infinitesimal refractory period might be added to modify the definition of a generalized solution without changing the classical solution, since it only makes a difference when the firing rate blows up. A refractory period ensures that one neuron can not spike more than once in an instant of the original timescale $t$. Thus it might exclude the eternal blow-up defined in \eqref{eternal-blow-up-def} and make it possible to have global well-posedness in timescale $t$ for arbitrary large connectivity parameter $b$.

To sum up, we have discussed two issues on defining a generalized solutions. We observe that the existing definitions in  \cite{delarue2015particle,taillefumier2022characterization,sadun2022global} and this work are different in those two aspects. Thus, it is appealing to study the uniqueness criteria in finding a scientifically relevant solution within a certain context.

\subsection{Other further directions}

We discuss possible future directions to explore the NNLIF dynamics \eqref{eq:NNLIF-intro} in the framework of the generalized solution in this paper.

The equation in the dilated timescale $\tau$ \eqref{eq:NNLIF-tau} might serve as an ideal platform to study the long time behavior of the original system, since the time-dilated transform  \eqref{eq:NNLIF-tau} is globally well-posed in the classical sense for all $b\in\mathbb{R}$ (Theorem \ref{thm:gwp-t}).

Moreover, \eqref{eq:NNLIF-tau} enriches the structure of the steady states of \eqref{eq:NNLIF-intro}. It not only has classical steady states with $N<+\infty$ but also allows steady states with $N=+\infty$, which is equivalent to $\tilde{N}=0$. These additional steady states appear if and only if $b\geq V_F-V_R$, and are exactly those given by Proposition \ref{prop:steady-excit} since they are also steady states to the limit equation \eqref{limit-linear}.

Interestingly, different solution behaviors in timescale $t$ might be unified as the long time convergence to a steady state in timescale $\tau$: If a solution converges to a steady state with $\tilde{N}>0$ (equivalent to $N<+\infty$), then it is global in timescale $t$. But if a solution converges to a steady state with $\tilde{N}=0$ (equivalent to $N=\infty$), it might have an eternal blow-up in timescale $t$ as defined in \eqref{eternal-blow-up-def}. Besides, the additional steady state might provide a way to understand the various phenomena including bistability observed in \cite{CiCP-30-820}.

To obtain the generalized solution, we first solve the equation in timescale $\tau$, then switch to the dynamics in timescale $t$ via the inverse change of time. This procedure also facilitates efficient numerical simulations. Indeed, by simulating the equation in $\tau$, the numerical scheme in \cite{hu2021structure} naturally becomes a computational approach for the generalized solution, which works for the dynamics with or without blow-ups.

Other interesting directions include to investigate how the blow-up times are distributed in the original timescale $t$ (for similar models this has been studied in \cite{delarue2022global,sadun2022global}), and to apply the idea of such time-dilation transforms to other integrate-and-fire models \cite{CACERES201481,Dumont2013Population,carrillo2022simplified}.

\section*{Acknowledgement}
Z. Zhou is supported by the National Key R\&D Program of China, Project Number 2020YFA0712000, 2021YFA1001200 and NSFC grant Number 12031013, 12171013. X.Dou is partially supported by The Elite Program of
Computational and Applied Mathematics for PhD Candidates in Peking University. We thanks Yantong Xie for his help in numerical simulations. We also thank
Jos\'e Carrillo and Pierre Roux for helpful discussions.

\begin{appendices}

\section{Postponed calculations and proofs from Section 3}\label{app:pf-sec3}

\subsection{Calculations and proofs for the transform to the free boundary problem}
\subsubsection{Calculations from (\ref{NNLIF-3.1}) to (\ref{reduce-q})}\label{app:pf-sec3-cal}

Here we present the calculations from \eqref{NNLIF-3.1} to \eqref{reduce-q}, via the change of variable \eqref{change-p-q}. Recall \eqref{change-p-q} is given by
\begin{equation*}
    \beta(\tau)q(\beta(\tau)v,S(\tau))=p(v,\tau),
\end{equation*} where $\beta(\tau)$ and $S(\tau)$ are given by \eqref{def-beta} and \eqref{def-S}.

For clarity, we split the transform \eqref{change-p-q} into two changes of variables, first in space \eqref{change-1} then in time \eqref{change-2}. Recall \eqref{NNLIF-3.1}
\begin{equation*}
    \p_\tau p+\p_v((-v\tilde{N}+\tilde{\mu})p)=\tilde{a}\p_{vv}p-\tilde{a}\p_vp(0,\tau)\delta_{v=V_R},\quad v<0,\, \tau>0,
\end{equation*}

\paragraph{Step 1} First we do a change of variable in the spatial variable $v$ to get rid of the $-v\tilde{N}(\tau)$ term in the drift. Let $\beta(\tau)$ as in \eqref{def-beta}
 \begin{equation*}
\beta(\tau):=\exp\left(\int_0^{\tau}\tilde{N}({\tau}')d{\tau}'\right),
 \end{equation*} which satisfies $\beta'(\tau)=\tilde{N}(\tau)\beta(\tau)$ and 
    $1\leq \beta(\tau)\leq e^{\tau/c}$,
 thanks to $\tilde{N}\in[0,1/c]$ \eqref{def-tildeN}. 
 
 Let $y=\beta(\tau)v$ be the new spatial variable and consider
\begin{equation}\label{change-1}
    p(v,\tau)=\beta(\tau)w(\beta (\tau)v,\tau).
\end{equation}
We calculate
\begin{align*}
    \p_{\tau}p=\beta \p_\tau w&+v\tilde{N}\beta^2 \p_yw+\tilde{N}\beta w,\\
    \tilde{N}\p_v(-vp)&=\tilde{N}(-v\beta^2\p_yw-\beta w),\\\p_{vv}p&=\beta^3\p_{yy}w,\quad \delta_{v=V_R}=\beta\delta_{y=\beta(\tau)V_R}.
\end{align*}
Substituting these into \eqref{NNLIF-3.1} we get
\begin{equation*}
    \beta \p_\tau w+\beta^2\tilde{\mu}\p_yw=\beta^3\tilde{a}\p_{yy}w-\beta^3\tilde{a}\p_yw(y=0,\tau)\delta_{y=\beta(\tau)V_R},\quad y<0,\,\tau>0,
\end{equation*} which simplifies to 
\begin{equation}\label{reduce-w}
     \p_\tau w+\beta\tilde{\mu}\p_yw=\beta^2\tilde{a}\p_{yy}w-\beta^2\tilde{a}\p_yw(y=0,\tau)\delta_{y=\beta(\tau)V_R},\quad y<0,\,\tau>0,
\end{equation}

 \paragraph{Step 2} Now we do another change of variable, in the time variable $\tau$. The aim is to simplify the diffusion coefficient in \eqref{reduce-q} to $1$. With a given $\tilde{N}(\tau)$, we define $S(\tau)$ as \eqref{def-S}
\begin{equation*}
    S(\tau):=\int_0^{\tau}\beta^2(\tau')\tilde{a}(\tau')d\tau',
\end{equation*} 
with $\beta(\tau)$ defined in \eqref{def-beta} and $\tilde{a}$ defined in \eqref{def-tilde-a}. We expect $\tilde{N}(\tau)$ is continuous in time as in Definition \ref{def:classical-tau}. Therefore $S(\tau)$ is a $C^1$ function, whose derivative has a positive lower bound since $\beta\geq 1$ thanks to $\tilde{N}\geq 0$ and $\tilde{a}\geq \min(a_0/c,a_1)>0$. Then we define its inverse map $T(s)$ as in \eqref{def-T'}, whose derivative is
\begin{equation*}
    T'(s)=\frac{1}{\beta^2(T(s))\tilde{a}(T(s))}.
\end{equation*}

We introduce the new time variable $s=S(\tau)$ and consider \begin{equation}\label{change-2}
    q(y,S(\tau))=w(y,\tau).
\end{equation} Then 
\begin{equation*}
    \beta^2\tilde{a} \p_sq=\p_{\tau}w,\quad \p_yq=\p_yw,\p_{yy}q=\p_{yy}w.
\end{equation*} Substituting these to \eqref{reduce-w} and we get \eqref{reduce-q}.

\subsubsection{Proofs of Lemma \ref{lemma:existence} and Lemma \ref{lemma:key-bound}}

Here we present the proof of Lemma \ref{lemma:existence} and Lemma \ref{lemma:key-bound}.

\begin{proof}[Proof of Lemma \ref{lemma:existence}]
Since $M$ is given, to be concise we drop the dependence on $M$ in the notations, i.e., writing $\beta(s)=\beta(s;M)$ and $\tilde{N}(s)=\tilde{N}(s;M)$, etc. 

Note that once $\beta(s)$ is known, along \eqref{312-N},\eqref{312-a-s-M} \eqref{312-mu-s-M} and \eqref{312-D-s-M} we can subsequently get $\tilde{N},\tilde{a},\tilde{\mu}$ and finally $D(s)$. Therefore we only need to solve for $\beta$. Substituting \eqref{312-a-s-M} into \eqref{312-beta}, we get

\begin{equation}\label{312-beta-2}
    \beta(s)=\exp\left(\int_0^{s}\frac{\tilde{N}(s)}{\beta^2(\bar{s})(a_1+(a_0-ca_1)\tilde{N})}d\bar{s}\right)
\end{equation}

 We have shown that $\beta$ satisfies the ODE \eqref{ode-beta-312}. But here it is more convenient to introduce $\gamma(s):=\beta^2(s)$ and to work with $\gamma(s)$. From \eqref{312-N} we write $\tilde{N}$ as a function of $\gamma$ and $M$, denoted as $\tilde{N}(\gamma,M)$ with abuse of notation
\begin{equation}\label{312-Ngamma}
    \tilde{N}(\gamma,M):=\frac{(1-a_1\gamma M)_+}{a_0\gamma M+c(1-a_1\gamma M)_+}.
\end{equation}
Taking the square for both sides of \eqref{312-beta-2}, we get an equation for $\gamma$
\begin{equation}\label{deq-gamma}
    \gamma(s)=\exp\left(\int_0^{s}\frac{2\tilde{N}(\gamma(\bar{s}),M(\bar{s}))}{\gamma(\bar{s})(a_1+(a_0-ca_1)\tilde{N})}d\bar{s}\right)
\end{equation}
From \eqref{312-beta} we know a priori that $\beta(s)\geq 0$ thus $\tilde{N}\in [0,1/c]$ by \eqref{312-N}. Hence, in turn by \eqref{312-beta} we know $\beta(s)\geq 1$ and is $C^1$. Therefore $\gamma=\beta^2 \geq 1$ and is in $C^1$. Then we differentiate \eqref{deq-gamma} with respect to $s$ to derive an ODE
\begin{align}\notag
    \gamma'(s)&=\frac{2\tilde{N}(\gamma({s}),M({s}))}{\gamma({s})(a_1+(a_0-ca_1)\tilde{N})}\exp\left(\int_0^{s}\frac{2\tilde{N}(\gamma(\bar{s}),M(\bar{s}))}{\gamma(\bar{s})(a_1+(a_0-ca_1)\tilde{N})}d\bar{s}\right)\\ &=\frac{2\tilde{N}(\gamma(s),M(s))}{a_1+(a_0-ca_1)\tilde{N}}=:F(\gamma(s),M(s))\label{ode-gamma}
\end{align}
The rest of the proof is just checking the well-posedness of this ODE \eqref{ode-gamma}. 

Thanks to that $\tilde{N}\in[0,1/c]$, we obtain
\begin{equation}\label{bound-F}
    0\leq F(\gamma,M)\leq \frac{2}{\min(a_0,a_1c)},\quad \forall \gamma\geq0, M\geq 0.
\end{equation}
Hence, with the initial condition $\gamma(0)=1$ by setting $s=0$ in \eqref{deq-gamma}, the solution of \eqref{ode-gamma} satisfies
\begin{equation}\label{bound-gamma}
    1\leq \gamma(s)\leq Cs+1.
\end{equation}
Therefore to show the well-posedness of the ODE \eqref{ode-gamma}, it suffices to show local well-posedness by checking that $F(\gamma,M)$ is Lipschitz in $\gamma$. In fact, we shall show
\begin{equation}\label{pF-pg}
    \left|\frac{\p F}{\p \gamma}(\gamma,M)\right|\leq C,\quad \forall \gamma\geq 1, M\geq 0.
\end{equation}
By a direct calculation from \eqref{ode-gamma} we have
\begin{align}\label{312-tmp-pf}
    \left|\frac{\p F}{\p \gamma}\right|&=\left|\frac{2a_1}{(a_1+(a_0-ca_1)\tilde{N})^2}\frac{\p\tilde{N}}{\p{\gamma}}\right|\leq C\left|\frac{\p\tilde{N}}{\p{\gamma}}\right|,
\end{align} where we have used $a_1+(a_0-ca_1)\tilde{N}\geq \min(a_0/c,a_1)$. Then from \eqref{312-Ngamma} we calculates
\begin{align}\label{pN-pg}
    \frac{\p\tilde{N}}{\p{\gamma}}(\gamma,M)=\mathbb{I}_{\{1-a_1\gamma M\geq 0\}}\frac{(-a_0M)}{(a_0\gamma M+c(1-a_1\gamma M)_+)^2},\quad\quad  \gamma,M\geq 0.
\end{align}
Note that when $0\leq a_1\gamma M\leq 1$ and $\gamma>0$, we have $M\leq C/\gamma$ and 
\begin{align*}
    a_0\gamma M+c(1-a_1\gamma M)_+\geq \min(1,\frac{a_0}{a_1})>0.
\end{align*}
Hence, from \eqref{pN-pg} we derive
\begin{equation}\label{pN-pg-est}
    \left|\frac{\p\tilde{N}}{\p{\gamma}}(\gamma,M)\right|\leq C/\gamma\leq C,\quad \forall \gamma\geq 1,M\geq 0,
\end{equation} combining which with \eqref{312-tmp-pf} we deduce \eqref{pF-pg}.

Therefore we can get the unique $\gamma(s)$ through the ODE \eqref{ode-gamma}. Then we obtain unique $\beta(s)=\sqrt{\gamma(s)}$ and subsequently $\tilde{N}(s)$,$\tilde{a}(s)$ $\tilde{\mu}(s)$ and $D(s)$ from \eqref{312-N},\eqref{312-a-s-M},\eqref{312-mu-s-M} and \eqref{312-D-s-M}. 
\end{proof}

Now we give the proof of Lemma \ref{lemma:key-bound}.
\begin{proof}[Proof of Lemma \ref{lemma:key-bound}]
We introduce $\gamma(s;M)=\beta^2(s;M)$ as in the proof of Lemma \ref{lemma:existence}. The bound on $\beta$ \eqref{bound-beta} follows from \eqref{bound-gamma}. For $D(s;M)$ from its definition \eqref{312-D-s-M} we have, thanks to $\tilde{N}\in[0,1/c]$,
\begin{align*}
    |D(s;M)|&=\left|\frac{b+(b_0-cb)\tilde{N}}{\beta(s;M)(a_1+(a_0-ca_1)\tilde{N})}\right|\\&\leq \frac{\max(|b_0|/c,|b|)}{\beta(s;M)\min(a_0/c,a_1)}\leq \frac{\max(|b_0|/c,|b|)}{\min(a_0/c,a_1)}<+\infty,
\end{align*} where in the second inequality we use $\beta(s;M)\geq 1$. Therefore the uniform bound \eqref{bound-D-key} holds.

For the Lipschitz continuity, we first work with the ODE \eqref{ode-gamma} to derive an estimate for $\gamma$. Similar to \eqref{312-tmp-pf} and \eqref{pN-pg} we calculate 
\begin{align}\label{tmp-pfpm}
        \left|\frac{\p F}{\p M}\right|&=\left|\frac{2a_1}{(a_1+(a_0-ca_1)\tilde{N})^2}\frac{\p\tilde{N}}{\p{M}}\right|\leq C\left|\frac{\p\tilde{N}}{\p{M}}\right|,\\
        \frac{\p\tilde{N}}{\p{M}}(\gamma,M)&=\mathbb{I}_{\{1-a_1\gamma M\geq 0\}}\frac{(-a_0\gamma)}{(a_0\gamma M+c(1-a_1\gamma M)_+)^2},\quad\quad \forall \gamma,M\geq 0.\notag
\end{align}
Dealing with the denominator in the same way as \eqref{pN-pg-est}, we obtain
\begin{equation}\label{pN-pM-est}
    \left|\frac{\p\tilde{N}}{\p{M}}(\gamma,M)\right|\leq C\gamma,\quad \forall \gamma\geq 1,M\geq 0.
\end{equation} And therefore from \eqref{tmp-pfpm} we have  
\begin{equation}\label{pF-pM}
    \left|\frac{\p F}{\p M}(\gamma,M)\right|\leq C\gamma,\quad \forall \gamma\geq 1, M\geq 0.
\end{equation}
Now we start the estimate, suppose $M_1\equiv M_2$ on $[0,S_0-\sigma]$ with some $\sigma\in[0,S_0]$. In other words, they only start to differ at $s_0:=S_0-\sigma$. From the ODE \eqref{ode-gamma} we deduce for $s_0\leq s\leq S_0$  
\begin{align*}
    |\gamma_1(s)-\gamma_2(s)|&=|\int_{s_0}^sF(\gamma_1,M_1)dt-\int_{s_0}^sF(\gamma_2,M_2)dt|\\&\leq \int_{s_0}^s|F(\gamma_1,M_1)dt-\int_{s_0}^sF(\gamma_2,M_2)|dt,
\end{align*}where we use the shorthand $\gamma_i(s):=\gamma(s,M_i),i=1,2$. By the triangle inequality and derivative bounds \eqref{pF-pg},\eqref{pF-pM} we get
\begin{align*}
    |\gamma_1(s)-\gamma_2(s)|&\leq C\int_{s_0}^s|\gamma_1(t)-\gamma_2(t)|dt+C\left(\max_{t\in[s_0,s],i=1,2}\gamma_i(t)\right)\int_{s_0}^s|M_1(t)-M_2(t)|dt.
\end{align*}
Then using the bound on $\gamma$ \eqref{bound-gamma} we deduce 
\begin{align*}
    |\gamma_1(s)-\gamma_2(s)|& \leq C\int_{s_0}^s|\gamma_1(t)-\gamma_2(t)|dt+ C(s+1)(s-s_0)\|M_1-M_2\|_{C[s_0,s]}\\& \leq C\int_{s_0}^s|\gamma_1(t)-\gamma_2(t)|dt+ C_S(S_0-s_0)\|M_1-M_2\|_{C[s_0,S_0]},
\end{align*} where we use $C_S$ for a constant depending on $S$, an upper bound for $S_0$. Applying Gronwall's inequality we conclude
\begin{align}\notag
    \|\gamma_1-\gamma_2\|_{C[s_0,S_0]}&\leq C_S(S_0-s_0)(1+(S_0-s_0)e^{C(S_0-s_0)})\|M_1-M_2\|_{C[s_0,S_0]}\\
    &\leq C_S(S_0-s_0)\|M_1-M_2\|_{C[s_0,S_0]}\label{estimate-gamma}
\end{align}Recall $s_0=S_0-\sigma$, and the estimate for $\beta$ \eqref{beta-pertub} follows from \eqref{estimate-gamma} since $\beta_i=\sqrt{\gamma_i}$ and $\beta_i\geq 1$, where $\beta_i$ is the shorthand for $\beta(\cdot;M_i)$, for $i=1,2$.

It remains to derive the estimate for $D$ \eqref{D-pertub}. By the triangle inequality we get
\begin{align*}
|D(s;M_1)-D(s;M_2)|&\leq \left|\frac{\tilde{N}_1b_0+(1-c\tilde{N}_1)b}{\beta_1(a_1+(a_0-a_1c)\tilde{N}_1)}-\frac{\tilde{N}_2b_0+(1-c\tilde{N}_2)b}{\beta_2(a_1+(a_0-a_1c)\tilde{N}_2)}\right|\\
&\leq \frac{1}{\beta_1}\left|\frac{\tilde{N}_1b_0+(1-c\tilde{N}_1)b}{(a_1+(a_0-a_1c)\tilde{N}_1)}-\frac{\tilde{N}_2b_0+(1-c\tilde{N}_2)b}{(a_1+(a_0-a_1c)\tilde{N}_2)})\right|\\
&\quad \quad \quad \quad +\left|(\frac{1}{\beta_1}-\frac{1}{\beta_2})\frac{\tilde{N}_2b_0+(1-c\tilde{N}_2)b}{(a_1+(a_0-a_1c)\tilde{N}_2)}\right|,
\end{align*} where $\tilde{N}_i$ is the shorthand for $\tilde{N}(\gamma_i,M_i)$ for $i=1,2$. Using the mean value theorem to estimate the above two terms, and noting that $\beta\geq 1,\tilde{N}\in[0,1/c]$, we obtain
\begin{align}\label{312-pf-2}
    |D(s;M_1)-D(s;M_2)|&\leq C|\tilde{N}(\gamma_1,M_1)-\tilde{N}(\gamma_2,M_2)|+C|\beta_1(s)-\beta_2(s)|.
\end{align}With bounds on derivatives of $\tilde{N}$ \eqref{pN-pg-est} and \eqref{pN-pM-est}, we derive
\begin{align*}
        |\tilde{N}(\gamma_1,M_1)-\tilde{N}(\gamma_2,M_2)|&\leq C|\gamma_1(s)-\gamma_2(s)|+C\left(\max_{i=1,2}\gamma_i(s)\right)|M_1(s)-M_2(s)|\\&\leq C|\gamma_1(s)-\gamma_2(s)|+C_S|M_1(s)-M_2(s)|,
\end{align*} where in the last inequality we use again \eqref{bound-gamma}, to bound $\gamma(s)$ with a constant depending on the upper bound of time $S$. Plugging the above estimate into \eqref{312-pf-2} we deduce
\begin{align*}
    |D(s;M_1)-D(s;M_2)|&\leq C_S(|\gamma_1(s)-\gamma_2(s)|+|\beta_1(s)-\beta_2(s)|+|M_1(s)-M_2(s)|)
\end{align*} Finally, applying the estimate for $\gamma$ \eqref{estimate-gamma} with $s_0=0$ and for $\beta$ \eqref{beta-pertub} with $\sigma=S_0$, we derive the bound for $D$ \eqref{D-pertub}.
\end{proof}

\subsection{Proof of Proposition \ref{cpde-4.3}}

Here we present the proof of Proposition \ref{cpde-4.3}. 
\begin{proof}[Proof of Proposition \ref{cpde-4.3}]
The proof is an adaption of  \cite[Proof of Proposition 4.3]{carrillo2013CPDEclassical}.

For some $\eps>0$ small to be determined below, we consider the integral formula \eqref{M-3.6-t0} for $M$ with starting time $s_0-\eps$
\begin{equation}
        \begin{aligned}
        M(s)&=-2\int_{-\infty}^{\ell(s_0-\eps)}G(\ell(s),s,\xi,0)u_x(\xi,s_0-\eps)d\xi+2\int_{s_0-\eps}^sM(\tau)G_x(\ell(s),s,\ell(\tau),\tau)d\tau\\\quad\quad &-2\int_{s_0-\eps}^sM(\tau)G_x(\ell(s),s,\ell_R(\tau),\tau)d\tau.\\
        &=:K_1+K_2+K_3,
    \end{aligned}
\end{equation} where $s\in(s_0-\eps,s_0)$ and $G$ is the heat kernel given in \eqref{heat-kernel}.

The estimates for $K_2$ and $K_3$ are different from that in \cite[Proof of Proposition 4.3]{carrillo2013CPDEclassical}. We shall use the global Lipschitz continuity of $\ell(s)$ by Lemma \ref{lemma:Lip-s}, while in the setting of \cite[Proposition 4.3]{carrillo2013CPDEclassical} some sign condition has been utilized to give controls.

We define
\begin{equation*}
    \Phi(s):=\sup_{s_0-\eps<\tau<s_0}M(\tau).
\end{equation*}
Then we get
\begin{equation}\label{pf43-tmp1}
    \begin{aligned}
    |K_2|&\leq 2\Phi(s)\int_{s_0-\eps}^s|G_x(\ell(s),s,\ell(\tau),\tau)|d\tau,\\
    |K_3|&\leq 2\Phi(s)\int_{s_0-\eps}^s|G_x(\ell(s),s,\ell_R(\tau),\tau)|d\tau.
\end{aligned}
\end{equation}We want to gain smallness of the two integrals in \eqref{pf43-tmp1} when $\eps$ is small enough. For $K_2$, we note that $|G_x(\ell(s),s,\ell(\tau),\tau)|\leq C\frac{1}{(s-\tau)^{1/2}}$ thanks to that $\ell(s)$ is Lipschitz continuous with a uniform constant by Lemma \ref{lemma:Lip-s}. Therefore one deduces
\begin{equation}\label{pf43-tmp2}
    2\int_{s_0-\eps}^s|G_x(\ell(s),s,\ell(\tau),\tau)|d\tau\leq C\int_{s_0-\eps}^s\frac{1}{(s-\tau)^{1/2}}d\tau\leq C\eps^{1/2}.
\end{equation} Here and in the following $C$ denotes a universal constant which does not depend on the solution or $s_0$.

For $K_3$, by the Lipschitz continuity from Lemma \ref{lemma:Lip-s} again we deduce
\begin{align*}
    \ell(s)-\ell_R(\tau)&= -V_R\beta(\tau)+\ell(s)-\ell(\tau)\geq -V_R\beta(\tau)-L(s-\tau)
\end{align*}
Using that $V_R<0$ (since we have set $V_F=0$) and $\beta(\tau)\geq 1$ by Lemma \ref{lemma:key-bound}, we obtain for $s_0-\eps\leq \tau\leq s\leq s_0$
\begin{align*}
    \ell(s)-\ell_R(\tau)\geq  |V_R|-L\eps>0,
\end{align*}if we choose $\eps$ small enough. Here $L$ is the Lipschitz constant of $\ell(s)$. Therefore, with the same calculation as in \cite[Equation (3.13)]{carrillo2013CPDEclassical}, we get
\begin{equation}\label{pf43-tmp3}
2\int_{s_0-\eps}^s|G_x(\ell(s),s,\ell_R(\tau),\tau)|d\tau\leq C\int_{{\frac{|V_R|-L\eps}{\sqrt{8\eps}}}}^{+\infty}\frac{1}{z}e^{-z^2}dz,
\end{equation} which can be made small by choosing $\eps$ small enough. 

In view of \eqref{pf43-tmp1},\eqref{pf43-tmp2} and \eqref{pf43-tmp3}, we can choose $\eps$ small enough, such that 
\begin{equation*}
    M(s)\leq \frac{1}{2}\Phi(s)+K_1,\quad \forall t\in(s_0-\eps,s_0).
\end{equation*} For $K_1$, by the integrability of the heat kernel one directly has $$|K_1|\leq 2\sup_{x\in(-\infty,\ell(s_0-\eps)]}|u_x(x,s_0-\eps)|=2\mathscr{M}_0.$$
Therefore for $s\in(s_0-\eps,s_0)$ it holds that $ M(s)\leq \frac{1}{2}\Phi(s)+2\mathscr{M}_0$.
Since $\Phi$ is an increasing function, we get that
\begin{equation}
            M(s)\leq \frac{1}{2}\Phi(t)+2\mathscr{M}_0,\quad s_0-\eps < s\leq t<s_0.
\end{equation} Taking the supremum for $s$ in $(s_0-\eps,t)$, we derive
\begin{equation}
    \Phi(t)\leq \frac{1}{2}\Phi(t)+2\mathscr{M}_0,\quad \forall t\in(s_0-\eps,s_0).
\end{equation} Note that $\Phi(t)<\infty$, since we are working with a classical solution on $[0,s_0)$. Therefore we conclude that $\Phi(t)\leq 4\mathscr{M}_0$ for all $t\in(s_0-\eps,s_0)$, which gives
\begin{equation*}
    \sup_{s_0-\eps<t<s_0}M(t)\leq 4\mathscr{M}_0<+\infty.
\end{equation*} Finally, since the solution is classical on $[0,s_0)$ we have $\sup_{0\leq t\leq  s_0-\eps}M(t)<+\infty$. Therefore
\begin{equation*}
        \sup_{0\leq t<s_0}M(t)<+\infty.
\end{equation*} The proof is complete.
\end{proof}

\section{Proof of a Poincar\'e inequality}\label{app:poincare}

In this appendix, we prove the following Poincar\'e inequality, which is similar to \cite[Appendix]{caceres2011analysis}.
\begin{lemma}\label{lemma-poincare}
Let $n(x)=K\min(x,e^{-x})$, where $K>0$ is the normalization constant such that $\int_0^{+\infty}n(x)dx=1$. Then for some $\alpha>0$, we have
\begin{equation}
    \alpha\int_0^{+\infty}|f(x)|^2n(x)dx\leq \int_0^{+\infty}|f'(x)|^2n(x)dx,
\end{equation} for any function $f\in H^1((0,+\infty);n(x)dx)$ such that
\begin{equation}
    \int_0^{+\infty}f(x)n(x)dx=0.
\end{equation}
\end{lemma}
Through a change of variable from $(-\infty,V_F)$ to $(0,+\infty)$, Lemma \ref{lemma-poincare} implies the Poincare inequality for $p_{\infty}$ in Proposition \ref{prop:poincare}, since $p_{\infty}$ has equivalent asymptotic behaviors at $x=0$ ($v=V_F$) and $x=+\infty$ ($v=-\infty$) with $n(x)$.

\begin{proof}[Proof of Lemma \ref{lemma-poincare}]
By the Muckenhoupt Criterion for Poincar\'e inequalities (see for example \cite[Theorem 4.5.1]{bakry2014analysis}), it suffices to check that the following two quantities 
\begin{equation}
\begin{aligned}
        B_+&=\sup_{x>1}\left(\int_x^{+\infty}n(y)dy\right)\left(\int_1^x\frac{1}{n(z)}dz\right),\\ B_-&=\sup_{0< x<1}\left(\int_0^{x}n(y)dy\right)\left(\int_x^1\frac{1}{n(z)}dz\right),
\end{aligned}
\end{equation}are finite. 

Indeed, noting the behavior of $n(x)$ at $x=+\infty$, we have
\begin{align*}
    B_+&\leq C\sup_{x>1}\left(\int_x^{+\infty}e^{-y}dy\right)\left(\int_1^xe^{z}dz\right)\\&=C\sup_{x>1}(e^{-x}(e^x-1))<\infty,
\end{align*}
And using the behavior of $n(x)$ at $x=0$, we obtain
\begin{align*}
    B_-&\leq C\sup_{0<x<1}\left(\int_0^{x}ydy\right)\left(\int_x^1\frac{1}{z}dz\right)\\&\leq C\sup_{0<x<1}(x^2(-\ln x))<\infty.
\end{align*}
Therefore the proof is complete.
\end{proof}

\end{appendices}

\bibliography{z_neuron}
\bibliographystyle{plain}
\end{document}